\documentclass{article}
\usepackage[letterpaper,left=1.5in, right= 1.5in, top=1in, bottom=1in]{geometry}
\newenvironment{cpf}{\begin{trivlist} \item[] {\em Proof of Claim.}}{\hspace*{\stretch{1}} $\diamond$ \end{trivlist}}

\usepackage{lmodern}
\usepackage{blkarray}
\usepackage{xcolor}
\usepackage{amssymb, xfrac}
\usepackage{bm}
\usepackage{nicematrix}

\usepackage{amsfonts,mathtools,nicefrac,amsmath} % mathematical writing
\usepackage{amsthm}
\usepackage{thmtools}
\usepackage{url}
\usepackage{thm-restate}
\usepackage[hidelinks]{hyperref}
\usepackage[capitalize]{cleveref}
\usepackage{todonotes}
\usepackage{caption,subcaption}
\usepackage[shortlabels]{enumitem}
\usepackage{tikz}
\usetikzlibrary{decorations.pathmorphing}
\usetikzlibrary{shapes.geometric}
\usetikzlibrary{arrows.meta}
\usetikzlibrary{arrows}
\usepackage{cite} % it allows more compact citations (e.g. from [1,2,3,4,5] to [1-5]) 

\DeclareMathOperator{\ocp}{ocp}
\DeclareMathOperator{\supp}{supp}
\DeclareMathOperator{\internal}{in}

\newcommand{\mbb}[1]{\mathbb{#1}}

\newcommand{\mbfs}[1]{\boldsymbol{#1}}
\newcommand{\pin}{\subseteq_P}
\newcommand{\qin}{\subseteq_Q}
\newcommand{\children}{\mathrm{children}}
\newcommand{\descendants}{\mathrm{descendants}}
\newcommand{\leaves}{\mathrm{leaves}}
\newcommand{\dis}{\mathrm{dis}}
\newcommand{\cupdot}{\mathbin{\dot{\cup}}}

\theoremstyle{theorem}
\newtheorem{theorem}{Theorem}[section]
\newtheorem{lemma}[theorem]{Lemma}
\newtheorem{corollary}[theorem]{Corollary}
\newtheorem{proposition}[theorem]{Proposition}
\newtheorem{conjecture}[theorem]{Conjecture}
\newtheorem{claim}[theorem]{Claim}
\newtheorem{observation}[theorem]{Observation}
\theoremstyle{definition}
\newtheorem{definition}[theorem]{Definition}

\title{Characterizing unimodular laminar hypergraphs via forbidden subhypergraphs
\thanks{M. Caoduro and J. Paat were supported by an NSERC Grant [RGPIN-2021-02475].
Part of the research for this project was done while M. Neuwohner was affiliated with the London School of Economics and Political Science. During that time, M. Neuwohner was supported by the Engineering and Physical Sciences Research Council, which is part of UK Research and Innovation, grant ref. EP/X030989/1.}}
\author{Marco Caoduro\thanks{Laboratoire d'Informatique de Paris Nord, Université Sorbonne Paris Nord, Paris, France (\href{mailto:caoduro@lipn.univ-paris13.fr}{caoduro@lipn.univ-paris13.fr})}
\and Meike Neuwohner\thanks{CNRS \& DIENS, ENS Paris - PSL Research University, Paris, France. (\href{mailto:meike.neuwohner@ens.fr}{meike.neuwohner@ens.fr})}
\and Joseph Paat\thanks{Sauder School of Business, University of British Columbia, Vancouver, Canada  (\href{mailtojoseph.paat@sauder.ubc.ca:}{joseph.paat@sauder.ubc.ca})} }

\date{}
\begin{document}	

\maketitle

%%%%%%%%%%%%%%%%%%%%%%%%%%%%%%
%%%%%%%%%%%%%%%%%%%%%%%%%%%%%%
\begin{abstract}
The incidence matrix of a graph is totally unimodular if and only if the graph is bipartite, i.e., it contains no odd cycles.
We extend the characterization of total unimodularity to hypergraphs whose hyperedges of size at least four form a laminar family. 
Such hypergraphs have been used to model problems with fairness constraints that ensure balanced representation, among other applications.
Our main result shows that total unimodularity for laminar hypergraphs is equivalent to forbidding odd cycles and structures that we call avocados and tree houses.
As a corollary, we resolve a special case of a conjecture on almost totally unimodular matrices, originally posed by Padberg and later modified by Cornu\'ejols and Zuluaga.
We discuss applications of laminar hypergraphs and connect our results to integer programming with bounded subdeterminants.
\end{abstract}

%%%%%%%%%%%%%%%%%%%%%%%%%%%%%%
%%%%%%%%%%%%%%%%%%%%%%%%%%%%%%
\section{Introduction}\label{sec:intro}
%%%%%%%%%%%%%%%%%%%%%%%%%%%%%%
%%%%%%%%%%%%%%%%%%%%%%%%%%%%%%

A matrix $\mbfs{A}\in \mbb{Z}^{m\times n}$ is {\it totally unimodular (TU)} if every square submatrix has an absolute determinant of at most $1$.
For decades, TU matrices have found applications in integer optimization, including the bipartite matching problem.
Indeed, it is well-known that a graph has a TU incidence matrix if and only if the graph contains no odd cycles, i.e., it is bipartite; see, e.g.,~\cite[Theorem 18.2]{S2003V1}.
The motivation for our work is to extend this well-known result to hypergraphs, that is, we aim to identify forbidden subhypergraphs that characterize when a hypergraph has a TU incidence matrix.

A {\it hypergraph} $G=(V,E)$ is a pair of a finite set $V =: V(G)$ of vertices and a finite multiset $E =: E(G) \subseteq 2^V$ of non-empty hyperedges.
We consider graphs as a special type of hypergraph and, as such, use the same notation for both.
The {\it incidence matrix} $\mbfs{M}(G)\in \{0,1\}^{V\times E}$ of $G$ is defined as $\mbfs{M}(G)_{v,e}=1$ if and only if $v\in e$.
With subdeterminants in mind, our notion of `subhypergraph' should correspond to taking a submatrix of an incidence matrix of a larger hypergraph.
For this purpose, we use the notion of `partial subhypergraph' from Schrijver~ \cite[\S77.1]{S2003V3}.
A {\it partial subhypergraph} is a hypergraph of the form $G[U,F]\coloneqq (U,F[U])$, where $U \subseteq V$, $F \subseteq E$, and
\[
F[U]\coloneqq \{f\cap U: f\in F~\text{and}~f \cap U \neq \emptyset\}.
\]
We say that a hypergraph $G$ \emph{contains} a hypergraph $H$ if $G$ has a partial subhypergraph isomorphic to $H$.\footnote{For technical reasons, we will later work with a slightly broader notion of partial subhypergraphs.}

Hypergraphs with TU incidence matrices are called {\it unimodular}.
Our goal is to identify a family of forbidden hypergraphs such that the following holds: a hypergraph $G$ is unimodular if and only if it contains no hypergraph $H$ in the family as a partial subhypergraph.
For graphs, this family consists precisely of the odd cycles.
Looking to general hypergraphs, a characterization of unimodular hypergraphs based on partial subhypergraphs would provide a complete list of almost TU (i.e., minimally non-TU) matrices with $\{0,1\}$-entries. 
Finding a complete list of almost TU matrices seems unlikely, as indicated by Truemper's~\cite{TRUEMPER1992302} involved recursive characterization of almost TU matrices.
Thus, we focus on structured families of hypergraphs.

For hypergraphs with no hyperedge larger than size $3$, Berge shows that forbidding odd cycles also suffices for unimodularity.
\begin{theorem}[Berge {\cite[\S 5.3]{Berge_1989}}]\label{prop:unimodular_for_rank_3}
Let $G$ be a hypergraph such that each hyperedge has size at most $3$. 
Then $G$ is unimodular if and only if $G$ contains no odd cycles.
\end{theorem}
Forbidding odd cycles as partial subhypergraphs is no longer sufficient for unimodularity when there are hyperedges with more than $3$ vertices.
Berge~\cite[\S 5.3]{Berge_1989} gives a counterexample that we reproduce in~\cref{figure:OTH}.
Note that the odd cycle $r, \{r, \ell_1\}, \ell_1, \{\ell_1, \ell_2\}, \ell_2, \{r, \ell_2\}, r$ is not a partial subhypergraph because $r$ appears as a vertex in the cycle, yet it has been removed from the hyperedge $\{r, \ell_1, \ell_2, \ell_3\}$ to create the edge $\{\ell_1, \ell_2\}$.
We use the term tree house to refer to hypergraphs that arise from the one shown in \cref{figure:OTH} by replacing the edges $\{r,\ell_i\}$ by odd paths.

\begin{definition}\label{def:strong_odd_tree}
	A {\bf tree house} is a hypergraph $T=(V,E\cup\{h\})$ such that
	\begin{enumerate}[{\it({T}1)}, leftmargin = *]
		\item \label{item:oth_1} $h=:\{r,\ell_1,\ell_2,\ell_3\}$ has size $4$.
		\item \label{item:oth_2} $(V,E)$ is a tree consisting of three odd paths $P_1,P_2,P_3$.
		\item \label{item:oth_3} The endpoints of $P_i$ are $r$ and $\ell_i$, and the vertex sets $V(P_i)\setminus\{r\}$, $i\in\{1,2,3\}$, are pairwise disjoint.
	\end{enumerate}
\end{definition}

\begin{figure}[h]
\footnotesize
	\centering
	\begin{tabular}{c@{\hskip 1.5 cm}c}
	\begin{tikzpicture}[scale = 1, mynode/.style={circle, draw, fill, minimum size = .15 cm, inner sep=0pt}, baseline = 10]
		\node[mynode, draw = black, fill = black,label={[left = 3 pt] \color{black}$r$}] (A) at (0,0){};
		\node[mynode, draw = black, fill = black,label={[right= 3 pt] \color{black}$\ell_1$}] (B) at (1,0){};
		\node[mynode, draw = black, fill = black,label={[right= 3pt] \color{black}$\ell_2$}] (C) at (1,1){};
		\node[mynode, draw = black, fill = black,label={[left= 3pt] \color{black}$\ell_3$}] (D) at (0,1){};
		\draw[] (A)--(B);
		\draw[] (A)--(C);
		\draw[] (A)--(D);
		\draw[black, rounded corners] (-0.65,-0.25) rectangle (1.65,1.35);
	\end{tikzpicture}
	&
	$
	\begin{bNiceArray}{cccc}[first-row, first-col]
	&\{r, \ell_1\}&\{r, \ell_2\}&\{r, \ell_3\}&\{r, \ell_1, \ell_2, \ell_3\}\\
	r & 1 & 1 & 1 & 1\\
	\ell_1 & 1 & 0 & 0 & 1\\
	\ell_2 & 0 & 1 & 0 & 1\\
	\ell_3 & 0 & 0 & 1 & 1 
	\end{bNiceArray}
	$
	\end{tabular}
	\caption{A tree house and its incidence matrix.	}\label{figure:OTH}
\end{figure}

Our first result shows that forbidding odd cycles and tree houses characterizes unimodular \emph{disjoint} hypergraphs.

\begin{definition}
A hypergraph is {\bf disjoint} if its hyperedges of size at least $4$ are pairwise disjoint.\label{def:disjoint_hypergraph}
\end{definition}
While we could define disjoint hypergraphs based on the hyperedges of size at least $3$, our results also hold for the weaker notion in \cref{def:disjoint_hypergraph}.

\begin{restatable}{theorem}{TheoremSizeFour}\label{theorem:size_four}
Let $G$ be a disjoint hypergraph.
Then $G$ is unimodular if and only if it neither contains an odd cycle nor a tree house.
\end{restatable}
We prove \cref{theorem:size_four} in \cref{secUnimodular} through \cref{subsec:no_odd_cycle}.
After the initial submission of this work, Grappe and Nicolas-Thouvenin~\cite{GRAPPE2026107426} gave a simple proof characterizing total unimodularity for incidence matrices of graphs, plus one additional row or column; their work recovers a special case of Berge's work and of our work, respectively.
\cref{theorem:size_four} generalizes naturally to the setting of disjoint \emph{mixed} hypergraphs, where we allow for $\{0,\pm 1\}$-entries in the incidence matrix.
We state our second main result, \cref{thm:DirectedUniEquivalence}, but defer the details, including definitions, to \cref{secSigned}.

\begin{theorem}\label{thm:DirectedUniEquivalence}
	For a disjoint mixed hypergraph $D$, $\mbfs{M}(D)$ is unimodular if and only if there is no mixed odd cycle or mixed tree house in $D$.
\end{theorem}

We remark that a generalization of \cref{thmOCP} to mixed graphs has been shown by Zaslavsky~\cite{ZASLAVSKY198247}.

After studying disjoint hypergraphs, an extension would involve overlapping hyperedges of size at least $4$.
Laminar hypergraphs form a natural first extension because the overlaps occur in a nested way.
\begin{definition}\label{def:Laminar}
	A {\bf laminar hypergraph} is a hypergraph $G=(V,E)$, where the set of hyperedges of size at least $4$ forms a laminar family.
\end{definition}

If the multiset of all hyperedges (of any size) forms a laminar family, then the incidence matrix is TU; see, e.g.~\cite{CCZ2014}. 
In contrast, laminar hypergraphs according to \cref{def:Laminar}, which allow for overlap in size-2 and size-3 hyperedges, are not always TU.
Laminar (and disjoint) hypergraphs have applications in scheduling and group fairness, and can be used to model problems such as the partial independent transversal problem. 
We discuss this in more detail in \cref{subsec:laminar}.

After applying elementary column operations, the incidence matrix of a laminar hypergraph can be transformed into the incidence matrix of a disjoint hypergraph.
Nevertheless, the family of laminar hypergraphs introduces one new forbidden substructure preventing unimodularity that we call an avocado.
We use $\cupdot$ to denote the disjoint union.
\begin{definition}\label{def:avocado}
	An {\bf avocado} is a hypergraph $A=(V\cup \{t\},E\cup \{s,f\})$ such that
	\begin{enumerate}[{\it({A}1)}, leftmargin = *]
		\item \label{item:avocado_1} We have $t\in s\subseteq f$, $|s|\ge 4$, and $t \not \in e$ for each $e \in E$.
		\item \label{item:avocado_2} The hypergraph $C\coloneqq (V,E\cupdot \{f\setminus s\})$ is an odd cycle.
		\item \label{item:avocado_3} For every pair of vertices $v, w \in f \cap V$, if the $v$-$w$-subpath $C[v,w]$ of $(V,E)$ only intersects $f$ at $v$ and $w$, then $C[v,w]$ is odd.
			\end{enumerate}
\end{definition}

\begin{figure}[h]
	\centering
	
	\begin{tikzpicture}[scale = .75, mynode/.style={circle, draw, fill, minimum size = .15 cm, inner sep=0pt}, baseline = 10]
		\node[mynode, draw = black, fill = black,label={[above = 2 pt] \color{black}$t$}] (A) at (0,0){};
		\node[mynode, draw = black, fill = black]  (B) at (2.35,.8){};
		\node[mynode, draw = black, fill = black] (C0) at (.75,.8){};
		\node[mynode, draw = black, fill = black] (C1) at (.75,-.2){};
		\node[mynode, draw = black, fill = black] (C2) at (.75,.3){};
		\node[mynode, draw = black, fill = black]  (F) at (2.35,-.2){};
		
		\draw[] (B)--(C0)--(C2) -- (C1) --(F);
		
		\draw[draw = brown!60!black, very thick, rounded corners] (-.5,-.5) rectangle (1.65,1);
		\draw[draw = green!50!black, very thick, rounded corners] (-.75,-.65) rectangle (2.75,1.25);
	\end{tikzpicture}	
	\caption{An avocado, where the `seed' $s$ is drawn in brown and the `flesh' $f$ is drawn in green. The vertex $t$ is the `tip of the avocado'.}\label{figure:AVO}
\end{figure}

\begin{restatable}{theorem}{TheoremNonUnimodularLaminar}\label{theorem:non_unimodular_laminar}
	Let $G$ be a laminar hypergraph. 
	Then $G$ is unimodular if and only if it does not contain an odd cycle, a tree house, or an avocado.
\end{restatable}
We prove \cref{theorem:non_unimodular_laminar} in \cref{sec:proof_main_thm}.

Theorems~\ref{theorem:size_four} and~\ref{theorem:non_unimodular_laminar} characterize a family of unimodular laminar hypergraphs using forbidden partial subhypergraphs.
As a consequence, we provide a complete list of almost TU incidence matrices for these families.
Our proof introduces the notion of quasi-subhypergraph, which, loosely speaking, generalizes partial subhypergraphs by allowing for more types of containment relationships.
We believe that quasi-subhypergraphs may help characterize other families of unimodular hypergraphs using subhypergraphs.
The remainder of the introduction serves as motivation of our work: \cref{subsec:laminar} discusses laminar hypergraphs, \cref{subsec:subdet} connects to integer programs with bounded subdeterminants, and \cref{subsec:CZ} discusses how our work partially resolves a conjecture by Cornu{\'e}jols and Zuluaga.
We provide preliminaries in~\cref{sec:Prelims} and introduce quasi-subhypergraphs in \cref{ssecQuasi}.
Following this, we prove our main results.
Finally, future directions are given in~\cref{secFuture}.
%

%%%%%%%%%%%%%%%%%%%%%%%%%%%%%%
%%%%%%%%%%%%%%%%%%%%%%%%%%%%%%
\subsection{Disjoint and laminar hypergraphs}\label{subsec:laminar}
%%%%%%%%%%%%%%%%%%%%%%%%%%%%%%
%%%%%%%%%%%%%%%%%%%%%%%%%%%%%%

%
A variety of problems can be modeled using disjoint or laminar hypergraphs.
Since this paper focuses on incidence matrices, we highlight examples that use the incidence matrix of a disjoint or laminar hypergraph.

Disjoint hypergraphs find applications in problems that consider vertex-colored graphs; here, the disjoint hyperedges are the color classes.
Examples can be found in job interval selection~\cite[\S 4.1]{2015_vanBevern}, (partial) independent transversals~\cite{H2011,HMZ2025}, train scheduling~\cite{ZKRS1996}, and group fairness~\cite{2017_Aharoni,2024_Hojny}.
Loosely speaking, these applications are `colorful' stable set problems, where each color class has lower and upper bounds on the number of representatives in the stable set. 
Moreover, recall that our definitions of disjoint and laminar hypergraphs place no restrictions on the hyperedges of size 3.
Thus, there are also applications in vertex-colored stable set problems on hypergraphs where each hyperedge has size 3; see, e.g.,~\cite[Conjecture 12.3]{AB2023}.

Incidence matrices of disjoint hypergraphs also appear in a classic integer programming problem: the uncapacitated facility location problem, specifically the disaggregated model; see, e.g., (2.19) in~\cite{CCZ2014}.
While this constraint matrix contains $\{0, -1\}$-valued columns, multiplying these by $-1$ yields a $\{0,1\}$-valued matrix, preserving whether or not the matrix is TU.

As laminar hypergraphs are generalizations of disjoint hypergraphs, applications of the corresponding incidence matrices generalize applications from the previous paragraphs.
The authors in~\cite{WMLF2023} consider the fair stable set problem in the context of group fairness; this problem imposes a global capacity constraint on all variables, along with quotas on disjoint subsets of variables.
The corresponding IP formulation features the incidence matrix of a laminar, but non-disjoint hypergraph.

%%%%%%%%%%%%%%%%%%%%%%%%%%%%%%
%%%%%%%%%%%%%%%%%%%%%%%%%%%%%%
\subsection{Connections to integer programs with bounded subdeterminants}\label{subsec:subdet}
%%%%%%%%%%%%%%%%%%%%%%%%%%%%%%
%%%%%%%%%%%%%%%%%%%%%%%%%%%%%%
TU matrices play an important role in integer programming for the following reason: An integer program with a TU constraint matrix can be solved in polynomial time using the linear programming relaxation; see, e.g.,~\cite[Chapter 5]{S2003V1}. 
Applications of TU matrices include problems in network flow, scheduling, and bipartite matching.

Recently, there has been a significant amount of research on bounded subdeterminants beyond the TU case~\cite{AWZ2017,AS2023,BFMR2014,CFHW2020,FFKY2023,FJWY2022,GSW2021,GSMP2022,JB2022,LPSX2022,NCSZ2023,NSZ2022,PSW2022}.
A matrix $\mbfs{A}\in\mathbb{Z}^{m\times n}$ is called \emph{totally $\Delta$-modular} if $\Delta(\mbfs{A})\le \Delta$, where
\[
\Delta(\mbfs{A}) \coloneqq \max \left\{|\det  \mbfs{B} |:\ \mbfs{B}~\text{is a square submatrix of}~\mbfs{A}\right\}.
\]
The value of $\Delta(\mbfs{A})$ is known for problems whose constraint matrix is the incidence matrix of a general graph. 
The following known result generalizes the classic result that graph incidence matrices are TU if and only if the graph is bipartite.

\begin{theorem}[Grossman et al.~\cite{GKS1995}]\label{thmOCP}
For a graph $G = (V,E)$, $\Delta(\mbfs{M}(G)) = 2^{\ocp(G)}$, where $\ocp(G)$ is the odd cycle packing number of $G$, i.e., the maximum number of vertex-disjoint odd cycles in $G$.
\end{theorem}
Now, it is conjectured that for fixed $\Delta$, integer programs with totally $\Delta$-modular constraint matrices can be solved in polynomial time.
If $\Delta=1$, then the matrix at hand is TU, and the conjecture is true.
While this conjecture remains open in general, affirmative results exist in special cases, namely if $\Delta=2$~\cite{AWZ2017} or if one restricts the structure of the matrix.
A seminal result by Fiorini et al.~\cite{FJWY2022} confirms the conjecture if the matrix has at most two non-zero entries per row or per column, respectively.
The authors use \cref{thmOCP} to reduce the problem to the independent set (or matching) problem in graphs of bounded odd cycle packing number.
Tackling the independent set problem in graphs of bounded odd cycle packing number requires deep structural results from the graph minors series (specifically~\cite{ROBERTSON200343}).
This approach highlights the intricate connection between understanding the structural properties of (certain classes of) totally $\Delta$-modular matrices and being able to solve the corresponding integer programs.

Kober~\cite{kober2025totally} extends~\cite{FJWY2022} by showing that totally $\Delta$-modular IPs with at most two non-zeros in all but a constant number of rows or columns can be solved in polynomial time.
To obtain this result, Kober establishes a class of forbidden substructures called $r$-decorated trees and then uses structural results for graphs avoiding these to solve the independent set problem.

Instead of adding a constant number of arbitrary rows/columns, another natural direction to extend the result in~\cite{FJWY2022} is to allow for an arbitrary number of additional rows/columns with restricted interaction.
This motivates our study on forbidden submatrices.
Indeed, in order to pursue a similar approach as in~\cite{FJWY2022,kober2025totally}, we find it important to understand the structure of $\Delta$-modular hypergraphs with bounded subdeterminants, and leverage it in order to obtain a (combinatorial) algorithm for the independent set and matching problem in these hypergraphs.
Besides \cref{thmOCP}, previous work on forbidden submatrices has resulted in partial lists of forbidden submatrices in the context of bounded determinants~\cite{CKW2024,OW2022,PSWX2024}.
In this paper, we take a step forward by characterizing the structure of TU disjoint and laminar hypergraphs.

%%%%%%%%%%%%%%%%%%%%%%%%%%%%%%
%%%%%%%%%%%%%%%%%%%%%%%%%%%%%%
\subsection{Partially resolving a conjecture by Cornu\'{e}jols and Zuluaga}\label{subsec:CZ}
%%%%%%%%%%%%%%%%%%%%%%%%%%%%%%
%%%%%%%%%%%%%%%%%%%%%%%%%%%%%%
%
Padberg~\cite{P1988} discusses different structural characterizations of TU matrices.
In his work, he writes {\it `To characterize TU matrices in terms of forbidden structures it will be convenient to characterize those matrices which are in a minimal sense non TU.'}
Such minimally non-TU matrices are called almost TU.

\begin{definition}
    A matrix $\mbfs{A}$ is \textbf{almost TU} if it is not TU, but every proper submatrix of $\mbfs{A}$ is TU.
\end{definition}

Padberg observed that many almost TU matrices, including the one in \cref{figure:OTH}, can be brought into a `normalized' form by multiplying with a TU matrix; this normalized form is referred to as an unbalanced hole. 

\begin{definition}
    A $\{0,\pm 1\}$-matrix is a \textbf{hole} if it has exactly two nonzero entries in each row and column, and no proper submatrix has this property. 
    A hole is {\bf unbalanced} if its entries sum to $2 \mod 4$.
\end{definition}
Note that a matrix is a hole if and only if switching all $-1$s to $1$s yields the incidence matrix of a cycle.

Padberg initially conjectured that every almost TU matrix can be normalized by a transformation preserving the TU property.
However, Cornu\'{e}jols and Zuluaga demonstrated this was false, and updated the conjecture as follows:

\begin{conjecture}[Padberg~\cite{P1988}, Cornu\'{e}jols and Zuluaga~\cite{CZ2000}]\label{conjCorn}
Given an almost TU matrix $\mbfs{A}$, there exists a TU matrix $\mbfs{L}$ such that $\mbfs{L}\mbfs{A}$ is an unbalanced hole.
\end{conjecture}

In~\cref{secCorn} we prove that odd cycles, tree houses, and avocados are almost TU.
We also show that all three satisfy~\cref{conjCorn}.
Thus, as a consequence of~\cref{,thm:DirectedUniEquivalence,theorem:non_unimodular_laminar}, we obtain the following corollary.

\begin{corollary}\label{cor:CZConj}
\cref{conjCorn} holds if $\mbfs{A}$ or $\mbfs{A}^\top$ is the incidence matrix of a mixed disjoint or a laminar hypergraph.
\end{corollary}
%

%%%%%%%%%%%%%%%%%%%%%%%%%%%%%%
%%%%%%%%%%%%%%%%%%%%%%%%%%%%%%
\section{Preliminaries}\label{sec:Prelims}
%%%%%%%%%%%%%%%%%%%%%%%%%%%%%%
%%%%%%%%%%%%%%%%%%%%%%%%%%%%%%
%
Let $G = (V,E)$ be a hypergraph.
For $k \in \mbb{Z}_{\ge 1}$, we use $E_{k}(G)$ and $E_{\ge k}(G)$ to denote the sub-multisets of $E$ consisting of hyperedges with exactly and at least $k$ vertices, respectively. 
The multiset $E_{\ge 3}(G)$ contains the {\it proper hyperedges} of $G$.
If $E = E_2(G)$, then $G$ is a {\it graph}.
For $v \in V$, the multiset $\delta_G(v) \subseteq E$ denotes the hyperedges containing $v$, and $d_G(v) \coloneqq |\delta_G(v)|$ is its {\it degree}.
We use common graph terminology like {\it walks}, {\it paths}, and {\it cycles}; see, e.g.,~\cite{korte2011combinatorial}.
Following~\cite{korte2011combinatorial}, a walk in a graph may contain distinct parallel copies of an edge, however, it is not allowed to use the same copy repeatedly.
A path has no repeated vertices.
A path, walk, or cycle is {\it odd} (respectively, {\it even}) if it has an odd (respectively, even) number of edges.
Given a path $P$ and vertices $v,w\in V(P)$, we use $P[v,w]$ to denote the {\it $v$-$w$-subpath} of $P$.
We use $V_{\internal}(P)$ to denote the set of internal vertices of $P$.
A walk, path, or cycle in $G$ is a partial subhypergraph $H$ of $G$ such that $H$ is a graph that is a walk, path, or cycle, respectively.

The {\it union} of $G$ and another hypergraph $H$ is
\[
G+H \coloneqq \left(V(G)\cup V(H), ~E(G) \cupdot E(H)\right).
\]
For $e \subseteq V(G)$, we write $G+e\coloneqq G+(e,\{e\})$.
The {\it deletion} of a multiset $F \subseteq E(G)$ is
\[
G-F \coloneqq (V(G), E(G) \setminus F).
\]
For $e \in E(G)$, we write $G-e \coloneqq G - \{e\}$.

We say $G$ is {\it Eulerian} if each hyperedge has even size and each vertex has even degree.
Equivalently, $G$ is Eulerian if each row and column of $\mbfs{M}(G)$ contains an even number of $1$s.
The {\it support} of $\mbfs{M}(G)$ is 
\[
\supp\mbfs{M}(G)\coloneqq \{(v,e)\in V\times E\colon \mbfs{M}(G)_{v,e}=1\}.
\]
Our proof of Theorem~\ref{theorem:size_four} uses a TU characterization due to Camion. 
\begin{theorem}[Corollary of {\cite[Theorem 2]{Camion}}]\label{thm:camion}
Let $G$ be a hypergraph.
The following are equivalent:
	\begin{enumerate}[(i), leftmargin = *]
		\item \label{camion:1} $G$ is unimodular.
		\item \label{camion:2}$|\supp\mbfs{M}(H)| \equiv 0 ~\mathrm{mod}~ 4$ for each Eulerian partial subhypergraph $H$ of $G$ with $|V(H)|=|E(H)|$.
		\item \label{camion:3}$|\supp\mbfs{M}(H)| \equiv 0 ~\mathrm{mod}~ 4$ for each Eulerian partial subhypergraph $H$ of $G$.
	\end{enumerate}
	\label{theorem:restrict_to_Eulerian}
\end{theorem}

Camion's original result only asserts the equivalence of {\it\ref{camion:1}} and {\it\ref{camion:2}}.
Necessity of the (stronger) `rectangular' version in {\it\ref{camion:3}} can be obtained by appending $|V(H)|$ and $|E(H)|$ many all-zero rows and columns, respectively: The enlarged matrix $A$ is TU if and only if $\mbfs{M}(G)$ is, and for every Eulerian partial subhypergraph $H$ of $G$, we obtain an Eulerian square submatrix of $A$ with the same number of $1$s by adding zero rows or columns.

We set $[k] \coloneqq \{1, \dotsc, k\}$ for $k \in \mbb{Z}_{\ge1}$.
%
%%%%%%%%%%%%%%%%%%%%%%%%%%%%%%
%%%%%%%%%%%%%%%%%%%%%%%%%%%%%%
\section{Quasi-subhypergraphs}\label{ssecQuasi}
%%%%%%%%%%%%%%%%%%%%%%%%%%%%%%
%%%%%%%%%%%%%%%%%%%%%%%%%%%%%%
 
Our main result shows that a non-unimodular laminar hypergraph $G$ must contain an odd cycle, a tree house, or an avocado.
In order to prove this containment, we use induction on the size of the support of $\mbfs{M}(G)$.
As containment of submatrices is transitive, we can assume that $\mbfs{M}(G)$ is almost TU, so every proper partial subhypergraph of $G$ is unimodular.
As such, we cannot use partial subhypergraphs to apply the induction hypothesis (that non-TU matrices of smaller support contain an odd cycle, a tree house, or an avocado).
To work around this, we introduce a broader notion of subhypergraph, which we call a {\it quasi-subhypergraph}, and apply the induction hypothesis to a non-unimodular quasi-subhypergraph of smaller support.

Quasi-subhypergraphs are obtained by removing a vertex from a single hyperedge, or splitting a hyperedge into a hypermatching; see \cref{def:quasi_subhypergraph}.
The vertex removal allows us to reduce the support, while the splitting allows us to redistribute support.
Neither of these operations yields a partial subhypergraph.
As such, there are `conflicts' that prevent a quasi-subhypergraph from being a partial subhypergraph.
Important for our purposes is that if these conflicts are resolved, then the quasi-subhypergraph can be turned into a partial subhypergraph.
See \cref{ssecQuasiConf} for more.

This inductive proof technique, based on quasi-subhypergraphs, is applied in \cref{secUnimodular} for disjoint hypergraphs and \cref{sec:proof_main_thm} for laminar hypergraphs.
We believe that quasi-subhypergraphs may be of independent interest, so we highlight them in this section. 
%

%%%%%%%%%%%%%%%%%%%%%%%%%%%%%%
%%%%%%%%%%%%%%%%%%%%%%%%%%%%%%
\subsection{A formal definition of quasi-subhypergraphs}\label{ssecQuasiDefn}
%%%%%%%%%%%%%%%%%%%%%%%%%%%%%%
%%%%%%%%%%%%%%%%%%%%%%%%%%%%%%

We recall the notion of a partial subhypergraph that was informally introduced in \cref{sec:intro}.

\begin{definition}\label{def:partial}
Let $G$ be a hypergraph. 
We call a hypergraph $H$  with $V(H)\subseteq V(G)$ a {\bf partial subhypergraph} of $G$ if there is an injection $\iota\colon E(H)\hookrightarrow E(G)$ with $f=\iota(f)\cap V(H)$ for every $f\in E(H)$.

\noindent We write $H\pin G$ to denote that $H$ is a partial subhypergraph of $G$, and $(H,\iota)\pin G$ to express that $\iota$ witnesses that $H$ is a partial subhypergraph of $G$. 
\end{definition}
Schrijver's definition~\cite[\S77.1]{S2003V3} of partial subhypergraph considers $G[U,F]$.
Using a simple embedding, we observe that \cref{def:partial} generalizes Schrijver's definition.
\begin{observation}
	Let $G=(V,E)$ be a hypergraph and let $U\subseteq V$ and $F\subseteq E$. Define $\iota\colon F[U]\hookrightarrow E,f\cap U\mapsto f$. Then $(H,\iota)\pin G$.
\end{observation}
Recall that we say that a hypergraph $G$ \emph{contains} a hypergraph $H$ or that $H$ is \emph{in} $G$ if $H$ is \emph{isomorphic} to a partial subhypergraph of $G$.

Next, we consider maps $\Phi$ that relax the properties of $\iota$ in \cref{def:partial}, and, as such, induce quasi-subhypergraphs of $G$.
First, rather than having $f = \Phi(f) \cap V(H) $ for each $f \in E(H)$, we only require $f \subseteq \Phi(f) \cap V(H)$; this can be rewritten as $f \subseteq \Phi(f)$ because $f \subseteq V(H)$.
This relaxation allows quasi-subhypergraphs to be formed from $G$ by removing vertices from a hyperedge.
Second, rather than requiring $\Phi$ to be injective, we allow $\Phi^{-1}(e)$ to contain disjoint hyperedges for each $e \in E(G)$.
This allows quasi-subhypergraphs to be formed from $G$ by partitioning a hyperedge $e$ into a hypermatching.
We are now ready to define quasi-subhypergraphs:

\begin{definition}
A {\bf quasi-subhypergraph} of $G$ is a pair $(H,\Phi)$, where $H$ is a hypergraph with $V(H)\subseteq V(G)$ and $\Phi\colon E(H)\rightarrow E(G)$ satisfies

\smallskip

\begin{enumerate}[{\it(Q1)}, leftmargin = *]
	\item \label{def:quasi_1} For each $f\in E(H)$, we have $f\subseteq \Phi(f)$.
	\item \label{def:quasi_2} For each $e\in E(G)$, the hyperedges in $\Phi^{-1}(e)$ are pairwise disjoint.
\end{enumerate}
\smallskip

\noindent We write $(H,\Phi) \qin G$ to denote that $(H,\Phi)$ is a quasi-subhypergraph of $G$.
\label{def:quasi_subhypergraph}
\end{definition}

\begin{figure}
	\centering
	\scalebox{0.95}{
	\begin{NiceTabular}{c@{\hskip 2 cm}c}
	\begin{tikzpicture}[scale = .5, mynode/.style={circle, fill, minimum size = .1 cm, inner sep=0pt}]
	\begin{scope}[shift={(0cm,0 cm)}]
		\draw[draw = none](-2, -2.5) rectangle (1.85,.75);
		\node[mynode, draw = red, fill = red,  label={[above = 0 pt] \color{red}$u_0$}] (v0) at (0,0){};
		\node[mynode, label={[right = 0 pt] $u_1$}] (v1) at (1,-1){};
		\node[mynode, draw = red, fill = red,  label={[below right = 0 pt] \color{red}$u_2$}] (v2) at (.5, -2){};
		\node[mynode, draw = red, fill = red,  label={[below left = 0 pt] \color{red}$u_3$}] (v3) at (-.5,-2){};
		\node[mynode, label={[left = 0 pt] $u_4$}] (v4) at (-1, -1){};
		\draw[](v2) to (v1) to (v0) to (v4) to (v3) to (v2);
		\draw[thick, red](v2)to (v3);
		\node[] at (0, -3.5) {$H'$};
	\end{scope}
	
	\draw[ ->] (2.5,-1) to node[below, pos = .5]{$\Phi'$} (4.5,-1);
	
	\begin{scope}[shift={(6.75 cm,0 cm)}]
		\draw[draw = none](-2, -2.5) rectangle (1.85,.75);
		\node[mynode, draw = red, fill = red, label={[above = 0 pt] \color{red}$u_0$}] (v0) at (0,0){};
		\node[mynode, label={[right = 0 pt] $u_1$}] (v1) at (1,-1){};
		\node[mynode, draw = red, fill = red, label={[below right = 0 pt] \color{red}$u_2$}] (v2) at (.5, -2){};
		\node[mynode, draw = red, fill = red,  label={[below left = 0 pt] \color{red}$u_3$}] (v3) at (-.5,-2){};
		\node[mynode, label={[left = 0 pt] $u_4$}] (v4) at (-1, -1){};
		\node[] at (0, -1){\color{red}$e'$};
		\draw[draw = red, thick, rounded corners = 2pt] (0, .2) to[out = 0, in = 90] (.65, -2.15) to (-.65, -2.15) to[out = 90, in = 180] (0, .2);
		\draw[](v2) to (v1) to (v0) to (v4) to (v3);
		\node[] at (0, -3.5) {$G'$};
	\end{scope}
	\end{tikzpicture}
	&
	\begin{tikzpicture}[scale = .5, mynode/.style={circle, fill, minimum size = .1 cm, inner sep=0pt}]

	\begin{scope}[shift={(0 cm,0 cm)}]
		\draw[draw = none](-2, -2.5) rectangle (1.85,.75);
		\node[mynode, draw = red, fill = red, label={[above = 0 pt] \color{red}$v_0$}] (v0) at (.5,0){};
		\node[mynode, label={[right = 0 pt] $v_1$}] (v1) at (1,-1){};
		\node[mynode,draw = red, fill = red,  label={[below right = 0 pt] \color{red}$v_2$}] (v2) at (.5, -2){};
		\node[mynode, draw = red, fill = red, label={[below left = 0 pt] \color{red}$v_3$}] (v3) at (-.5,-2){};
		\node[mynode, label={[left = 0 pt] $v_4$}] (v4) at (-1, -1){};
		\node[mynode,draw = red, fill = red,  label={[above = 0 pt] \color{red}$v_5$}] (v5) at (-.5, 0){};
		\draw[](v0) to (v1) to (v2) to (v3) to (v4) to (v5) to (v0);
		\draw[thick, red](v0) to (v5);
		\draw[thick, red](v2) to (v3);
		\node[] at (0, -3.5) {$H''$};
	\end{scope}
	
	\draw[ ->] (2.5,-1) to node[below, pos = .5]{$\Phi''$} (4.5,-1);
	
	\begin{scope}[shift={(6.75 cm,0 cm)}]
		\draw[draw = none](-2, -2.5) rectangle (1.85,.75);
		\node[mynode, draw = red, fill = red, label={[above = 0 pt] \color{red}$v_0$}] (v0) at (.5,0){};
		\node[mynode, label={[right = 0 pt] $v_1$}] (v1) at (1,-1){};
		\node[mynode, draw = red, fill = red, label={[below right = 0 pt] \color{red}$v_2$}] (v2) at (.5, -2){};
		\node[mynode, draw = red, fill = red, label={[below left = 0 pt] \color{red}$v_3$}] (v3) at (-.5,-2){};
		\node[mynode, label={[left = 0 pt] $v_4$}] (v4) at (-1, -1){};
		\node[mynode,draw = red, fill = red,  label={[above = 0 pt] \color{red}$v_5$}] (v5) at (-.5, 0){};
		\node[] at (0, -1){\color{red}$e''$};
		\draw[thick, draw = red, rounded corners = 2pt] (.65, -2.15) to (-.65, -2.15) to (-.65, .2) to (.65, .2) to (.65, -2.15);
		\draw[](v0) to (v1) to (v2);
		\draw[](v3) to (v4) to (v5);
		\node[] at (0, -3.5) {$G''$};
	\end{scope}
	
	\end{tikzpicture}\\[.15cm]
	$\Phi'(f) \coloneqq 
	\left\{\begin{array}{l@{\hskip .1 cm}l}
	e' & \text{if}~f = \{u_2, u_3\}\\
	f &\text{otherwise}
	\end{array}\right.$
	&
	$\Phi''(f) \coloneqq 
	\left\{\begin{array}{l@{\hskip .1 cm}l}
	e'' & \text{if}~f = \{v_2, v_3\}\\
	e'' & \text{if}~f = \{v_0, v_5\}\\
	f &\text{otherwise}
	\end{array}\right.$
	\end{NiceTabular}}
\caption{Quasi-subhypergraphs $(H', \Phi') \qin G'$ (left) and $(H'', \Phi'') \qin G''$ (right).}
\label{fig:quasi_sh}
\end{figure}

%%%%%%%%%%%%%%%%%%%%%%%%%%%%%%
%%%%%%%%%%%%%%%%%%%%%%%%%%%%%%
\subsection{Conflicts preventing a quasi- from being a partial-subhypergraph}\label{ssecQuasiConf}
%%%%%%%%%%%%%%%%%%%%%%%%%%%%%%
%%%%%%%%%%%%%%%%%%%%%%%%%%%%%%

Every partial subhypergraph of $G$ constitutes a quasi-subhypergraph of $G$.
However, not every quasi-subhypergraph corresponds to a partial subhypergraph. 
Two examples of quasi-subhypergraphs are in~\cref{fig:quasi_sh}: $G'$ does not contain a 5-cycle as a partial subhypergraph, but a 5-cycle quasi-subhypergraph $H'$ can be obtained by removing $u_0$ from $e'$. 
Similarly, there is no 6-cycle in $G''$, but a 6-cycle quasi-subhypergraph $H''$ can be constructed by splitting $e''$ into edges $\{v_0,v_5\}$ and $\{v_2,v_3\}$.
The removal of vertices from a hyperedge or the splitting into a hypermatching conflicts with the definition of a partial subhypergraph.
Such conflicts are the difference between a quasi-subhypergraph and a partial subhypergraph.
\begin{definition}\label{def:Conflict}
Let $(H,\Phi) \qin G$. 
We say $e\in E(G)$ is a {\bf conflict} for $(H, \Phi)$ if there exists $f\in E(H)$ with $\Phi(f)=e$ and $f\subsetneq e\cap V(H)$.
We say $e$ is a conflict when $(H, \Phi)$ is clear from the context.
If there is no conflict, then $(H,\Phi)$ is {\bf conflict-free}.
\end{definition}          

Hyperedges are non-empty and, by \ref{def:quasi_2}, $\Phi^{-1}(e)$ is a hypermatching for each $e \in E(G)$.
Thus, we have the following characterization of a conflict:
\begin{observation}\label{prop:conflict_superset}
	Let $(H,\Phi) \qin G$.
	A hyperedge $e \in E(G)$ is a conflict if and only if $\Phi^{-1}(e)\neq \emptyset$ and $f\subsetneq e\cap V(H)$ for each $f\in \Phi^{-1}(e)$.
\end{observation}

Conflict-free quasi-subhypergraphs are precisely partial subhypergraphs:

\begin{lemma} \label{lemma:conflict_free_isomorphic_to_subhypergraph}
A quasi-subhypergraph $(H,\Phi)$ of $G$ is conflict-free if and only if $(H,\Phi)\pin G$.
\end{lemma}
\begin{proof}
Assume $(H, \Phi) \qin G$ is conflict-free.
By \cref{def:Conflict}, we have $f = \Phi(f) \cap V(H)$ for each $f \in E(H)$.
It remains to show that $\Phi$ is injective.
Assume for a contradiction that $\Phi(f) = \Phi(f')$ for some distinct $f, f' \in E(H)$.
Then $f' \subseteq \Phi(f) \cap V(H) \setminus f = \emptyset$ by Properties \ref{def:quasi_1} and \ref{def:quasi_2}.
However, $\emptyset \neq f'$ because hyperedges are non-empty.
This is a contradiction.

Conversely, if $(H,\Phi) \pin G$, then $\Phi(f) \cap V(H) = f$ for every $f \in E(H)$.
Thus, $(H, \Phi) \qin G$ is conflict-free by \cref{def:Conflict}.
\end{proof}

The relationship $(H, \Phi) \qin G$ is inherited by partial subhypergraphs of $H$ by restricting the function $\Phi$.
\begin{definition}
Let $(H,\Phi) \qin G$ and $H[U,F] \pin H$.
The {\bf restriction} of $\Phi$ to $H[U,F]$ is
$\Phi\upharpoonright_{H[U,F]}:E(H[U,F]) \to E(G)$, where $U\cap f \mapsto \Phi(f)$.
\end{definition}

Properties {\it\ref{def:quasi_1}} and {\it\ref{def:quasi_2}} for $\Phi$ carry over to its restriction $\Phi\upharpoonright_{H[U,F]}$.
Moreover, if $e \in E(G)$ is a conflict for $(H[U,F],\Phi\upharpoonright_{H[U,F]})$, then there exists $f \in E(H)$ with $\Phi(f) = e$ such that $f \cap U \subsetneq e \cap U$.
As $U \subseteq V(H)$, we have $f \subsetneq e \cap V(H)$.
Hence, $e$ is a conflict for $(H,\Phi)$.
Therefore, we have the following:

\begin{observation}
Let $(H,\Phi) \qin G$ and $H[U,F] \pin H$.
Then we have \\$(H[U,F],\Phi\upharpoonright_{H[U,F]}) \qin G$. 
Also, if $e\in E(G)$ is a conflict for $(H[U,F],\Phi\upharpoonright_{H[U,F]})$, then $e$ is a conflict for $(H,\Phi)$. \label{lem:restriction}
\end{observation}
\cref{lem:restriction} implies that a partial subhypergraph avoiding conflicts is conflict-free.

\begin{corollary}
Let $E_C\subseteq E(G)$ be the set of conflicts for $(H,\Phi) \qin G$. 
If $U \subseteq V(H)$ and $F \subseteq E(H)$ satisfies $F \cap \Phi^{-1}(E_C)=\emptyset$, then $H[U,F] \pin G$. 
\label{cor:restriction}	
\end{corollary}

As a partial inverse of the restriction operation, which removes vertices and hyperedges, we introduce an operation that introduces new hyperedges but not vertices.
\begin{definition}
Let $(H,\Phi) \qin G$ and $e\in E(G)\setminus \Phi(E(H))$ with $e\cap V(H) \neq \emptyset$. 
The {\bf addition of $e$ to $H$} is the pair $(H \oplus e,\Phi \oplus e)$, where 
\[
H \oplus e \coloneqq (V(H),~ E(H) \cupdot  \{e\cap V(H)\}),
\]
and $\Phi \oplus e \colon E(H \oplus e) \rightarrow E(G)$ with $(\Phi \oplus e)\upharpoonright_{E(H)}  = \Phi$ and $(\Phi \oplus e)(e\cap V(H))=e$.\label{def:addition}
\end{definition}

By~\cref{def:addition}, the hyperedge $e$ added to $H$ is not a conflict in $H \oplus e$.
Hence, any conflict in $H \oplus e$ originates from $H$.

\begin{observation}
Let $(H,\Phi) \qin G$ and $e\in E(G)\setminus \Phi(E(H))$ with $e\cap V(H)\neq \emptyset$. 
Then $(H \oplus e,\Phi \oplus e) \qin G$. 
Moreover, $g\in E(G)$ is a conflict for $(H \oplus e,\Phi \oplus e)$ if and only if $g$ is a conflict for $(H,\Phi)$. \label{lem:addition}	
\end{observation}

If $G$ has no odd cycles, then we can use the addition operation to infer structure within a conflict-free quasi-subhypergraph. 
We see this, for instance, in the next lemma, which we often use in arguments to determine the parity of paths and cycles.
\begin{lemma}
Assume $G$ does not contain an odd cycle, and let $(H,\Phi) \pin G$.\label{lemma:parities}
\begin{enumerate}[{\it (i)}, leftmargin = *]
\item \label{lemma:parities:1} Assume $H$ is an $a$-$b$-walk with $a\neq b$. 
If there exists $e\in E(G)\setminus \Phi(E(H))$ such that $e\cap V(H)=\{a,b\}$, then $|E(H)|$ is odd.
\item \label{lemma:parities:2} Assume $H$ is the edge-disjoint union of an $a$-$b$-walk $P$ and a $c$-$d$-walk $Q$, where $a,b,c,d$ are distinct. 
If there exist $e,f\in E(G)\setminus \Phi(E(H))$ such that $e\cap V(H)=\{a,c\}$ and $f\cap V(H)=\{b,d\}$, then $|E(H)|=|E(P)|+|E(Q)|$ is even.
\end{enumerate}
\end{lemma}
\begin{proof}
{\it \ref{lemma:parities:1}}: $H \oplus e\pin G$ by \cref{lem:addition} and \cref{lemma:conflict_free_isomorphic_to_subhypergraph}.
Moreover, $H\oplus e$ is a closed walk. 
As $G$ contains no odd cycle, $|E(H\oplus e)|=|E(H)|+1$ is even.

\noindent{\it\ref{lemma:parities:2}}: $H\oplus e\oplus f \pin G$ by \cref{lem:addition} and \cref{lemma:conflict_free_isomorphic_to_subhypergraph}.
Also, $H\oplus e\oplus f$ is a closed walk. 
As $G$ contains no odd cycle, $|E(H\oplus e\oplus f)|=|E(P)|+|E(Q)|+2$ is even. 
\end{proof}

%%%%%%%%%%%%%%%%%%%%%%%%%%%%%%
%%%%%%%%%%%%%%%%%%%%%%%%%%%%%%
\section[A proof of Theorem~\ref{theorem:size_four}]{A proof of \cref{theorem:size_four}}\label{secUnimodular}
%%%%%%%%%%%%%%%%%%%%%%%%%%%%%%
%%%%%%%%%%%%%%%%%%%%%%%%%%%%%%
%
We repeat the theorem statement for convenience. 

\TheoremSizeFour*
The necessary condition for unimodularity in \cref{theorem:size_four} is easy to prove.
Indeed, if $H$ is an odd cycle of length $2k+1$ for $k \in \mbb{Z}_{\ge 1}$, then $|\supp\mbfs{M}(H)| = 4k+2\equiv 2\mod 4$.
Similarly, if $H$ is a tree house with hyperedge $\{r, \ell_1, \ell_2, \ell_3\}$ of size $4$ and odd paths $P_1$, $P_2$, and $P_3$, then $|\supp\mbfs{M}(H)| = 4+2 \cdot (|E(P_1)|+|E(P_2)|+|E(P_3)|)\equiv 2\mod 4$.
In both cases, $H$ is Eulerian.
Thus, if $G$ is a unimodular hypergraph, then it can neither contain an odd cycle nor a tree house by \cref{theorem:restrict_to_Eulerian}.

It remains to prove the sufficient condition, that is, if $G$ is a non-unimodular disjoint hypergraph, then it contains an odd cycle or a tree house.
Proving this direction is more involved than proving the necessary condition.
The remainder of this section gives an overview, with technical details deferred to \cref{sec:existence_almost_nice_cycle,ssecquasiG*,subsubsec:no_odd_tree_house,subsec:no_odd_cycle}.

We prove the sufficient condition by induction on the support size $|\supp \mbfs{M}(G)|$.
The result holds vacuously if $|\supp \mbfs{M}(G)| = 0$ because $E(G)$ contains only non-empty hyperedges by our definition of a hypergraph.
Therefore, assume that for $t \ge 0$, any non-unimodular disjoint hypergraph $G$ with $|\supp \mbfs{M}(G)| \le t$ contains an odd cycle or a tree house.

Let $G^*$ be a non-unimodular disjoint hypergraph satisfying $|\supp \mbfs{M}(G^*)| = t+1$.
As $G^*$ is non-unimodular, \cref{thm:camion}\ref{camion:2} guarantees an Eulerian
partial subhypergraph $H$ of $ G^*$ with $|V(H)| = |E(H)|$ and
$|\supp \mbfs{M}(H)| \equiv 2 \bmod 4$.
By choosing $H$ to be minimal with this property, we may assume that $H$ contains no isolated vertices.
If $|\supp \mbfs{M}(H)| \le t$, then the induction hypothesis applied to $H$ implies that $H$ contains an odd cycle or a tree house.
The partial subhypergraph relation is transitive, so this substructure is also a partial subhypergraph of $G^*$, giving the result.
Otherwise, $H$ is a partial subhypergraph of $G^*$ with $|\supp \mbfs{M}(H)| = t+1 = |\supp \mbfs{M}(G^*)|$, 
so we may replace $G^*$ by $H$ and assume
\begin{equation}\label{eqG*1}
\begin{array}{l}
\text{$G^*$ is Eulerian with $|V(G^*)| = |E(G^*)|$, $|\supp \mbfs{M}(G^*)| \equiv 2\mod 4$,}\\
\text{and $G^*$ does not have isolated vertices.}
\end{array}
\end{equation}
The inductive step is complete if $G^*$ contains an odd cycle.
Thus, we assume that 
\begin{equation}\label{eqG*3}
\text{$G^*$ does not contain an odd cycle.}
\end{equation}
To use the induction hypothesis, we need to reduce $|\supp \mbfs{M}(G^*)|$ while preserving disjointness and non-unimodularity.
The latter can be achieved by preserving Camion's conditions in~\eqref{eqG*1}; a natural way of doing this is by removing an even cycle.
Indeed, consider the following.
Suppose 
\[
C^* \coloneqq G^*[U_{C^*},F_{C^*}]
\]
is an even cycle in $G^*$.
Let $H^*$ be the hypergraph whose incidence matrix $\mbfs{M}(H^*)$ is obtained by deleting the support of $\mbfs{M}(C^*)$ from $\mbfs{M}(G^*)$, and removing all-zero rows and columns:

\begin{equation}\label{eq:def_H}
\begin{array}{rcl}
V(H^*)&\coloneqq & V(G^*) \setminus \{v\in U_{C^*}: \delta_{G^*}(v) \subseteq F_{C^*}\} \\[.1cm]
E(H^*)&\coloneqq & ( E(G^*)\setminus F_{C^*} ) \cupdot \{e\setminus U_{C^*}:  e\in F_{C^*} \cap E_{\geq4}(G^*)\}.
\end{array}
\end{equation}

For each row and column of $\mbfs{M}(G^*)$, we change either zero or two 1s to 0s in order to obtain $\mbfs{M}(H^*)$.
Hence, $H^*$ is Eulerian and disjoint.
Also, as $|\supp\mbfs{M}(C^*)| \equiv 0 \mod 4$, we have $|\supp \mbfs{M}(H^*)| \equiv |\supp \mbfs{M}(G^*)| \equiv 2 \mod 4$.
Moreover, 
\[
|\supp \mbfs{M}(H^*)| = |\supp \mbfs{M}(G^*)| - |\supp \mbfs{M}(C^*)| \le |\supp \mbfs{M}(G^*)| - 4.
\]

Hence, $H^*$ is non-unimodular and $|\supp \mbfs{M}(H^*)|  < |\supp \mbfs{M}(G^*)| $.
If we can find $C^*$ such that $H^*$ is a partial subhypergraph of $G^*$, then we can apply the induction hypothesis to conclude that $H^*$, and, hence, $G^*$, contains an odd cycle or a tree house; this would complete the induction step.

We observe that 
\begin{equation}\label{eqSpecialCond}
\begin{array}{l}
\text{if $\delta_{G^*}(v)\subseteq F_{C^*}$ for each $e \in F_{C^*} \cap E_{\ge 4}(G^*)$ and $v\in e\cap U_{C^*}$,}\\
\text{then $H^*$ is a partial subhypergraph of $G^*$}.
\end{array}
\end{equation}
If the prerequisite of \eqref{eqSpecialCond} holds, then $\mbfs{M}(H^*)$ is obtained by deleting the rows and columns of $\mbfs{M}(G^*)$ corresponding to $\cup_{e \in F_{C^*} \cap E_{\ge 4}(G^*)} e \cap U_{C^*}$ and $F_{C^*} \cap E_{2}(G^*)$, respectively.
The prerequisite of \eqref{eqSpecialCond} holds, e.g., when $C^*$ only uses edges from $G^*$, i.e., if $F_{C^*} \cap E_{\ge 4}(G^*)= \emptyset$.
Thus, we assume
\begin{equation}\label{eqG*2}
\text{the graph $(V(G^*), E_2(G^*))$ is a forest.}
\end{equation}

Although $G^*$ may not satisfy the prerequisite of \eqref{eqSpecialCond} for each $e \in F_{C^*} \cap E_{\ge 4}(G^*)$, e.g., if $G^*$ is a tree house, we can choose $C^*$ such that at most one hyperedge $g^* \in F_{C^*} \cap E_{\ge 4}(G^*)$ violates this condition.
We prove the following result in~\cref{sec:existence_almost_nice_cycle}.
The additional property \ref{eq:nice_and_almost_nice} disallows shortcutting the cycle using proper hyperedges.

\begin{lemma}\label{lemma:existence_almost_nice_cycle}
$G^*$ contains an even cycle $G^*[U,F]$ satisfying:
\begin{enumerate}[{\it (N1)}, leftmargin = *]
\item\label{eq:almost_nice} There exists $g^*\in F$ such that for each $e\in F\cap E_{\geq 4}(G^*)\setminus\{g^*\}$ and each $v\in e\cap U$, we have $\delta_{G^*}(v)\subseteq F$. 
\item\label{eq:nice_and_almost_nice} For each $e\in E_{\geq 4}(G^*)\setminus F$, we have $|e\cap U|\leq 1$. 
\end{enumerate}
\end{lemma}

We now assume the cycle $C^* = G^*[U_{C^*},F_{C^*}]$ meets Condition \ref{eq:almost_nice}, with a hyperedge $g^* \in F_{C^*}$, and Condition \ref{eq:nice_and_almost_nice}.
If either $g^* \in F_{C^*} \cap E_{2}(G^*)$ or if $\delta_{G^*}(v) \subseteq F_{C^*}$ for all $v\in g^*\cap U_{C^*}$, then the prerequisite of \eqref{eqSpecialCond} is met and $H^*$ is a partial subhypergraph of $G^*$, which completes the proof by induction.
Thus, assume 
\begin{equation} g^* \in F_{C^*} \cap E_{\ge 4}(G^*)\text{ and }\delta_{G^*}(v^*) \not\subseteq F_{C^*}\text{ for some }v^*\in g^*\cap U_{C^*}.\label{eq:g_star}\end{equation}
Then $\mbfs{M}(H^*)_{v^*, g^*} = 0$ but $\mbfs{M}(G^*)_{v^*, g^* } =1$, so $\mbfs{M}(H^*)$ is not a submatrix of $\mbfs{M}(G^*)$ because $v^*$ has been removed from $g^*$.
In \cref{ssecquasiG*}, we show $H^*$ is a quasi-subhypergraph of $G^*$ with a unique conflict, namely $g^*$.

As $|\supp \mbfs{M}(H^*)| < |\supp \mbfs{M}(G^*)|$, we can apply induction on $H^*$ to conclude it contains an odd cycle or a tree house.
By \eqref{eqG*3}, $G^*$ does not contain an odd cycle.
Thus, our final step is to prove the following lemmas:

\begin{restatable}{lemma}{LemmaNoOddTreeHouse}\label{lemma:no_odd_tree_house}
If $H^*$ contains a tree house, then $G^*$ contains a tree house.
\end{restatable}
\begin{restatable}{lemma}{LemmaNoOddCycleInH}\label{lemma:no_odd_cycle_in_H}
If $H^*$ contains an odd cycle, then $G^*$ contains a tree house. 
\end{restatable}

Proving \cref{lemma:no_odd_tree_house,lemma:no_odd_cycle_in_H} completes the inductive step, and thus finishes the proof of \cref{theorem:size_four}.
\cref{lemma:no_odd_tree_house,lemma:no_odd_cycle_in_H} are proved in~\cref{subsubsec:no_odd_tree_house,subsec:no_odd_cycle}, respectively.
%

%%%%%%%%%%%%%%%%%%%%%%%%%%%%%%
%%%%%%%%%%%%%%%%%%%%%%%%%%%%%%
\section[A proof of Lemma~\ref{lemma:existence_almost_nice_cycle}]{A proof of~\cref{lemma:existence_almost_nice_cycle}}\label{sec:existence_almost_nice_cycle}
%%%%%%%%%%%%%%%%%%%%%%%%%%%%%%
%%%%%%%%%%%%%%%%%%%%%%%%%%%%%%
%
The hypergraph $G^*$ is Eulerian and $|V(G^*)| = |E(G^*)|$ by~\eqref{eqG*1}.
Moreover, $G^*$ does not contain empty hyperedges.
Thus,
\begin{equation}
|E_2(G^*)|=|E(G^*)|-|E_{\geq 4}(G^*)|=|V(G^*)|-|E_{\geq 4}(G^*)|. \label{eq:num_conn_comp}
\end{equation}

Set $V^*_{\ge 4}\coloneqq \cup_{e \in E_{\geq 4}(G^*)} e$. 
As $G^*$ is Eulerian and disjoint, each vertex in $V^*_{\ge 4}$ has odd degree in $(V(G^*),E_2(G^*))$, whereas each vertex in $V(G^*)\setminus V^*_{\ge 4}$ has even degree in $(V(G^*),E_2(G^*))$.
The graph $(V(G^*),E_2(G^*))$ is a forest by~\eqref{eqG*2}. 
Let $L^*\subseteq V^*_{\ge 4}$ be the set of leaves (i.e., vertices of degree $1$) in $(V(G^*),E_2(G^*))$.
Each vertex in $V^*_{\ge 4}\setminus L^*$ has degree at least $3$ in $(V(G^*),E_2(G^*))$.
Each vertex in $V(G^*)\setminus V^*_{\ge 4}$ has degree at least $2$ in $(V(G^*),E_2(G^*))$ because $G^*$ is Eulerian and $\mbfs{M}(G^*)$ has no all-zero rows by \eqref{eqG*1}.
Thus,
\begin{align*}
2\cdot (|V(G^*)|-|E_{\geq 4}(G^*)|)& =2\cdot |E_2(G^*)|=\sum_{v\in V(G^*)} d_{(V(G^*),E_2(G^*))}(v) \\
&\geq |L^*|+3\cdot (|V^*_{\ge 4}|-|L^*|)+2\cdot (|V(G^*)|-|V^*_{\ge 4}|)\\
&=2\cdot (|V(G^*)|-|L^*|)+|V^*_{\ge 4}|.
\end{align*}

Hence, $2\cdot (|L^*|-|E_{\geq 4}(G^*)|) \geq |V^*_{\ge 4}|$.
The hyperedges in $E_{\geq 4}(G^*)$ are pairwise disjoint and contain at least $4$ vertices each, so
\begin{equation}
\frac{|L^*|-|E_{\geq 4}(G^*)|}{2}\geq \frac{|V^*_{\ge 4}|}{4}\geq |E_{\geq 4}(G^*)|\label{eq:lower_bound_leaves}.
\end{equation}

We now construct an auxiliary graph $\overline{G}$ with vertex set $V(\overline{G})=V(G^*)$ as follows: 
Add $E_2(G^*)$ to $E(\overline{G})$; refer to these edges as \emph{forest edges}. 
For each $e\in E_{\geq 4}(G^*)$, add to $E(\overline{G})$ a perfect matching (respectively, a near-perfect matching)
on the leaves $L^* \cap e$ if $|L^* \cap e|$ is even (respectively, odd).
Refer to these edges as \emph{matching edges}. 
Define $\Phi: E(\overline{G}) \to E(G^*)$ as follows:
\[
\Phi(f) := \left\{
\begin{array}{ll}
e, & \text{if $f$ is a matching edge and $f \subseteq e \in E_{\ge 4}(G^*)$}\\
f, & \text{if $f$ is a forest edge}
\end{array}\right..
\]
For each $e \in E(G^*)$, the set $\Phi^{-1}(e)$ is a matching in $\overline{G}$. 
Thus, $(\overline{G}, \Phi) \qin G^*$.

The matching edges cover all but at most $|E_{\geq 4}(G^*)|$ leaves.
Thus, there are at least $(|L^*|-|E_{\geq 4}(G^*)|)/2$ matching edges.
Hence,
\begin{align*}
|E(\overline{G})|\geq |E_2(G^*)|+\frac{|L^*|-|E_{\geq 4}(G^*)|}{2}& \stackrel{\cref{eq:num_conn_comp}}{=}|V(G^*)|-|E_{\geq 4}(G^*)|+\frac{|L^*|-|E_{\geq 4}(G^*)|}{2}\\
&\stackrel{\cref{eq:lower_bound_leaves}}{\geq} |V(G^*)|=|V(\overline{G})|.
\end{align*}

The previous inequalities imply there exists a cycle $\overline{C} \pin \overline{G}$. 
We now show how to turn $\overline{C}$ into an even cycle $C^* \pin G^*$ satisfying Properties \ref{eq:almost_nice} and \ref{eq:nice_and_almost_nice}.
To this end, we say that distinct vertices $v,w\in V(\overline{C})$ are \emph{crossable} if there exists $e\in E_{\geq 4}(G^*)$ such that $v,w\in e$, but $\{v,w\}$ is not a matching edge contained in $\overline{C}$.
If there is a crossable pair, let $v$ and $w$ be two crossable vertices with a shortest $v$-$w$-path in $\overline{C}$, and let $\overline{P}$ be a shortest $v$-$w$-path in $\overline{C}$.
If there is no pair of crossable vertices in $V(\overline{C})$, then select an arbitrary pair $v,w$ of vertices that are connected in $\overline{C}$ by a matching edge $e_0$ (such a pair exists because $(V(G^*),E_2(G^*))$ is a forest), and let $\overline{P}\coloneqq \overline{C}-e_0$.
In both cases, let $g^*  \in E_{\geq 4}(G^*)$ be the proper hyperedge containing $v$ and $w$, which is unique by disjointness of $E_{\geq 4}(G^*)$.
Whether a crossable pair exists or not, $\overline{P}$ is a $v$-$w$-subpath of $\overline{C}$ satisfying
\begin{equation}\label{barP:3} 
\begin{array}{l}
	\text{$V(\overline{P})$ contains no pair of crossable vertices other than possibly $\{v,w\}$.}
	\end{array}
\end{equation}

Statement~\eqref{barP:3} implies that 
\begin{equation}\label{eq:leq2}
\text{$|e \cap V(\overline{P})|\le 2$ for each $e \in E_{\ge 4}(G^*)$ and $g^* \notin \Phi(E(\overline{P}))$.}
\end{equation}
Indeed, if $e \cap V(\overline{P})$ were to contain three vertices, then the middle one can only be contained in a matching edge with one of the outer two, i.e., we would get a crossable pair other than $\{v,w\}$.
Applying this bound to $g^*$, together with $v, w \in g^* \cap V(\overline{P})$, gives $g^* \cap V(\overline{P}) = \{v, w\}$.
Any matching edge $f \in E(\overline{P})$ with $\Phi(f) = g^*$ would therefore satisfy $f \subseteq \{v, w\}$, forcing $f = \{v, w\}$.
However, $\{v, w\} \notin E(\overline{P})$ because in the crossable case $\{v,w\}$ is not a matching edge of $\overline{C}$ by definition, and in the other case $\{v,w\} = e_0 \notin E(\overline{P}$) and $e_0$ is the unique matching edge connecting $v$ to $w$.
Hence, $g^* \notin \Phi(E(\overline{P}))$.

In fact, we claim that the following holds for all $e \in E_{\ge 4}(G^*)$:
\begin{equation}\label{eq:equality_2}
\text{$|e \cap V(\overline{P})| = 2$ if and only if $e \in \Phi(E(\overline{P})) \cup \{g^*\}$}.
\end{equation}

To prove \eqref{eq:equality_2}, let $e \in E_{\ge 4}(G^*)$.
By~\eqref{eq:leq2}, we have $|e \cap V(\overline{P})|\le 2$.
If $e=\Phi(f)\in\Phi(E(\bar{P}))$, then $f\subseteq e\cap V(\bar{P})$, so $|e \cap V(\overline{P})|= 2$. Moreover, $\{v,w\}\subseteq g^*\cap V(\bar{P})$, so also $|g^* \cap V(\overline{P})|= 2$.
If $|e \cap V(\overline{P})| \le 1$, then $e \cap V(\overline{P})$ does not form an edge of $\overline{P}$ (as it does not contain enough vertices) nor is it equal to $g^*$ (as $g^*$ contains the two end vertices of $\overline{P}$).
Hence, $e \not\in \Phi(E(\overline{P})) \cup \{g^*\}$.
Conversely, suppose $e \cap V(\overline{P}) = \{x,y\}$ has size $2$ and $e \neq g^*$.
The vertices $x$ and $y$ must form a matching edge in $\overline{C}$, otherwise we have a crossable pair contradicting \eqref{barP:3}.
As $G^*$ is disjoint, $e$ is the only proper hyperedge containing $x$ and $y$; thus, $e = \Phi(\{x,y\})$.
Moreover, the assumption $e \neq g^*$ implies $\{x,y\} \cap \{v,w\} = \emptyset$ because $G^*$ is disjoint.
Hence, $x, y\in V_{\internal}(\overline{P})$, so $\{x,y\}\in E(\overline{P})$, implying $e=\Phi(\{x,y\})\in \Phi(E(\overline{P}))$. 
This completes the proof of~\eqref{eq:equality_2}.

We use $\overline{P}$ to construct the even cycle $C^*$ satisfying Properties \ref{eq:almost_nice} and \ref{eq:nice_and_almost_nice}.
We start by showing $(\overline{P},\Phi\upharpoonright_{\overline{P}})$ is conflict-free.
Assume to the contrary that $e \in E(G^*)$ is a conflict.
Let $f \in E(\overline{P})$ be such that $\Phi\upharpoonright_{\overline{P}}(f) = e$ and $f \subsetneq e \cap V(\overline{P})$.
As $f$ is an edge, we have $|e \cap V(\overline{P})| \ge 3$ and $e \in E_{\ge 4}(G^*)$, which contradicts~\eqref{eq:leq2}.

Finally, using that $g^* \notin \Phi(E(\overline{P}))$ by~\eqref{eq:leq2}, define the cycle $C^* := \overline{P}\oplus g^*$.
As $\overline{P}$ is conflict-free,~\cref{lem:addition} and \cref{lemma:conflict_free_isomorphic_to_subhypergraph} imply that $C^* \pin G^*$.
Moreover, $C^*$ is isomorphic to $G^*[U,F]$ with $U \coloneqq V(\overline{P})$ and $F \coloneqq \Phi(E(\overline{P})) \cup \{g^*\}$.
The cycle $C^* $ is even by~\eqref{eqG*3}.
For Property \ref{eq:almost_nice}, observe that each $e\in F\cap E_{\geq 4}(G^*)\setminus \{g^*\}$ is of the form $\Phi(f)$ for a matching edge $f\in E(\overline{P})$. 
Each endpoint $u$ of $f$ lies in $L^* \cap e$, so $u$ has degree $1$ in $(V(G^*),E_2(G^*))$ and belongs to exactly one proper hyperedge.
This implies that $u$ has degree $2$ in $G^*$.
Thus, $\delta_{G^*}(u)\subseteq F$ because $C^*$ is a cycle, establishing Property~\ref{eq:almost_nice}.
Property \ref{eq:nice_and_almost_nice} is satisfied by~\eqref{eq:leq2} and~\eqref{eq:equality_2}.
%

%%%%%%%%%%%%%%%%%%%%%%%%%%%%%%
%%%%%%%%%%%%%%%%%%%%%%%%%%%%%%
\section{Showing $H^* $ is a quasi-subhypergraph of $G^*$}\label{ssecquasiG*}
%%%%%%%%%%%%%%%%%%%%%%%%%%%%%%
%%%%%%%%%%%%%%%%%%%%%%%%%%%%%%
Recall $G^*$, $H^*$, and $C^* = G^*[U_{C^*}, F_{C^*}]$ from \cref{eq:def_H}. 
Let $g^*\in F_{C^*}$ be as in \cref{lemma:existence_almost_nice_cycle}.
Define $\Phi^*:E(H^*)\rightarrow E(G^*)$ by
\begin{equation} \label{eq:def_phi_H}
\Phi^*(f)\coloneqq \left\{\begin{array}{ll}
f & \text{if}~f\in E(G^*)\setminus F_{C^*}\\
e & \text{if}~f=e\setminus U_{C^*}~\text{for some}~e\in F_{C^*} \cap E_{\geq 4}(G^*).
\end{array}\right.
\end{equation} 
By the definitions in \cref{eq:def_H} and \cref{eq:def_phi_H}, we have, for each $e \in E(G^*)$:
\begin{equation}\label{eq:phi_H_preimage}
	(\Phi^*)^{-1}(e)=\left\{\begin{array}{cl}
		\{e\} & \text{if}~e\in E(G^*)\setminus F_{C^*}\\
		\{e\setminus U_{C^*}\} & \text{if}~e\in F_{C^*} \cap E_{\geq 4}(G^*)\\
		\emptyset &\text{if}~e \in F_{C^*} \cap E_{2}(G^*).
	\end{array}\right.
\end{equation}
Equations~\cref{eq:def_phi_H,eq:phi_H_preimage} imply Properties~{\it\ref{def:quasi_1}} and~{\it\ref{def:quasi_2}} for $(H^*,\Phi^*)$, respectively. 

\begin{observation}\label{lem:H_plus_C_quasi}
   We have $(H^*,\Phi^*) \qin G^*$. 
\end{observation}
\begin{lemma}\label{lem:only_conflict_g^*}
For each $e \in E(G^*) \setminus \{g^*\}$, we have
\begin{equation}\label{eq:def_phi_H_conflict}
(\Phi^*)^{-1}(e)=\left\{\begin{array}{cl}
\emptyset & \text{if}~e\in F_{C^*}\cap E_2(G^*)\\
\{e\cap V(H^*)\} & \text{otherwise}.
\end{array}\right.
\end{equation}
In particular, $g^*$ is the only conflict for $(H^*, \Phi^*)$.
\end{lemma}
\begin{proof}
Let $e \in E(G^*) \setminus \{g^*\}$.
First, consider the case $e\in E(G^*)\setminus F_{C^*}$.
Then $e \in E(H^*)$ by \eqref{eq:def_H}. 
Hence, \cref{eq:phi_H_preimage} implies that $(\Phi^*)^{-1}(e)=\{e\}=\{e\cap V(H^*)\}$.

Next, assume $e\in F_{C^*}$.
As $C^*$ is a cycle, $|e\cap U_{C^*}|=2$. 
Note that $E(G^*)=E_2(G^*)\cup E_{\geq 4}(G^*)$ because $G^*$ is Eulerian.
If $e\in E_2(G^*)$, then $(\Phi^*)^{-1}(e)=\emptyset$ by \cref{eq:phi_H_preimage}.
If $e\in E_{\geq 4}(G^*)$, then $e\setminus U_{C^*}\neq \emptyset$.
Also, $(\Phi^*)^{-1}(e)=\{e\setminus U_{C^*}\}$ by \cref{eq:phi_H_preimage}.
By \cref{eq:def_H}, we have $e\setminus U_{C^*}\subseteq e\cap V(H^*)$. 
By Property \ref{eq:almost_nice}, we further know $\delta_{G^*}(v)\subseteq F_{C^*}$ for every $v\in e\cap U_{C^*}$, so $v \notin V(H^*)$ by \cref{eq:def_H}.
Thus, $e\cap V(H^*)\subseteq e\setminus U_{C^*}$.
Hence, $(\Phi^*)^{-1}(e)=\{e\setminus U_{C^*}\}=\{e\cap V(H^*)\}$.

As outlined in \cref{secUnimodular},  $H^* \not\pin G^*$ because we removed non-zero entries from 
$\mbfs{M}(G^*)_{\cdot, g^*}$.
Thus, $g^*$ is a conflict for $(H^*, \Phi^*)$, and it is the only conflict by \eqref{eq:def_phi_H_conflict}.
\end{proof}

%%%%%%%%%%%%%%%%%%%%%%%%%%%%%%
%%%%%%%%%%%%%%%%%%%%%%%%%%%%%%
\section[A proof of Lemma~\ref{lemma:no_odd_tree_house}]{A proof of \cref{lemma:no_odd_tree_house}}\label{subsubsec:no_odd_tree_house}
%%%%%%%%%%%%%%%%%%%%%%%%%%%%%%
%%%%%%%%%%%%%%%%%%%%%%%%%%%%%%
%
We repeat the statement for convenience. 

\LemmaNoOddTreeHouse*
Let $G^*$ and $H^*$ be as in~\eqref{eqG*1} and \eqref{eq:def_H}, and let $\Phi^*$ be as in~\cref{eq:def_phi_H}.
By the assumption in~\cref{lemma:no_odd_tree_house}, there is a tree house $T^* \pin H^*$.
By \cref{lem:only_conflict_g^*,lem:restriction}, $(T^*,\Phi^*\upharpoonright_{T^*}) \qin G^*$ has at most one conflict.
We may assume that $g^*$ actually constitutes a conflict for $(T^*,\Phi^*\upharpoonright_{T^*})$  as $T^*\pin G^*$ otherwise.
Moreover, $\Phi^*\upharpoonright_{T^*}$ is injective by \cref{eq:phi_H_preimage}.
Using \eqref{eqG*3}, \cref{lemma:no_odd_tree_house} follows from \cref{lemma:odd_tree_house}, which applies to the more general laminar setting and will be reused in the proof of \cref{theorem:non_unimodular_laminar}.
We remark that \ref{item:laminar_odd_tree_house_3} holds in the setting of \cref{lemma:no_odd_tree_house} by our assumption that $g^*\in E_{\ge 4}(G^*)$ is the (unique) conflict and because $E_{\ge 4}(G^*)$ is a disjoint set family, so every set is maximal.
For \ref{item:laminar_odd_tree_house_4}, we note that since $\Phi^*$ is injective, $f_c\ne h$ implies $g^*=e_c=\Phi^*(f_c)\ne \Phi(h)$. 
As $|h|=4$, we have $\Phi(h)\in E_{\ge 4}(G^*)$. 
As $g^*\in E_{\ge 4}(G^*)\setminus \{\Phi(h)\}$, we have that $\Phi(h)$ and $g^*$ are disjoint. 
In particular, $f_c\cap h\subseteq g^*\cap \Phi(h)=\emptyset$.

\begin{restatable}{lemma}{LemmaOddTreeHouse}\label[lemma]{lemma:odd_tree_house}
	Let $G$ be an Eulerian laminar hypergraph and let $\mathcal{L}\coloneqq E_{\ge 4}(G)$. Suppose that there is a quasi-subhypergraph $(T,\Phi)\qin G$ with the following properties:
	\begin{enumerate}[{\it({TQ}1)}, leftmargin = *]
		\item \label{item:laminar_odd_tree_house_1}$T$ is a tree house with paths $(P_i)_{i\in [3]}$ and proper hyperedge $h=\{r,\ell_1,\ell_2,\ell_3\}$.
		\item \label{item:laminar_odd_tree_house_2}$\Phi$ is injective.
		\item \label{item:laminar_odd_tree_house_3}There is a unique conflict $e_c\in E(G)$, which constitutes a maximal set in $\mathcal{L}$.
		\item \label{item:laminar_odd_tree_house_4}Let $f_c\in E(T)$ with $\Phi(f_c)=e_c$. If $f_c\ne h$, then $f_c\cap h=\emptyset$.
	\end{enumerate}
	Then $G$ contains an odd cycle or a tree house.
\end{restatable}
\begin{proof}
We may assume without loss of generality that

\begin{equation}
	\text{$G$ does not contain an odd cycle} \label{eq:non_odd_cycle_2}
\end{equation}
as we are done otherwise.

We distinguish two cases, depending on the interaction between $e_c$ and $\Phi(h)$. Note that $\Phi(h)\in E_{\ge 4}(G)=\mathcal{L}$ because $h\subseteq \Phi(h)$ and $|h|=4$. As $e_c$ is a maximal set in $\mathcal{L}$, $\Phi(h)$ can either be equal to or a descendant of $e_c$, or we have $\Phi(h)\cap e_c=\emptyset$.

\subsubsection*{Case 1: $\Phi(h)\in \{e_c\}\cup \descendants_{\mathcal{L}}(e_c)$} 
In this case, we have $h\subseteq \Phi(h)\subseteq e_c$. For $i\in [3]$, let $v_i$ be the first vertex on $P_i$ after $r$ that is contained in $e_c$. Note that $v_i$ exists because $\ell_i$ is contained in $e_c$.  For $i\in [3]$, let $P'_i\coloneqq P_i[r,v_i]$. 
\begin{claim}
	For $i\in [3]$, we have $f_c\notin E(P'_i)$. In particular, $e_c\notin \Phi(E(P'_i))$.
\end{claim}
\begin{cpf}
	If $f_c=h$, then clearly $f_c\notin E(P'_i)$. 
	Otherwise, \ref{item:laminar_odd_tree_house_4} tells us that $f_c\cap h=\emptyset$, so, in particular, $r\notin f_c$. 
	As $r$ and $v_i$ are the only two vertices in $V(P'_i)\cap e_c$, $f_c$ cannot appear on $P'_i$. By injectivity of $\Phi$, $e_c\notin \Phi(E(P'_i))$ for $i\in [3]$.
\end{cpf}
Let $J\coloneqq P'_1+P'_2+P'_3\pin T$. As the unique conflict $e_c$ is not contained in $\Phi(E(J))$, we know that $(J,\Phi\upharpoonright J)\pin G$ by \cref{cor:restriction}. 
Analogously, we know that $(P'_i,\Phi\upharpoonright P'_i)\pin G$ for every $i\in [3]$ and $e_c\notin \Phi(E(P'_i))$. Finally, $e_c\cap V(P'_i)=\{r,v_i\}$. 
By \eqref{eq:non_odd_cycle_2}, we can apply \cref{lemma:parities}~\ref{lemma:parities:1} to conclude that $P'_i$ is odd for $i\in [3]$.
As $V(J)\cap e_c=\{r,v_1,v_2,v_3\}$, by \cref{lem:addition}, we can conclude that $J\oplus e_c\pin G$ is a tree house.

\subsubsection*{Case 2: $\Phi(h)\cap e_c=\emptyset$}
In this case, $f_c\ne h$, so $f_c$ is an edge appearing on one of the three paths $(P_i)_{i\in [3]}$. 
Without loss of generality, we may assume that $f_c\in E(P_1)$. Let $s$ and $t$ be the first and the last vertex of $P_1$ contained in $e_c$ when traversing it from $r$ to $\ell_1$. 
Note that $f_c\in E(P_1[s,t])$ because both endpoints of $f_c$ are contained in $e_c$. In particular, $s\ne t$. 
Further note that $r\ne s$ and $t\ne \ell_1$ because $\Phi(h)\cap e_c=\emptyset$.

As $\Phi$ is injective, we know that $e_c,\Phi(h)\notin \Phi(E(P_1[r,s])\cup E(P_1[t,\ell_1]))$. 
In particular, $(P_1[r,s]+P_1[t,\ell_1],\Phi\upharpoonright(P_1[r,s]+P_1[t,\ell_1]))\pin G$ by \cref{cor:restriction}.
Finally, we have $e_c\cap (V(P_1[r,s])\cup V(P_1[t,\ell_1]))=\{s,t\}$ and $\Phi(h)\cap (V(P_1[r,s])\cup V(P_1[t,\ell_1]))=h\cap(V(P_1[r,s])\cup V(P_1[t,\ell_1]))= \{r,\ell_1\}$ because $\Phi(h)$ is not a conflict. 
By \eqref{eq:non_odd_cycle_2}, we can apply \cref{lemma:parities}~\ref{lemma:parities:2} to derive
\begin{equation}
	|E(P_1[r,s])|+|E(P_1[t,\ell_1])|\equiv 0 \mod 2\label{eq:subpath_even}.
\end{equation}
\begin{claim}
	We have $V(P_i)\cap e_c=\emptyset$ for $i\in \{2,3\}$.
\end{claim}
\begin{cpf}
	Assume towards a contradiction that $V(P_i)\cap e_c\ne \emptyset$ for some $i\in \{2,3\}$. Let $v$ be the first vertex on $P_i$ contained in $e_c$ when traversing it from $r$ to $\ell_i$. 
	We note that $v\ne r$ because $e_c\cap \Phi(h)=\emptyset$.
	
	By injectivity of $\Phi$, we have $e_c,\Phi(h)\notin \Phi(E(P_1[r,s])\cup P_1[t,\ell_1]\cup E(P_i[r,v]))$. 
	In particular, $(P_1[r,s]+P_i[r,v],\Phi\upharpoonright (P_1[r,s]+P_i[r,v]))\pin G$ and $(P_1[t,\ell_1]+P_i[r,v],\Phi\upharpoonright (P_1[t,\ell_1]+P_i[r,v]))\pin G$ by \cref{cor:restriction}. 
	Moreover, we have $e_c\cap (V(P_1[r,s])\cup V(P_i[r,v]))=\{s,v\}$, $e_c\cap (V(P_1[t,\ell_1])\cup V(P_i[r,v]))=\{t,v\}$ and $\Phi(h)\cap (V(P_1[t,\ell_1])\cup V(P_i[r,v]))=\{\ell_1,r\}$ because $\Phi(h)$ is not a conflict. 
	By \eqref{eq:non_odd_cycle_2} and \cref{lemma:parities}, we have
	\[|E(P_1[r,s])|+|E(P_i[r,v])|\equiv 1 \mod 2 \text{ and } |E(P_1[t,\ell_1])|+|E(P_i[r,v])|\equiv 0 \mod 2.\]
	However, adding these two identities contradicts \eqref{eq:subpath_even}.
\end{cpf}

\smallskip

Let $J\coloneqq P_1[r,s]+P_1[t,\ell_1]+P_2+P_3+h\pin T$. We have $f_c\notin E(J)$, so $e_c\notin \Phi(E(J))$ by injectivity of $\Phi$. Hence, $(J,\Phi\upharpoonright J)\pin G$ by \cref{cor:restriction}. By definition of $s$ and $t$ and the previous claim, we have $s\ne t$ and $V(J)\cap e_c=\{s,t\}$. By \eqref{eq:subpath_even} and \cref{lem:addition}, $J\oplus e_c\pin G$ is a tree house.
\end{proof}
%
%%%%%%%%%%%%%%%%%%%%%%%%%%%%%%
%%%%%%%%%%%%%%%%%%%%%%%%%%%%%%
\section[A proof of Lemma~\ref{lemma:no_odd_cycle_in_H}]{A proof of \cref{lemma:no_odd_cycle_in_H}}\label{subsec:no_odd_cycle}
%%%%%%%%%%%%%%%%%%%%%%%%%%%%%%
%%%%%%%%%%%%%%%%%%%%%%%%%%%%%%

We repeat the statement for convenience. 

\LemmaNoOddCycleInH*
We prove \cref{lemma:no_odd_cycle_in_H}, deferring some technical details to later subsections.
Our goal is to find a tree house in $G^*$.
A tree house consists of three odd paths and one size-$4$ hyperedge.
Our proof starts by identifying paths $P'_1, P'_2, P'_3$, and a hyperedge $h'$ to form a candidate for a tree house in $G^*$.
To this end, recall the hypergraph $H^*$ defined in \eqref{eq:def_H} and the map $\Phi^*$ from \eqref{eq:def_phi_H}. Further recall the even cycle $C^* = G^*[U_{C^*}, F_{C^*}]$ satisfying \ref{eq:almost_nice} and \ref{eq:nice_and_almost_nice}, and the hyperedge $g^*$ defined in \eqref{eq:g_star}.

The first candidate path is the partial subhypergraph of $G^*$ induced by the vertices in $U_{C^*}$ and the hyperedges in $F_{C^*} \setminus \{g^*\}$.
The hyperedge $g^*$ induces an edge in $F_{C^*}$, say $g^* \cap U_{C^*} \eqqcolon \{r', \ell_1'\}$, so the first candidate path is the $r'$-$\ell'_1$-path $C^* - \{r', \ell'_1\}$:
\begin{equation}\label{eq:Candidate1}
P'_1 \coloneqq C^*-\{r', \ell'_1\} = G^*[U_{C^*}, F_{C^*} \setminus \{g^*\}].
\end{equation}

In order to define the remaining two paths, we use the assumption of \cref{lemma:no_odd_cycle_in_H} that $H^*$ contains an odd cycle.
Among the minimum length odd cycles in $H^*$, let $K^* \coloneqq H^*[U_{K^*},F_{K^*}]$ be one that maximizes $|E_{\ge 4}(G^*)\cap \Phi^*(F_{K^*})|$, i.e., $K^*$ is an optimal solution to the following problem:
\begin{align*}
\max &~~ |E_{\ge 4}(G^*)\cap \Phi^*(F_{K})|\\
\text{ such that} &~~ K = H^*[U_K, F_K] ~\text{is an odd cycle minimizing}~|E(K)|=|F_K|.
\end{align*}

If $(K^*,\Phi^*\upharpoonright_{K^*})$ were conflict-free, then $K^*\pin G^*$ by \cref{lemma:conflict_free_isomorphic_to_subhypergraph}, contradicting \cref{eqG*3}.
Thus, $(K^*,\Phi^*\upharpoonright_{K^*}) \qin G^*$ has a conflict.
By \cref{lem:restriction}, every conflict for $(K^*,\Phi^*\upharpoonright_{K^*})$ is also a conflict for $(H^*,\Phi^*)$.
Thus, by \cref{lem:only_conflict_g^*}, $g^*$ is the unique conflict for $(K^*,\Phi^*\upharpoonright_{K^*})$.
By the definition of a conflict,
\begin{equation}\label{eq:f1}
\text{there exists $f^* \in E(K^*)$ such that $\Phi^*\upharpoonright_{K^*}(f^*) = g^*$ and $f^* \subsetneq g^* \cap U_{K^*}$.}
\end{equation}

In \cref{subsubsec:structure_of_conflict}, we prove the following result, which states that $g^* \cap U_{K^*}$ consists of the vertices in $f^*$ and one additional vertex in $U_{C^*}$:
\begin{restatable}{lemma}{LemmaThreeVertexConflict}
We have $f^* =  (g^* \setminus U_{C^*}) \cap U_{K^*} $ and $|g^* \cap U_{C^*} \cap U_{K^*}| = 1$.
\label{lemma:3_vertex_conflict}
\end{restatable}
We also prove the following result, which describes the interaction of $K^*$ with other hyperedges.
\begin{restatable}{lemma}{LemmaInteractionWithOtherHyperedges}\label{lemma:interaction_with_other_hyperedges}
For each $e\in E_{\ge 4}(G^*)\setminus \{g^*\}$, we have $|e\cap U_{K^*}|\le 2$. Moreover, if $|e\cap U_{K^*}|=2$, then there is $f\in F_{K^*}$ with $\Phi^*(f)=e$.
\end{restatable}

Using \cref{lemma:3_vertex_conflict}, we know  $\{\ell'_2, \ell'_3\}\coloneqq f^* $ is disjoint from $g^* \cap U_{C^*} = \{r', \ell'_1\}$.
After possibly relabeling $r'$ and $\ell'_1$, we assume that 
\begin{equation}\label{eq:PathOverlap}
g^* \cap U_{C^*} \cap U_{K^*} = \{r'\}.
\end{equation}
We now define our second and third candidate paths:
\begin{align}
P'_2 &\coloneqq  (K^*-\{\ell'_2, \ell'_3\})[r', \ell'_2]\phantom{,}\label{eq:Candidate2}\\
P'_3 &\coloneqq (K^*-\{\ell'_2, \ell'_3\})[r', \ell'_3],\label{eq:Candidate3}
\intertext{and our candidate hyperedge}
h' &\coloneqq \{r', \ell'_1, \ell'_2, \ell'_3\}.\label{eq:Candidateh}
\intertext{When our candidates are added together, they yield the candidate hypergraph}
W' &\coloneqq P'_1+P'_2+P'_3+h'.\label{eq:CandidateW}
\end{align}
\begin{figure}
\footnotesize
\begin{center}
\begin{tabular}{c@{\hskip .5 cm}c@{\hskip .5 cm}c}
\begin{tikzpicture}[scale = .5, mynode/.style={circle, fill, minimum size = .1 cm, inner sep=0pt}, baseline = 0, every node/.style={scale=1}]
		\draw[draw = white](-1, -2.3) rectangle (5.75,2);
		\node[circle, draw = red, fill = red, minimum size = .1 cm, inner sep=0pt, label={[left = 0 pt, color = red] $r'$}] (r) at (-.5,0){};
		\node[circle, draw = red, fill = red, minimum size = .1 cm, inner sep=0pt, label={[right = 0 pt, color = red] $\ell'_1$}] (l1) at (5,.75){};

		\node[circle, draw = black, fill = black, minimum size = .1 cm, inner sep=0pt] (v3) at (.75,0){};
		\node[circle, draw = black, fill = black, minimum size = .1 cm, inner sep=0pt] (v7) at (2,.75){};
		\node[circle, draw = black, fill = black, minimum size = .1 cm, inner sep=0pt] (v8) at (3.25,.75){};
		\node[circle, draw = black, fill = black, minimum size = .1 cm, inner sep=0pt] (v9) at (4.5,.75){};
		\node[circle, draw = black, fill = black, minimum size = .1 cm, inner sep=0pt] (v10) at (.75,.75){};
		\node[circle, draw = black, fill = black, minimum size = .1 cm, inner sep=0pt] (v11) at (0,.75){};		
		\draw[](r) to (v11) to (v10) to (v3) to (v7) to (v8) to node[below, pos = .5]{$P'_1$}(v9) to (l1);
		\draw[red](r) to [out = 85, in = 155] (l1);
\end{tikzpicture}
&
\begin{tikzpicture}[scale = .5, mynode/.style={circle, fill, minimum size = .1 cm, inner sep=0pt}, baseline = 0, every node/.style={scale=1}]
		\draw[draw = white](-1, -2.3) rectangle (5.75,2);
		\node[circle, draw = red, fill = red, minimum size = .1 cm, inner sep=0pt, label={[left = 0 pt, color = red] $r'$}] (r) at (-.5,0){};
		\node[circle, draw = red, fill = red, minimum size = .1 cm, inner sep=0pt, label={[right = 0 pt, color = red] $\ell'_2$}] (l2) at (5,-1.5){};
		\node[circle, draw = red, fill = red, minimum size = .1 cm, inner sep=0pt, label={[below = 1 pt, color = red] $\ell'_3$}] (l3) at (3,-1.5){};

		\node[circle, draw = black, fill = black, minimum size = .1 cm, inner sep=0pt] (v1) at (.75,-1){};
		\node[circle, draw = black, fill = black, minimum size = .1 cm, inner sep=0pt] (v2) at (2,-1.5){};
		\node[circle, draw = black, fill = black, minimum size = .1 cm, inner sep=0pt] (v3) at (.75,0){};
		\node[circle, draw = black, fill = black, minimum size = .1 cm, inner sep=0pt] (v4) at (2,-.25){};
		\node[circle, draw = black, fill = black, minimum size = .1 cm, inner sep=0pt] (v5) at (3.25,-.35){};
		\node[circle, draw = black, fill = black, minimum size = .1 cm, inner sep=0pt] (v6) at (4.5,-.75){};
		
		\draw[red](l2) to (l3);
		\draw[](r) to (v1) to node[below = 2pt, pos = .5]{$P'_3$} (v2) to (l3);
		\draw[](r) to (v3) to (v4) to node[above, pos = .5]{$P'_2$} (v5) to (v6) to (l2);
\end{tikzpicture}
&
\begin{tikzpicture}[scale = .5, mynode/.style={circle, fill, minimum size = .1 cm, inner sep=0pt}, baseline = 0, every node/.style={scale=1}]
		\draw[draw = white](-1, -2.3) rectangle (5.75,2);
		\node[circle, draw = red, fill = red, minimum size = .1 cm, inner sep=0pt, label={[left = 0 pt, color = red] $r'$}] (r) at (-.5,0){};
		\node[circle, draw = red, fill = red, minimum size = .1 cm, inner sep=0pt, label={[right = 0 pt, color = red] $\ell'_1$}] (l1) at (5,.75){};
		\node[circle, draw = red, fill = red, minimum size = .1 cm, inner sep=0pt, label={[right = 0 pt, color = red] $\ell'_2$}] (l2) at (5,-1.5){};
		\node[circle, draw = red, fill = red, minimum size = .1 cm, inner sep=0pt, label={[below = 1 pt, color = red] $\ell'_3$}] (l3) at (3,-1.5){};

		\node[circle, draw = black, fill = black, minimum size = .1 cm, inner sep=0pt] (v1) at (.75,-1){};
		\node[circle, draw = black, fill = black, minimum size = .1 cm, inner sep=0pt] (v2) at (2,-1.5){};
		\node[circle, draw = black, fill = black, minimum size = .1 cm, inner sep=0pt] (v3) at (.75,0){};
		\node[circle, draw = black, fill = black, minimum size = .1 cm, inner sep=0pt] (v4) at (2,-.25){};
		\node[circle, draw = black, fill = black, minimum size = .1 cm, inner sep=0pt] (v5) at (3.25,-.35){};
		\node[circle, draw = black, fill = black, minimum size = .1 cm, inner sep=0pt] (v6) at (4.5,-.75){};
		\node[circle, draw = black, fill = black, minimum size = .1 cm, inner sep=0pt] (v7) at (2,.75){};
		\node[circle, draw = black, fill = black, minimum size = .1 cm, inner sep=0pt] (v8) at (3.25,.75){};
		\node[circle, draw = black, fill = black, minimum size = .1 cm, inner sep=0pt] (v9) at (4.5,.75){};
		\node[circle, draw = black, fill = black, minimum size = .1 cm, inner sep=0pt] (v10) at (.75,.75){};
		\node[circle, draw = black, fill = black, minimum size = .1 cm, inner sep=0pt] (v11) at (0,.75){};
		\draw[](r) to (v1) to  (v2) to (l3);
		\draw[](r) to (v3) to (v4) to (v5) to (v6) to (l2);
		\draw[](r) to (v11) to (v10) to (v3) to (v7) to (v8) to (v9) to (l1);
\end{tikzpicture}\\[1 cm]
$C^*$ & $K^*$ & $W'$
\end{tabular}
\end{center}
\caption{The left figure illustrates the cycle $C^* \pin G^*$.
The red vertices are in $g^*$.
The red line denotes that $g^*$ connects $r'$ and $\ell'_1$ in $C^*$.
The middle figure illustrates the cycle $K^* \pin H^*$. 
The right figure illustrates $P'_1+P'_2+P'_3+h'$.
For simplicity, we do not draw $h'$ in the right figure.}\label{figure:odd_cycle_K}
\end{figure}

The candidate hypergraph $W'$ is defined by three paths that start and end at the same proper hyperedge; see \cref{figure:odd_cycle_K}.
In order for $W'$ to be a tree house, the three paths must be odd and not share any vertex other than $r'$.
Next, we show that to meet these conditions, it suffices that $V_{\internal}(P'_i) \cap h' = \emptyset$ for each $i \in [3]$ and that $W' \pin G^*$.
We prove \cref{lemma:existence_CF_WOTH} in \cref{subsec:constructing_the_odd_tree_house}.

\begin{restatable}{lemma}{LemmaExistenceCFWOTH}\label{lemma:existence_CF_WOTH}
Let $W \coloneqq P_1 + P_2 + P_3 + h$ be a hypergraph, where $P_1$, $P_2$, and $P_3$ are paths and $h \coloneqq \{r, \ell_1, \ell_2, \ell_3\}$ has size $4$.
Suppose $W$ satisfies the following property:

\smallskip

\begin{enumerate}[{\it (R1)}, leftmargin = *]
%\item\label{property:Paths1} $P_1$, $P_2$, and $P_3$ are edge-disjoint paths.
%
\item\label{property:Paths1} For each $i \in [3]$, the path $P_i$ is an $r$-$\ell_i$-path with $V_{\internal}(P_i) \cap h = \emptyset$.
\end{enumerate}

\smallskip

\noindent 
If $W \pin G^*$, then $W$ is a tree house.
\end{restatable}

It may be the case that $W'\not\pin G^*$.
For instance, there may be some hyperedge $e \in E_{\ge 4}(G^*)$ that forms an edge $f \coloneqq e \cap V(P'_1)$ in $P'_1$ and that contains vertices of $P'_2$ that are not on $P'_1$.
However, we can choose a function $\Phi'$ such that $(W', \Phi') \qin G^*$ and conflicts exhibit the following structure: each conflict cannot intersect both $P'_2$ and $P'_3$, and, moreover, it intersects each path at an edge, at a vertex, or not at all.

\begin{restatable}{lemma}{LemmaExistenceTStar}\label{lemma:existence_T_star}
There exists a hypergraph $W \coloneqq P_1 + P_2 + P_3 + h$, where $P_1$, $P_2$, and $P_3$ are paths and $h \coloneqq \{r, \ell_1, \ell_2, \ell_3\}$ has size $4$, and a function $\Phi : E(W) \to E(G^*)$ such that $W$ satisfies {\it\ref{property:Paths1}} and $(W, \Phi)$ satisfies the following properties:

\smallskip
\begin{enumerate}[{\it (R1)}, leftmargin = *]
\setcounter{enumi}{1}
\item\label{property:Paths2} $(W, \Phi) \qin G^*$.
\item\label{property:Paths3} Each conflict $e\in E_{\geq 4}(G^*)$ for $(W, \Phi)$ satisfies $e\cap V(P_2)=\emptyset$ or $e\cap V(P_3)=\emptyset$.
\item\label{property:Paths4} For each $i\in[3]$ and conflict $e\in E_{\geq 4}(G^*)$ for $(W, \Phi) $, we have $|e\cap V(P_i)|\leq 2$.
Moreover, if $|e\cap V(P_i)|= 2$, then the two vertices in $e \cap V(P_i)$ form an edge $f \in E(P_i)$ such that $\Phi(f)=e$.
\end{enumerate}
\end{restatable}

In \cref{subsec:constructing_the_odd_tree_house}, we use $W'$ to prove \cref{lemma:existence_T_star}.
However, as it may be the case that $W' \not\pin G^*$, we choose another hypergraph satisfying the properties in the lemma.
More precisely, among the pairs $(W, \Phi)$ that satisfy the properties in \cref{lemma:existence_T_star}, let $(\overline{W},\overline{\Phi})$ be one minimizing $|E(\overline{W})|$.
We use {\it \ref{property:Paths1}}-{\it \ref{property:Paths4}} to prove \cref{lem:MinWOTH}; see \cref{subsec:minWOTH}.

\begin{restatable}{lemma}{LemMinWOTH}\label{lem:MinWOTH}
We have $(\overline{W}, \overline{\Phi}) \pin G^*$.
\end{restatable}

Our proofs of \cref{lemma:3_vertex_conflict,lemma:interaction_with_other_hyperedges,lemma:existence_CF_WOTH,lemma:existence_T_star,lem:MinWOTH} use arguments involving the parities of the paths defining $W'$ and $\overline{W}$.
In our proof of \cref{lem:MinWOTH}, we additionally introduce a `crossover' operation between two paths defining $\overline{W}$.
Special attention is needed to define this operation so that it preserves Properties \ref{property:Paths1}-\ref{property:Paths4}.
Hence, the proof of \cref{lem:MinWOTH} requires the most detail among the proofs of \cref{lemma:3_vertex_conflict,lemma:interaction_with_other_hyperedges,lemma:existence_CF_WOTH,lemma:existence_T_star,lem:MinWOTH}.
%.

By \cref{lemma:existence_CF_WOTH,lem:MinWOTH}, $\overline{W}$ is a tree house in $G^*$.
This completes the proof of \cref{lemma:no_odd_cycle_in_H}.
In the remainder of the section, we prove \cref{lemma:3_vertex_conflict,lemma:interaction_with_other_hyperedges,lemma:existence_CF_WOTH,lemma:existence_T_star,lem:MinWOTH}.

%%%%%%%%%%%%%%%%%%%%%%%%%%%%%%
%%%%%%%%%%%%%%%%%%%%%%%%%%%%%%
\subsection[A proof of Lemma~\ref{lemma:3_vertex_conflict}]{A proof of \cref{lemma:3_vertex_conflict}}\label{subsubsec:structure_of_conflict}
%%%%%%%%%%%%%%%%%%%%%%%%%%%%%%
%%%%%%%%%%%%%%%%%%%%%%%%%%%%%%
We repeat the statement for convenience. 

 \LemmaThreeVertexConflict*
As $g^* \in F_{C^*} \cap E_{\ge 4}(G^*)$, Equation \eqref{eq:phi_H_preimage} implies that $(\Phi^*)^{-1}(g^*) = \{g^* \setminus U_{C^*}\}$.
Thus, $f^* = (g^* \setminus U_{C^*}) \cap U_{K^*} $.
It remains to prove that $|g^* \cap U_{C^*} \cap U_{K^*}| = 1$.
We have
\[
 |g^* \cap U_{C^*} \cap U_{K^*}| =  |g^*\cap U_{K^*}| - |(g^* \setminus U_{C^*}) \cap U_{K^*}| =  |g^*\cap U_{K^*}| - |f^*|.
\]
As $g^*$ is a conflict for $(K^*,\Phi^*\upharpoonright_{K^*})$, \cref{prop:conflict_superset} implies $ |g^*\cap U_{K^*}| - |f^*| \ge 1$.
Furthermore, as $C^*$ is a cycle and $g^* \in F_{C^*}$, we have $ |g^* \cap U_{C^*} \cap U_{K^*}|  \le |g^*\cap U_{C^*}| = 2$.
Hence, $ |g^* \cap U_{C^*} \cap U_{K^*}| \in \{1,2\}$.

Assume to the contrary that $ |g^* \cap U_{C^*} \cap U_{K^*}| = 2$.
Thus, $|g^*\cap U_{K^*}|= 4$.
Denote the vertices in $g^*\cap U_{K^*}$ by $v_0,v_1,v_2$, and $v_3$ in such a way that $\{v_2,v_3\} = f^*$ and $v_0$, $v_1$, $v_2$, and $v_3$ appear in this order along $K^*$. 
The path $P \coloneqq K^* - \{v_2, v_3\}$ is even, so one of the subpaths $P[v_2,v_1]$, $P[v_1,v_0]$, or $P[v_0,v_3]$ in $K^*$ must be even; denote this subpath by $S$ and its endpoints by $x$ and $y$.
We have removed $f^*$ from $K^*$ to create $S$, so $g^*$ is not in the image of $\Phi^*\upharpoonright_{S}$ because $\Phi^*$ is injective.
As $g^*$ is the unique conflict in $K^*$, we have $(S,\Phi^*\upharpoonright_{S}) \pin G^*$ by \cref{cor:restriction}.
We apply \cref{lemma:parities}~{\it\ref{lemma:parities:1}} with  $(G,H, a,b,e) = (G^*, S, x,y, g^*)$ to conclude $S$ is odd, contradicting that $S$ is even.

\subsection[A proof of Lemma~\ref{lemma:interaction_with_other_hyperedges}]{A proof of \cref{lemma:interaction_with_other_hyperedges}}\label{subsec:interaction_with_other_hyperedges}
We repeat the statement for convenience. 

\LemmaInteractionWithOtherHyperedges*
 Let $e\in E_{\ge 4}(G^*)\setminus \{g^*\}$. By \cref{lem:only_conflict_g^*}, we have $e\cap V(H^*)\in E(H^*)$. If $|e\cap U_{K^*}|\le 1$, there is nothing to show. Next, assume that $|e\cap U_{K^*}|\ge 2$. 
 By \cref{lemma:3_vertex_conflict}, let $f^*=\{v,w\}$ and let $U_{K^*}\cap g^*=\{v,w,x\}$. 
 Let $P\coloneqq K^*-f^*$. 
 As $\Phi^*$ is injective and $g^*=\Phi^*(f^*)$ is the unique conflict with respect to $(H^*,\Phi^*)$, \cref{cor:restriction} tells us that $P\pin G^*$.
 In particular, \cref{lemma:parities}~\ref{lemma:parities:1} implies
 \begin{equation}
 |E(P[x,v])|\equiv |E(P[x,w])|\equiv 1\mod 2. \label{eq:both_subpaths_odd}
 \end{equation} 
 We further observe that
 \begin{equation}
\text{any subpath $Q$ of $P$ with endpoints in $e$ and $V_{in}(Q)\cap e =\emptyset$ is odd.}\label{eq:every_subpath_odd}	
 \end{equation}
This is because either $|E(Q)|=1$, or no edge of $Q$ has both endpoints in $e$, so $e\notin \Phi^*(E(Q))$ and \cref{lemma:parities}~\ref{lemma:parities:1} applies.

 In particular, if $|e\cap U_{K^*}|=2$, then we know that $P[u_1,u_2]$ has odd length, where $\{u_1,u_2\}\coloneqq e\cap U_{K^*}$. Let $Q$ be the $u_1$-$u_2$-path in $K^*$ containing $f^*$. Then $Q=H^*[U_Q,F_Q]\pin H^*$ is even and intersects $e$ exactly in its endpoints. In particular, $e\cap V(H^*)\notin F_Q$. Hence, $K^{**}\coloneqq H^*[U_{K^{**}},F_{K^{**}}]\coloneqq  Q\oplus (e\cap V(H^*))\pin H^*$ is an odd cycle by \cref{lem:addition}. This odd cycle is shorter than $K^*$ unless $\{u_1,u_2\}=f$ forms an edge of $K^*$.
  As $e$ is the unique proper hyperedge of $G^*$ containing both $u_1$ and $u_2$, $|\Phi^*(F_{K^{**}})\cap E_{\ge 4}(G^*)|$ is larger than $|\Phi^*(F_{K^*})\cap E_{\ge 4}(G^*)|$ unless $e\cap V(H^*)\in F_{K^*}$, i.e., the edge $\{u_1,u_2\}\in E(K^*)$ is induced by $e\cap V(H^*)$ (rather than by a proper edge). This implies $\Phi^*(f)=e$.
 
 It remains to handle the case where $|e\cap U_{K^*}|\ge 3$. Let $P_1\coloneqq P[x,v]$ and let $P_2\coloneqq P[x,w]$. 
 Then at least one of $P_1$ or $P_2$ contains at least two vertices from $e$, assume without loss of generality that this is $P_1$. 
 Let $u_1$ and $u_2$ be the vertices in $e\cap V(P_1)$ closest to $x$ and $v$, respectively. 
 Then $\{u_1,u_2\}\cap \{x,y\}=\emptyset$ because $G^*$ is disjoint. Note that neither $g^*$ nor $e$ can appear in $\Phi^*(E(P[u_1,x])\cup E(P[u_2,v]))$ because each of the paths has at most one vertex in the respective hyperedge. By \cref{lemma:parities}~\ref{lemma:parities:2}, we have
 \begin{equation}
 |E(P[u_1,x])|+|E(P[u_2,v])|\equiv 0 \mod 2. \label{eq:parities_K_star_1}
 \end{equation}
By \eqref{eq:both_subpaths_odd} and \eqref{eq:parities_K_star_1}, $P[u_1,u_2]$ must be odd. If $e\cap V(P_2)=\emptyset$, then we can proceed as in the case $|e\cap U_{K^*}|= 2$. Hence, we may assume that $e\cap V(P_2)\ne \emptyset$. Let $u_3\in e\cap V(P_2)$ be closest to $x$. By \eqref{eq:every_subpath_odd}, we have
 \begin{equation}
 |E(P[u_1,x])|+|E(P[x,u_3])|=|E(P[u_1,u_3])|\equiv 1\mod 2. \label{eq:parities_K_star_2}
\end{equation}
 Again, as no edge of $P[x,u_3]$ and $P[u_2,v]$ has both endpoints in $e$ or $g^*$, respectively, neither $g^*$ nor $e$ can appear in $\Phi^*(E(P[x,u_3])\cup E(P[u_2,v]))$. Hence, \cref{lemma:parities}~\ref{lemma:parities:2} tells us that \begin{equation}
  |E(P[u_2,v])|+|E(P[x,u_3])|\equiv 0 \mod 2. \label{eq:parities_K_star_3}
 \end{equation} Adding \eqref{eq:parities_K_star_1}, \eqref{eq:parities_K_star_2}, and \eqref{eq:parities_K_star_3} results in a contradiction as the left-hand side is even, but the right-hand side is odd.

%%%%%%%%%%%%%%%%%%%%%%%%%%%%%%
%%%%%%%%%%%%%%%%%%%%%%%%%%%%%%
\subsection[A proof of Lemma~\ref{lemma:existence_CF_WOTH}]{A proof of \cref{lemma:existence_CF_WOTH}}\label{subsec:constructing_the_odd_tree_house}
%%%%%%%%%%%%%%%%%%%%%%%%%%%%%%
%%%%%%%%%%%%%%%%%%%%%%%%%%%%%%
%
We repeat the statement for convenience. 

\LemmaExistenceCFWOTH*
As $W \pin G^*$, there is an injection $\Phi$ such that $(W, \Phi) \qin G^*$ is conflict-free by \cref{lemma:conflict_free_isomorphic_to_subhypergraph}.
Moreover, $W$ cannot contain an odd cycle by transitivity of $\pin$ and \eqref{eqG*3}.
Next, we prove that $W$ is a tree house.

Let $i \in [3]$.
The fact that $P_i$ is an $r$-$\ell_i$-path follows from Property \ref{property:Paths1}.
As $h \notin E(P_i)$, the injectivity of $\Phi$ implies that $\Phi(h) \notin \Phi(E(P_i))$.
Moreover, $V(P_i) \cap \Phi(h) = \{r, \ell_i\}$ by Property \ref{property:Paths1} and because $(W, \Phi)$ is conflict-free.
\cref{lemma:parities} \ref{lemma:parities:1} with $(G,H, a,b,e) = (G^*, P_i, r, \ell_i, \Phi(h))$ implies that $P_i$ is odd.

It is left to show that the sets $V({P}_i) \setminus\{{r}\}$ and $V({P}_j) \setminus\{{r}\}$ are disjoint for distinct $i,j \in [3]$.
Assume to the contrary that this is false.
Without loss of generality, assume there exists a vertex $v\in V({P}_1)\cap V({P}_2)\setminus \{{r}\}$.

The hypergraph $ H'  \coloneqq P_1[r,v] + P_2[r,v]$ is a closed walk in ${W}$.
As $W$ does not contain an odd cycle, $|E(H')| = |E(P_1[r,v])| + |E(P_2[r,v])|$ is even.

The hypergraph $H'' \coloneqq P_1[v, \ell_1] + P_2[v,\ell_2]$ is an $\ell_1$-$\ell_2$-walk in $W$.
Also, $\Phi(h)\notin \Phi(E(H''))$ because $\Phi(h) \notin \Phi(E(P_i))$ for each $i \in [1,2]$.
Moreover, $\Phi(h) \cap V(H'') = h \cap V(H'') = \{\ell_1, \ell_2\}$ because $(W, \Phi) $ is conflict-free and by Property \ref{property:Paths1}.
Hence, \cref{lemma:parities} \ref{lemma:parities:1} with $(G,H, a,b,e) = (G^*, H'', \ell_1, \ell_2, \Phi(h))$ allows us to conclude that $|E(H'')| = |E(P_1[v, \ell_1])| + |E(P_2[v, \ell_2])|$ is odd.

In summary, $ |E({P}_1)|+|E({P}_2)| = |E(H')|+|E(H'')|$ is odd, contradicting that $P_1$ and $P_2$ are both odd.
%

%%%%%%%%%%%%%%%%%%%%%%%%%%%%%%
%%%%%%%%%%%%%%%%%%%%%%%%%%%%%%
\subsection[A proof of Lemma~\ref{lemma:existence_T_star}]{A proof of \cref{lemma:existence_T_star}}\label{subsubsec:conflict_free}
%%%%%%%%%%%%%%%%%%%%%%%%%%%%%%
%%%%%%%%%%%%%%%%%%%%%%%%%%%%%%
%
We repeat the statement for convenience. 

\LemmaExistenceTStar*

We show that $W'$ defined in \eqref{eq:CandidateW}, along with $P'_1, P'_2$, $P'_3$, and $h' = \{\ell_1', \ell_2', \ell_3', r'\}$, satisfies Property~\ref{property:Paths1}.
Afterwards, we construct a function $\Phi'$ such that $(W', \Phi')$ satisfies Properties \ref{property:Paths2}-\ref{property:Paths4}.

\begin{claim}\label{lemma:property_weak_odd_tree_house}
The hypergraph $W'$ satisfies Property~\ref{property:Paths1}. 
\end{claim}
\begin{cpf}
For each $i \in [3]$, the path $P'_i$ is an $r'$-$\ell'_i$-path by \eqref{eq:Candidate1}, \eqref{eq:Candidate2}, and \eqref{eq:Candidate3}.
It remains to check $V_{\internal}(P'_i) \cap h' = \emptyset$ for each $i \in [3]$.
The vertices $r'$ and $\ell'_i$ do not appear as internal vertices of $P'_i$ by the definition of a path.
The vertices $\ell'_2$ and $\ell'_3$ do not appear as internal vertices of $P'_2$ and $P'_3$ because $P'_2+P'_3$ forms the path $P=K^*-\{\ell'_2,\ell'_3\}$.
Moreover, $\ell'_1$ does not appear as an internal vertex of $P'_2$ or $P'_3$ by \eqref{eq:PathOverlap}.
As $g^*\cap U_{C^*}=\{r',\ell'_1\}$ and $P'_1=C^*-\{r',\ell'_1\}$, neither $\ell'_2$ nor $\ell'_3$ appears as an internal vertex of $P'_1$.
\end{cpf}
\smallskip

Using $\Phi^*$ from \eqref{eq:def_phi_H}, we define the function $ \Phi'\colon E(W') \rightarrow E(G^*)$ by
\begin{equation}\label{eq:def_phi_T*}
 \Phi'(f)\coloneqq
\left\{\begin{array}{ll}
e &\text{if}~f=e\cap U_{C^*}\in E(P'_1)~\text{for some}~ e\in F_{C^*} \setminus \{g^*\}\\
\Phi^*\upharpoonright_{K^*}(f)&\text{if}~ f \in E(P'_2)\cup E(P'_3)\\
g^*& \text{if}~ f=h'
\end{array}\right..
\end{equation}
We have $P'_1\pin C^*\pin G^*$ and $\Phi'\upharpoonright P'_1$ equals the corresponding inclusion map. Moreover, we have $P'_2+P'_3\pin K^*$ and $\Phi'\upharpoonright (P'_2+P'_3)=\Phi^*\upharpoonright (P'_2+P'_3)$. This implies
\begin{equation}
(P'_1,\Phi'\upharpoonright P'_1)\pin G^* \text{ and } (P'_2+P'_3,\Phi'\upharpoonright (P'_2+P'_3))\qin G^*. \label{eq:individual_paths_quasi}
	\end{equation}
The following claim is used to prove Properties \ref{property:Paths2}-\ref{property:Paths4} for $(W', \Phi')$. 

\begin{claim}\label{lemma:h_no_conflict} 
We have $(\Phi')^{-1}(g^*)=\{h'\}=\{g^* \cap V(W')\}$.
\end{claim}
\begin{cpf}
By \eqref{eq:CandidateW}, we have $E(W') = E(P'_1) \cupdot E(P'_2) \cupdot E(P'_3) \cupdot \{h'\}$.
Each $f \in E(P'_1)$ is of the form $f = e \cap U_{C^*}$ for some $e \in F_{C^*} \setminus \{g^*\}$ by the definition of $P'_1$; see \eqref{eq:Candidate1}.
Hence, by~\eqref{eq:def_phi_T*}, $\Phi'(f) \neq g^*$ for $f \in E(P'_1)$.
Furthermore, $\Phi^*\upharpoonright_{K^*}$ is injective because $\Phi^*$ is injective by \eqref{eq:phi_H_preimage}.
Thus, $\Phi^*\upharpoonright_{K^*}(f) \neq g^*$ for $f \in E(P'_2)\cup E(P'_3)$ because $\Phi^*\upharpoonright_{K^*}(f^*) = g^*$ and $f^* \notin E(P'_2)\cup E(P'_3)$, where $f^*$ is from \eqref{eq:f1}.
Hence, $(\Phi')^{-1}(g^*)=\{h'\}$.

We have $g^* \cap V(W') = g^* \cap (U_{C^*} \cup U_{K^*}) = (g^* \cap U_{C^*}) \cup ((g^* \setminus U_{C^*}) \cap U_{K^*})$.
Note $g^* \cap U_{C^*} = \{r', \ell'_1\}$ by the definitions of $r'$ and $\ell'_1$.
Also, $(g^* \setminus U_{C^*}) \cap U_{K^*} = \{\ell'_2, \ell'_3\}$ by \cref{lemma:3_vertex_conflict}. 
Thus, $g^* \cap V(W') = \{r', \ell'_1, \ell'_2, \ell'_3\} = h'$.
\end{cpf}

\smallskip

We now prove that $(W', \Phi')$ satisfies Properties \ref{property:Paths2}-\ref{property:Paths4}.

\begin{claim}
The pair $(W', \Phi')$ satisfies Property \ref{property:Paths2}. \label{lem:T_star_quasi_subhypergraph}
\end{claim}
\begin{cpf}
By \cref{lemma:h_no_conflict} and \eqref{eq:individual_paths_quasi}, Property {\it\ref{def:quasi_1}} holds for $(W',\Phi')$.

Next, we establish Property {\it\ref{def:quasi_2}}.
According to \cref{lemma:h_no_conflict}, $h'$ is the only hyperedge that maps to $g^*$ under $\Phi'$.
Thus, Property {\it\ref{def:quasi_2}} holds for $g^*$.
It remains to check Property {\it\ref{def:quasi_2}} for $e\in E(G^*)\setminus \{g^*\}$.
Assume that there are two distinct edges $e_1,e_2\in E(W')\setminus \{h'\}$ with $\Phi'(e_1)=\Phi'(e_2)=e$ and $e_1\cap e_2\ne \emptyset$. 
By \eqref{eq:individual_paths_quasi} and because $\Phi'\upharpoonright (P'_2+P'_3)=\Phi^*\upharpoonright (P'_2+P'_3)$ is injective, exactly one of $e_1$ and $e_2$, say $e_1$, is contained in $E(P'_1)$, and the other, say $e_2$, is contained in $E(P'_2)\cup E(P'_3)$. 

As $e_1 \in E(P'_1)$, we have $e\in F_{C^*}$ and $e_1=e\cap U_{C^*}$. 
As $e_2\in E(P'_2)\cup E(P'_3)$, we apply \eqref{eq:def_phi_T*} to conclude that $e = \Phi'(e_2) = \Phi^*\upharpoonright_{K^*}(e_2) = \Phi^*(g ) $ for some $g \in F_{K^*} \subseteq E(H^*)$ such that $e_2 = g \cap V(K^*)$.
As $e\in F_{C^*}$ and $g \in (\Phi^*)^{-1}(e)$, it follows from \eqref{eq:phi_H_preimage} that $g = e \setminus U_{C^*}$.
In conclusion, we have $e_1=e\cap U_{C^*}$ and $e_2\subseteq e\setminus U_{C^*}$, which establishes that $e_1\cap e_2=\emptyset$.
\end{cpf}

\begin{claim}\label{lemma:intersect_only_one}
	The pair $(W', \Phi')$ satisfies Property \ref{property:Paths3}.  
\end{claim}
\begin{cpf}
	Let $e\in E_{\geq 4}(G^*)$ be a conflict for $(W', \Phi')$.
	\cref{lemma:h_no_conflict} implies that $g^*$ is not a conflict, so $e \neq g^* = \Phi'(h')$.
	By \cref{lemma:interaction_with_other_hyperedges}, we know that $e$ intersects $K^*$ in at most two vertices, and if it intersects $K^*$ in two vertices, these are connected by an edge. As $P'_2$ and $P'_3$ form edge-disjoint subpaths of $K^*$ that only share the vertex $r'\in g^*$, and $g^*\cap e = \emptyset$, it follows that $e$ can intersect only one of $P'_2$ or $P'_3$.
\end{cpf}

Property \ref{property:Paths4} for $i=2,3$ follows immediately from \cref{lemma:interaction_with_other_hyperedges} and our choice of $P'_2$, $P'_3$, and $\Phi'$. Hence, it remains to establish the following claim.

\begin{claim}
The pair $(W', \Phi')$ satisfies Property \ref{property:Paths4} for $i=1$.\label{lem:structure_of_P1}
\end{claim}
%

% %
\begin{cpf}
Let $e \in E_{\geq4}(G^*)$ be a conflict for $(W', \Phi')$.
If $e\in E_{\geq 4}(G^*)\setminus F_{C^*}$, then $|e\cap V(P'_1)|=|e\cap U_{C^*}|\leq 1$ by Property {\it\ref{eq:nice_and_almost_nice}} of $C^*$.
Thus, $e$ satisfies \ref{property:Paths4}.
If $e\in E_{\geq 4}(G^*) \cap F_{C^*}$, then $e \neq g^*$ because \cref{lemma:h_no_conflict} implies that $g^*$ is not a conflict.
We have $P'_1=G^*[U_{C^*},F_{C^*}\setminus \{g^*\}]$ by \cref{eq:Candidate1}.
As $P'_1$ is a path and $e \in F_{C^*} \setminus\{g^*\}$, we have $|e\cap V(P'_1)| = 2$.
Moreover, $e\cap V(P'_1) \in E(P'_1)$.
Hence, \cref{eq:def_phi_T*} implies that $\Phi'(e\cap V(P'_1))=e$.
Again, $e$ satisfies \ref{property:Paths4}.
\end{cpf}
%
%%%%%%%%%%%%%%%%%%%%%%%%%%%%%%
%%%%%%%%%%%%%%%%%%%%%%%%%%%%%%
\subsection[A proof of Lemma~\ref{lem:MinWOTH}]{A proof of \cref{lem:MinWOTH}}\label{subsec:minWOTH}
%%%%%%%%%%%%%%%%%%%%%%%%%%%%%%
%%%%%%%%%%%%%%%%%%%%%%%%%%%%%%
We repeat the statement for convenience. 

\LemMinWOTH*
Denote the paths defining $\overline{W}$ by $\overline{P}_1, \overline{P}_2$, $\overline{P}_3$, and the proper hyperedge by $\overline{h} = \{\overline{r}, \overline{\ell}_1, \overline{\ell}_2, \overline{\ell}_3\}$.
For $i \in [3]$, suppose $\overline{P}_i$ is given by the vertex-edge sequence
\[
v_0^i\coloneqq \overline{r},e^i_1,v^i_1,\dots,e^i_{n_i},v^i_{n_i}\coloneqq \overline{\ell}_i.
\]

 Assume to the contrary that $(\overline{W}, \overline{\Phi})\not\pin G^*$.
 By \cref{lemma:conflict_free_isomorphic_to_subhypergraph}, there exists a conflict for $(\overline{W}, \overline{\Phi})$. 
 We will use this conflict in order to obtain a shorter tree house $(W^\times,\Phi^\times)\qin G^*$ subject to \ref{property:Paths1}-\ref{property:Paths4}, contradicting our choice of $(\overline{W}, \overline{\Phi})$. The high-level strategy can be described as follows:
 Using that each path $\overline{P}_i$ is `conflict-free on its own' (by \ref{property:Paths4}), we argue that each conflict has to interact with at least two different paths. Using \cref{lemma:parities} and the fact that $G^*$ does not contain an odd cycle, we argue that there must be a pair of interleaving conflicts on two paths (see \cref{fig:cross_operation}). We use these to `cross over' the two paths in order to eliminate the conflicts and shorten the tree house.
We begin with some claims illustrating how a conflict interacts with $\overline{W}$.

 \begin{claim}\label{lem:conflict_between_P1_P2}
If $e\in E(G^*)$ is a conflict for $(\overline{W}, \overline{\Phi})$, then it satisfies the following:
 \begin{enumerate}[{\it (CF1)}, leftmargin = *]
 \item \label{prop:bar1} $e \in E_{\ge 4}(G^*) \setminus \{\overline{\Phi}(\overline{h})\}$ and $e \cap \overline{h} = \emptyset$.
 \item\label{prop:bar2}  $(\overline{P}_i, \overline{\Phi}\upharpoonright_{\overline{P}_i}) \pin G^*$  for each $i \in [3]$. 
Thus, $e$ is not a conflict for $(\overline{P}_i, \overline{\Phi}\upharpoonright_{\overline{P}_i}) $.
 \item \label{prop:bar3}$e$ intersects $V(\overline{P}_1)$ and exactly one of $V(\overline{P}_2)$ or $V(\overline{P}_3)$.
 \item \label{prop:bar4}$e\cap V(\overline{P}_i)\cap V(\overline{P}_j)=\emptyset$ for all $i,j\in [3]$ with $i\neq j$.
 \end{enumerate}
 \end{claim}
 \begin{cpf}
 Let $e \in E(G^*)$ be a conflict.
 The hyperedge $\overline{h}$ intersects both $\overline{P}_2$ and $\overline{P}_3$ by definition.
As $(\overline{W}, \overline{\Phi}) \qin G^*$ by Property \ref{property:Paths2}, Property \ref{def:quasi_1} implies that $\overline{h} \subseteq \overline{\Phi}(\overline{h})$.
Therefore, $\overline{\Phi}(\overline{h})$ is not a conflict by Property \ref{property:Paths3}.
Hence, $e \neq \overline{\Phi}(\overline{h})$.
Moreover, by the definition of a conflict, $e \supsetneq f$ for some $f  \in E(\overline{W}) = E_{\ge 2}(\overline{W})$.
 Hence, $e \in E_{\ge 4}(G^*)$. 
 The disjointness of $G^*$ implies that $e \cap \overline{h} \subseteq e \cap \overline{\Phi}(\overline{h}) = \emptyset$.
 This proves \ref{prop:bar1}.
 
  Let $i \in [3]$.
Assume to the contrary that there exists a conflict $g \in E(G^*)$ for $(\overline{P}_i, \overline{\Phi}\upharpoonright_{\overline{P}_i})$.
As $g$ is a conflict, there exists an edge $f \in E(\overline{P}_i)$ such that $\overline{\Phi}(f) = g$ and $f \subsetneq g \cap V(\overline{P}_i) $.
Hence, $|g \cap V(\overline{P}_i)| \ge 3$, and $g \in E_{\ge 4}(G^*)$.
We have that $g$ is a conflict for $(\overline{W}, \overline{\Phi})$ by \cref{lem:restriction}.
By Property \ref{property:Paths4}, it follows that $|g \cap V(\overline{P}_i)| \le 2$, which is a contradiction.
This proves \ref{prop:bar2}.

Each pair $(\overline{P}_i, \overline{\Phi}\upharpoonright_{\overline{P}_i})$ is conflict-free by \ref{prop:bar2}, so $e$ intersects two distinct paths.
By Property \ref{property:Paths3}, $e$ does not intersect both $V(\overline{P}_2)$ and $V(\overline{P}_3)$.
This proves \ref{prop:bar3}.
 
Finally, we prove \ref{prop:bar4}.
As $e$ is a conflict and $e \neq \overline{\Phi}(\overline{h})$ by \ref{prop:bar1}, there exists an edge $f \in E(\overline{W}) \setminus \{\overline{h}\}$ such that $\overline{\Phi}(f) = e$.
Furthermore, by \ref{prop:bar3}, $e$ intersects $V(\overline{P}_1)$ and exactly one of $V(\overline{P}_2)$ or $V(\overline{P}_3)$.
Without loss of generality, suppose $e$ intersects $V(\overline{P}_2)$ and is disjoint from $V(\overline{P}_3)$.
Hence, $f \in E(\overline{P}_1) \cup E(\overline{P}_2)$.
Without loss of generality, suppose $f \in  E(\overline{P}_1) $.
As $e$ is not a conflict for $(\overline{P}_1, \overline{\Phi}\upharpoonright_{\overline{P}_1}) $ by Property \ref{prop:bar2}, we know that $e \cap V(\overline{P}_1) = f$.

Assume by contradiction that there exists $u \in e \cap V(\overline{P}_1) \cap V(\overline{P}_2)$.
As $e \cap V(\overline{P}_1) = f$, we have $ u \in f$.
Moreover, as $e$ is a conflict for $(\overline{W}, \overline{\Phi})$, we know that there exists some vertex $v \in V(\overline{P}_2) \cap e \setminus f$.
In particular, $v \neq u$.
By Property \ref{property:Paths4}, there is an edge $f' \in E(\overline{P}_2)$ such that $f' = \{v, u\}$ and $\overline{\Phi}(f') = e$.
Hence, $(\overline{\Phi}{}^{-1})(e)$ contains $f$ and $f'$, which both contain $u$ yet are distinct because $v \in f' \setminus f$.
However, this contradicts Property \ref{def:quasi_2}, proving \ref{prop:bar4}.
 \end{cpf}
\smallskip

We now take a step towards defining a crossover operation that uses two conflicts that cross over between two paths in $\overline{W}$.
This crossover will be used to create $(W^{\times}, \Phi^{\times})$ satisfying~\ref{property:Paths1}\,-\,\ref{property:Paths4} with $|E(W^{\times})| < |E(\overline{W})|$, contradicting the minimality of $(\overline{W}, \overline{\Phi})$.
We start by identifying the two conflicts.

 Let $s_1 \in [n_1-1]$ be the smallest index such that $v^1_{s_1}$ is in a conflict.
 Such an index exists because every conflict is disjoint from $\overline{h}$ by Property \ref{prop:bar1}, and intersects $V(\overline{P}_1) $ by Property  \ref{prop:bar3}.
 There is a unique conflict $g_1 \in E(G^*)$ containing $v^1_{s_1}$ because $G^*$ is disjoint and every conflict is in $E_{\ge 4}(G^*)$ by Property \ref{prop:bar1}.
Moreover, $g_1$ intersects exactly one of $V(\overline{P}_2)$ or $V(\overline{P}_3)$ by Property \ref{prop:bar3}.
 Without loss of generality, suppose $g_1\cap V(\overline{P}_2)\neq \emptyset$.
 Let $t_1 \in [n_2-1]$ be the largest index such that $v^2_{t_1} \in g_1$.

 Let $s_2 \in [n_2-1]$ be the smallest index such that $v^2_{s_2}$ is in a conflict; as $v^1_{t_1} \in g_1$, we have $s_2\le t_1$.
 As is the case with $v^1_{s_1}$, there is a unique conflict $g_2\in E(G^*)$ such that $v^2_{s_2}\in g_2$.
 We have $g_2 \cap \overline{h} = \emptyset$ by Property \ref{prop:bar1}, and $g_2 \cap V(\overline{P}_1) \neq \emptyset$ by Property \ref{prop:bar3}.
 Let $t_2 \in [n_1-1]$ be the largest index such that $v^1_{t_2} \in g_2$.
See \cref{fig:cross_operation}.

\begin{figure}
\footnotesize
\centering
\begin{tabular}{@{\hskip 0 cm}c@{\hskip .675 cm}c@{\hskip .675 cm}c}
\begin{tikzpicture}[scale = .5, mynode/.style={circle, fill, minimum size = .1 cm, inner sep=0pt}, baseline = 0, every node/.style={scale=1}]
\node[circle, draw = red, fill = red, minimum size = .1 cm, inner sep=0pt, label={[left = 0 pt, color = red] $\overline{r}$}] (r) at (0,0){};
\node[circle, draw = red, fill = red, minimum size = .1 cm, inner sep=0pt, label={[right = 0 pt, color = red] $\overline{\ell}_1$}] (l1) at (6,2){};
\node[circle, draw = red, fill = red, minimum size = .1 cm, inner sep=0pt, label={[right = 0 pt, color = red] $\overline{\ell}_2$}] (l2) at (6,0){};
\node[circle, draw = red, fill = red, minimum size = .1 cm, inner sep=0pt, label={[right = 0 pt, color = red] $\overline{\ell}_3$}] (l3) at (6,-1.5){};

\node[circle, draw = black, fill = black, minimum size = .1 cm, inner sep=0pt, label={[below = 2 pt, color = black] $v^2_{x_2}$}] (x2) at (1.75,0){};

\node[rectangle, draw = green!50!black!80, fill = green!50!black!80, minimum size = .15 cm, inner sep=0pt, label={[above left = -2 pt, color = green!50!black!80] $v^1_{s_1}$}] (s1) at (.5, 1.5){};
\node[rectangle, draw = green!50!black!80, fill = green!50!black!80, minimum size = .15 cm, inner sep=0pt, label={[above right = -3 pt, color = green!50!black!80] $v^1_{s_1+1}$}] (s1p1) at (1.25, 1.5){};
\node[shape border rotate=0, regular polygon, regular polygon sides=3, draw = orange, fill = orange, minimum size = .2 cm, inner sep=0pt, label={[above right= -5 pt, color = orange] $v^1_{t_2}$}] (t2) at (1.5, .5){};

\node[shape border rotate=0, regular polygon, regular polygon sides=3, draw = orange, fill = orange, minimum size = .2 cm, inner sep=0pt, label={[below = 3 pt, color = orange] $v^2_{s_2}$}] (s2) at (2.75,0){};
\node[shape border rotate=0, regular polygon, regular polygon sides=3, draw = orange, fill = orange, minimum size = .2 cm, inner sep=0pt, label={[above = -1 pt, color = orange] $v^2_{s_2+1}$}] (s2p1) at (3.75,0){};
\node[rectangle, draw = green!50!black!80, fill = green!50!black!80, minimum size = .15 cm, inner sep=0pt, label={[below = 3 pt, color = green!50!black!80] $v^2_{t_1}$}] (t1) at (5,0){};

%P1
\draw[decorate, decoration={snake, amplitude=.3mm, post length=0mm, pre length=0mm}](r) to [out = 90, in = -90] (s1); 
\draw[green!50!black!80, thick](s1) to (s1p1); 
\draw[decorate, decoration={snake, amplitude=.3mm, post length=0mm, pre length=0mm}](s1p1) to [out = -90, in = 180] (t2); 
\draw[decorate, decoration={snake, amplitude=.3mm, post length=0mm, pre length=0mm}](t2) to [out = -90, in = 90] (x2); 
\draw[decorate, decoration={snake, amplitude=.3mm, post length=0mm, pre length=0mm}](x2) to [out = 45, in = 180] (l1);

%P2
\draw[decorate, decoration={snake, amplitude=.3mm, post length=0mm, pre length=0mm}](r) to [out = 0, in = 180] (x2);
\draw[decorate, decoration={snake, amplitude=.3mm, post length=0mm, pre length=0mm}](x2) to [out = 0, in = 180] (s2);
\draw[orange, thick](s2) to [out = 0, in = 180] (s2p1);
\draw[decorate, decoration={snake, amplitude=.3mm, post length=0mm, pre length=0mm}](s2p1) to [out = 0, in = 180] (t1);
\draw[decorate, decoration={snake, amplitude=.3mm, post length=0mm, pre length=0mm}](t1) to [out = 0, in = 180] (l2);

%P3
\draw[decorate, decoration={snake, amplitude=.3mm, post length=0mm, pre length=0mm}](r) to [out = -90, in = 180] (l3);

\end{tikzpicture}
&
\begin{tikzpicture}[scale = .5, mynode/.style={circle, fill, minimum size = .1 cm, inner sep=0pt}, baseline = 0, every node/.style={scale=1}]
\node[circle, draw = red, fill = red, minimum size = .1 cm, inner sep=0pt, label={[left = 0 pt, color = red] $\overline{r}$}] (r) at (0,0){};
\node[circle, draw = red, fill = red, minimum size = .1 cm, inner sep=0pt, label={[right = 0 pt, color = red] $\overline{\ell}_1$}] (l1) at (6,2){};
\node[circle, draw = red, fill = red, minimum size = .1 cm, inner sep=0pt, label={[right = 0 pt, color = red] $\overline{\ell}_2$}] (l2) at (6,0){};
\node[circle, draw = red, fill = red, minimum size = .1 cm, inner sep=0pt, label={[right = 0 pt, color = red] $\overline{\ell}_3$}] (l3) at (6,-1.5){};

\node[circle, draw = black, fill = black, minimum size = .1 cm, inner sep=0pt, label={[below = 2 pt, color = black] $v^2_{x_2}$}] (x2) at (1.75,0){};

\node[rectangle, draw = green!50!black!80, fill = green!50!black!80, minimum size = .15 cm, inner sep=0pt, label={[above left = -2 pt, color = green!50!black!80] $v^1_{s_1}$}] (s1) at (.5, 1.5){};
\node[shape border rotate=0, regular polygon, regular polygon sides=3, draw = orange, fill = orange, minimum size = .2 cm, inner sep=0pt, label={[above right= -5 pt, color = orange] $v^1_{t_2}$}] (t2) at (1.5, .5){};

\node[shape border rotate=0, regular polygon, regular polygon sides=3, draw = orange, fill = orange, minimum size = .2 cm, inner sep=0pt, label={[below = 3 pt, color = orange] $v^2_{s_2}$}] (s2) at (2.75,0){};
\node[rectangle, draw = green!50!black!80, fill = green!50!black!80, minimum size = .15 cm, inner sep=0pt, label={[below =3 pt, color = green!50!black!80] $v^2_{t_1}$}] (t1) at (5,0){};

%P1
\draw[decorate, decoration={snake, amplitude=.3mm, post length=0mm, pre length=0mm}](r) to [out = 0, in = 180] (x2);
\draw[decorate, decoration={snake, amplitude=.3mm, post length=0mm, pre length=0mm}](x2) to [out = 0, in = 180] (s2);
\draw[dashed, thick, draw = orange](s2) to [out = 90, in = 0] (t2);
\draw[decorate, decoration={snake, amplitude=.3mm, post length=0mm, pre length=0mm}](t2) to [out = -90, in = 90] (x2); 
\draw[decorate, decoration={snake, amplitude=.3mm, post length=0mm, pre length=0mm}](x2) to [out = 45, in = 180] (l1);

%P2
\draw[decorate, decoration={snake, amplitude=.3mm, post length=0mm, pre length=0mm}](r) to [out = 90, in = -90] (s1); 
\draw[dashed, thick, draw = green!50!black!80](s1) to [out = 45, in = -180] (t1);
\draw[decorate, decoration={snake, amplitude=.3mm, post length=0mm, pre length=0mm}](t1) to [out = 0, in = 180] (l2);

%P3
\draw[decorate, decoration={snake, amplitude=.3mm, post length=0mm, pre length=0mm}](r) to [out = -90, in = 180] (l3);

\end{tikzpicture}
&
\begin{tikzpicture}[scale = .5, mynode/.style={circle, fill, minimum size = .1 cm, inner sep=0pt}, baseline = 0, every node/.style={scale=1}]
\node[circle, draw = red, fill = red, minimum size = .1 cm, inner sep=0pt, label={[left = 0 pt, color = red] $\overline{r}$}] (r) at (0,0){};
\node[circle, draw = red, fill = red, minimum size = .1 cm, inner sep=0pt, label={[right = 0 pt, color = red] $\overline{\ell}_1$}] (l1) at (6,2){};
\node[circle, draw = red, fill = red, minimum size = .1 cm, inner sep=0pt, label={[right = 0 pt, color = red] $\overline{\ell}_2$}] (l2) at (6,0){};
\node[circle, draw = red, fill = red, minimum size = .1 cm, inner sep=0pt, label={[right = 0 pt, color = red] $\overline{\ell}_3$}] (l3) at (6,-1.5){};

\node[circle, draw = black, fill = black, minimum size = .1 cm, inner sep=0pt, label={[below = 2 pt, color = black] $v^2_{x_2}$}] (x2) at (1.75,0){};

\node[rectangle, draw = green!50!black!80, fill = green!50!black!80, minimum size = .15 cm, inner sep=0pt, label={[above left = -2 pt, color = green!50!black!80] $v^1_{s_1}$}] (s1) at (.5, 1.5){};

\node[rectangle, draw = green!50!black!80, fill = green!50!black!80, minimum size = .15 cm, inner sep=0pt, label={[below = 3 pt, color = green!50!black!80] $v^2_{t_1}$}] (t1) at (5,0){};

%P1
\draw[decorate, decoration={snake, amplitude=.3mm, post length=0mm, pre length=0mm}](r) to [out = 0, in = 180] (x2);
\draw[decorate, decoration={snake, amplitude=.3mm, post length=0mm, pre length=0mm}](x2) to [out = 45, in = 180] (l1);

%P2
\draw[decorate, decoration={snake, amplitude=.3mm, post length=0mm, pre length=0mm}](r) to [out = 90, in = -90] (s1); 
\draw[dashed, thick, draw = green!50!black!80](s1) to [out = 45, in = -180] (t1);
\draw[decorate, decoration={snake, amplitude=.3mm, post length=0mm, pre length=0mm}](t1) to [out = 0, in = 180] (l2);

%P3
\draw[decorate, decoration={snake, amplitude=.3mm, post length=0mm, pre length=0mm}](r) to [out = -90, in = 180] (l3);

\end{tikzpicture}
    \end{tabular}
    \caption{The left figure shows $\overline{W}$ before the crossover operation.
    The green vertices are in $g_1$, and the orange vertices are in $g_2$. 
    The middle figure shows the walks $\Omega_1$ and $\Omega_2$; see \eqref{eq:Omega1Defn}.
    The right figure shows $W^{\times}$, the result of the crossover operation.
    For simplicity, we do not draw the hyperedge $\overline{h}$ in the figures.}
    \label{fig:cross_operation}
\end{figure}

The conflicts $g_1$ and $g_2$ interact with the paths $\overline{P}_1, \overline{P}_2$, and $\overline{P}_3$ in particular ways. 
For instance, Property \ref{prop:bar3} implies
\begin{equation}\label{eq:g_1_2_dont_intersect_P_3}
g_i\cap V(\overline{P}_3)=\emptyset \text{ for $i\in [2]$.} 
\end{equation}
Furthermore, we have the following:
\begin{claim}\label{claim:conflict_edge_f_i}
Let $i \in [2]$.
We have 
\[
\text{$g_i \cap V(\overline{P}_i) \subseteq \{v^i_{s_i}, v^i_{s_i + 1}\}$ and $g_i \cap V(\overline{P}_{3-i}) \subseteq \{v^{3-i}_{t_i - 1}, v^{3-i}_{t_i}\}$.}
\]
Moreover, there exists $f_i \in E(\overline{P}_1) \cup E(\overline{P}_2)$ with $\overline{\Phi}(f_i) = g_i$, where
\[
f_i = \left\{
\begin{array}{ll}
 \{v^i_{s_i}, v^i_{s_i + 1}\}, & \text{if $f_i \in E(\overline{P}_i)$}\\[.1cm]
 \{v^{3-i}_{t_i - 1}, v^{3-i}_{t_i}\} , & \text{if $f_i \in E(\overline{P}_{3-i})$}
\end{array}
\right..
\]
\end{claim}
\begin{cpf}
    We prove the result only for $i = 1$ as the proof for $i = 2$ is analogous.
    By the definition of a conflict, there exists $f_1 \in E(\overline{W})$ such that $\overline{\Phi}(f_1) = g_1$.
    Property \ref{prop:bar1} implies that $f_1\neq \overline{h}$. 
    Furthermore, Property \ref{def:quasi_1} for $(\overline{W}, \overline{\Phi})$ and \eqref{eq:g_1_2_dont_intersect_P_3} imply that $f_1 \in E(\overline{P}_1) \cup E(\overline{P}_2)$.

    If $f_1 \in E(\overline{P}_1)$, then Property \ref{property:Paths4} for $(\overline{W}, \overline{\Phi})$ implies $f_1 = g_1 \cap V(\overline{P}_1)$.
Moreover, the minimality of $s_1$ implies that $f_1= \{v^1_{s_1}, v^1_{s_1+1}\}$.
Similarly, Property \ref{property:Paths4} and the maximality of $t_1$ imply that $g_1 \cap V(\overline{P}_2) \subseteq \{v^2_{t_1-1}, v^2_{t_1}\}$.

If $f_1\in E(\overline{P}_2)$, then Property \ref{property:Paths4} for $(\overline{W}, \overline{\Phi})$ implies that $f_1 = g_1 \cap V(\overline{P}_2)$.
Moreover, the maximality of $t_1$ implies that $f_1= \{v^2_{t_1 -1}, v^2_{t_1}\}$.
Similarly, Property \ref{property:Paths4} and the minimality of $s_1$ imply that $g_1 \cap V(\overline{P}_1) \subseteq \{v^1_{s_1}, v^1_{s_1+1}\}$.
\end{cpf}
\smallskip

We can now show that the conflicts $g_1$ and $g_2$ are distinct.
\begin{claim}\label{claim:g1_neq_g2}
We have $g_1\ne g_2$.
\end{claim}
\begin{cpf}
Assume by contradiction that $g_1=g_2\eqqcolon g$. 
By \cref{claim:conflict_edge_f_i}, there exists $f \in E(\overline{P}_1) \cup E(\overline{P}_2)$ such that $\overline{\Phi}(f) = g$.
Without loss of generality, assume $f \in E(\overline{P}_1)$.
\cref{claim:conflict_edge_f_i} states $f = g\cap V(\overline{P}_1)= \{v^1_{s_1}, v^1_{s_1+1}\}$ and implies that $g \cap V(\overline{P}_2[\overline{r},v^2_{s_2}])=\{v^2_{s_2}\}$.
By Property \ref{prop:bar4}, $v^2_{s_2}$ is distinct from $v^1_{s_1}$ and $v^1_{s_1+1}$.

Consider the even number
\begin{equation}\label{eq:g1notg2Parity}
2\big ( |E(\overline{P}_1[\overline{r},v^1_{s_1}])| + |E(\overline{P}_1[v^1_{s_1+1},\overline{\ell}_1])| + |E(\overline{P}_2[\overline{r}, v^2_{s_2}])| \big ).
\end{equation}
We arrive at a contradiction by showing that \eqref{eq:g1notg2Parity} is also odd.
To this end, it suffices to prove the following relations because the sum of the left hand sides is \eqref{eq:g1notg2Parity} while the sum of the right hand sides is odd:
\begin{align}
|E(\overline{P}_1[\overline{r},v^1_{s_1}])| + |E(\overline{P}_1[v^1_{s_1+1}, \overline{\ell}_1])| &\equiv 0 \mod 2 \label{eq:P1_P2_iii}\\
|E(\overline{P}_1[\overline{r},v^1_{s_1}])| + |E(\overline{P}_2[\overline{r}, v^2_{s_2}])| &\equiv 1 \mod 2
	\label{eq:P1_P2_i}\\
|E(\overline{P}_1[v^1_{s_1+1},\overline{\ell_1}])| + |E(\overline{P}_2[\overline{r},v^2_{s_2}])| &\equiv 0 \mod 2.	\label{eq:P1_P2_ii}
\end{align}

By Property \ref{property:Paths1}, we have $V(\overline{P}_1) \cap \Phi(\overline{h}) = \{\overline{r}, \overline{\ell}_1\}$.
Property \ref{prop:bar2} tell us that $(\overline{P}_1, \overline{\Phi}\upharpoonright_{P_1}) \pin G^*$.
We have $(\overline{\Phi}{}^{-1})(\overline{\Phi}(\overline{h})) = \{\overline{h}\}$ because $\overline{\Phi}(\overline{h})$ is not a conflict by Property \ref{prop:bar1}.
\cref{lemma:parities} \ref{lemma:parities:1} with $(G,H,a,b,e) = (G^*, \overline{P}_1, \overline{r}, \overline{\ell}_1,\overline{\Phi}(\overline{h}))$ implies that $\overline{P}_1$ is odd.
This proves \eqref{eq:P1_P2_iii}.

In order to prove \eqref{eq:P1_P2_i} and \eqref{eq:P1_P2_ii}, consider
\[
J_1 \coloneqq \overline{P}_1[\overline{r},v^1_{s_1}] + \overline{P}_2[\overline{r}, v^2_{s_2}] \quad \text{and} \quad
 J_2 \coloneqq \overline{P}_1[v^1_{s_1+1},\overline{\ell_1}] + \overline{P}_2[\overline{r},v^2_{s_2}].
\]
We will apply \cref{lemma:parities} to $(J_1, \overline{\Phi}\upharpoonright_{J_1})$ and $(J_2,\overline{\Phi}\upharpoonright_{J_2})$.
To this end, we show $(J_i, \overline{\Phi}\upharpoonright_{J_i})\pin G^*$ for each $i \in [2]$.

For each $i \in [2]$, note that $V(J_i)\cap g \subseteq \{v^1_{s_1}, v^1_{s_1+1}, v^2_{s_2}\}$ and each of the subpaths of $\overline{P}_1$ and $\overline{P}_2$ forming $J_i$ contains at most one of these vertices.
Thus, there is no $f\in E(J_i)$ with $f\subseteq g$.
Hence, by Property \ref{def:quasi_1} we have
\begin{equation}\label{eqgnotinbarphi}
g\notin \overline{\Phi}(E(J_i)).
\end{equation}

For each $i \in [2]$, we have $J_i \pin \overline{W}$, so $(J_i,\overline{\Phi}\upharpoonright_{J})\qin G^*$ by \cref{lem:restriction}.
Assume towards a contradiction that $e\in E(G^*)$ is a conflict for $(J_i,\overline{\Phi}\upharpoonright_{J_i})$ for some $i \in [2]$. 
We have $e \neq g$ by \eqref{eqgnotinbarphi}.
By \cref{lem:restriction}, $e$ is a conflict for $(\overline{W}, \overline{\Phi})$.
As $(\overline{P}_1,\overline{\Phi}\upharpoonright_{\overline{P}_1})$ is conflict-free, we have $e \cap V(\overline{P}_2) \neq \emptyset$.
The minimality of $s_2$ implies that $v^2_{s_2}\in e$.
Hence, $e=g$ because $G^*$ is disjoint and $e \in E_{\ge 4}(G^*)$ by Property \ref{prop:bar1}.
However, this contradicts that $e \neq g$.
Thus, $(J_i,\overline{\Phi}\upharpoonright_{J_i})\pin G^*$.

We have $g \notin  \overline{\Phi}(E(J_i))$ for each $i \in [2]$ by \eqref{eqgnotinbarphi}.
Furthermore, we have $\overline{\Phi}(\overline{h})\notin \overline{\Phi}(E(J_2))$.
Indeed, $V(J_2)\cap \overline{h} = \{\overline{r}, \overline{\ell}_1\}$ by Property \ref{property:Paths1}.
As $E(J_2)$ does not contain an edge from $\overline{r}$ to $\overline{\ell}_1$, Property \ref{def:quasi_1} implies that $\overline{\Phi}(\overline{h})\notin \overline{\Phi}(E(J_2))$.

We can apply \cref{lemma:parities}~{\it\ref{lemma:parities:1}} with $(G,H,a,b,e) = (G^*, J_1, v^1_{s_1}, v^2_{s_2}, g)$ to conclude \eqref{eq:P1_P2_i}.
We conclude \eqref{eq:P1_P2_ii} by applying \cref{lemma:parities}~{\it\ref{lemma:parities:2}} with \\
$(G,H,a,b,c,d,e,f) = (G^*, J_2, v^1_{s_1+1}, \overline{\ell}_1, v^2_{s_2}, \overline{r}, g, \overline{\Phi}(\overline{h}))$. 
\end{cpf}

Now, we describe the aforementioned crossover operation: we walk along $\overline{P}_2$ from $\overline{r}$ to $v^2_{s_2}$, then cross to $\overline{P}_1$ using the edge $\{v^2_{s_2}, v^1_{t_2}\} \subseteq g_2$, then continue to $\overline{\ell}_1$ along $\overline{P}_1$;  see~\cref{fig:cross_operation}.
Crossing over in this way creates a walk
\begin{equation}\label{eq:Omega1Defn}
\Omega_1 \coloneqq \overline{P}_2[\overline{r},v^2_{s_2}]+\{v^2_{s_2},v^1_{t_2}\}+\overline{P}_1[v^1_{t_2},\overline{\ell}_1].
\end{equation}

By Property \ref{prop:bar4}, we have $v^2_{s_2} \neq v^1_{t_2}$.
This, together with the fact that $\overline{P}_1$ and $\overline{P}_2$ are paths, implies the walk $\Omega_1$ is a path if and only if $\overline{P}_2[\overline{r},v^2_{s_2}]$ and $\overline{P}_1[v^1_{t_2},\overline{\ell}_1]$ are vertex-disjoint.
In the case that $\overline{P}_2[\overline{r},v^2_{s_2}]$ and $\overline{P}_1[v^1_{t_2},\overline{\ell}_1]$ are not vertex-disjoint, we can identify a path in $\Omega_1$ by walking from $\overline{r}$ along $\overline{P}_2$ to the first vertex shared by $\overline{P}_2[\overline{r},v^2_{s_2}]$ and $\overline{P}_1[v^1_{t_2},\overline{\ell}_1]$, and then continuing to $\overline{\ell}_1$ in $\overline{P}_1$.
Formally, our crossover operation is:
\begin{equation}\label{eq:CrossOverPath1}
P^{\times}_1\coloneqq
\left\{
\begin{array}{@{\hskip 0 cm}c@{\hskip .2 cm}l@{\hskip 0 cm}}
\Omega_1, &\text{if}~V(\overline{P}_2[\overline{r},v^2_{s_2}]) \cap V(\overline{P}_1[v^1_{t_2},\overline{\ell}_1]) = \emptyset\\[.1cm]
\overline{P}_2[\overline{r},v^2_{x_2}]+\overline{P}_1[v^2_{x_2}, \overline{\ell}_1], &
\text{otherwise, where $x_2 \in [s_2]$ is the smallest}\\
&\text{index such that $v^2_{x_2} \in V(\overline{P}_1[v^1_{t_2},\overline{\ell}_1])$.}
\end{array}
\right.
\end{equation}

The cases defining $P^{\times}_1$ ensure that it is an $\overline{r}$-$\overline{\ell}_1$-path.
Also, the edge $\{v^2_{s_2},v^1_{t_2}\} \in E(\Omega_1)$ appears in $E(P^{\times}_1)$ if and only if $\Omega_1 = P^{\times}_1$.
We can similarly define an $\overline{r}$-$\overline{\ell}_2$-walk by crossing over from $\overline{P}_1$ to $\overline{P}_2$ using the edge $\{v^1_{s_1}, v^2_{t_1}\} \subseteq g_1$:
\[
\Omega_2 \coloneqq \overline{P}_1[\overline{r},v^1_{s_1}]+\{v^1_{s_1},v^2_{t_1}\}+\overline{P}_2[v^2_{t_1},\overline{\ell}_2].
\]
We define the analogous $\overline{r}$-$\overline{\ell}_2$-path in $\Omega_2$:
\begin{equation}\label{eq:CrossOverPath2}
P^{\times}_2\coloneqq
\left\{
\begin{array}{@{\hskip 0 cm}c@{\hskip .2 cm}l@{\hskip 0 cm}}
\Omega_2, &\text{if}~V(\overline{P}_1[\overline{r},v^1_{s_1}]) \cap V(\overline{P}_2[v^2_{t_1},\overline{\ell}_2]) = \emptyset\\[.1cm]
\overline{P}_1[\overline{r},v^1_{y_1}]+\overline{P}_2[v^1_{y_1}, \overline{\ell}_2], &
\text{otherwise, where $y_1 \in [s_1]$ is the smallest}\\
&\text{index such that $v^1_{y_1} \in V(\overline{P}_2[v^2_{t_1},\overline{\ell}_2]) $.}
\end{array}
\right.
\end{equation}
The edge $\{v^1_{s_1},v^2_{t_1}\} \in E(\Omega_2)$ appears in $E(P^{\times}_2)$ if and only if $\Omega_2 = P^{\times}_2$.

Given the tuple $(P^{\times}_1, P^{\times}_2, \overline{P}_3, \overline{h})$, we define a new hypergraph
\begin{equation}\label{eqWxdefn}
W^{\times}\coloneqq P^{\times}_1 + P^{\times}_2 + \overline{P}_3 + \overline{h}
\end{equation}
and a function $\Phi^{\times}\colon E({W}^{\times})\rightarrow E(G^*)$ by
\begin{align}\label{eq:def_Phi_prime}
\Phi^{\times} (f) \coloneqq 
\left\{
\begin{array}{cl}
\overline{\Phi}\upharpoonright_{E(\overline{W})\cap E(W^{\times})}(f), &\text{if}~f \in E(\overline{W})\cap E(W^{\times})\\
g_1, &\text{if $f = \{v^1_{s_1}, v^2_{t_1}\} \in E(W^{\times})$}\\
g_2, &\text{if $f = \{v^2_{s_2}, v^1_{t_2}\} \in E(W^{\times})$}.
\end{array}
\right.
\end{align}

We will show that $W^{\times}$ is smaller than $\overline{W}$.
To this end, note that we remove (at least) the subpaths $\overline{P}_1[v^1_{s_1}, v^1_{t_2}]$ and $\overline{P}_2[v^2_{s_2}, v^2_{t_1}]$ when creating $W^{\times}$.
The next result bounds how many edges are contained in these subpaths.

\begin{claim}\label{lem:indices}
We have $s_1 < t_2$ and $s_2 < t_1$.
Moreover, $s_1 + 1 < t_2$ or $s_2 + 1 < t_1$.
\end{claim}
\begin{cpf}
The vertices $v^1_{s_1}$ and $v^1_{t_2}$ are distinct because $g_1 \neq g_2$ by \cref{claim:g1_neq_g2} and because $G^*$ is disjoint.
Thus, $s_1 \neq t_2$.
The minimality of $s_1$ implies $s_1 < t_2$.
Similarly, we can show that $s_2 < t_1$.

Now, consider the conflict $g_1$.
By \cref{claim:conflict_edge_f_i}, there exists $f_1 \in E(\overline{P}_1) \cup E(\overline{P}_2)$ such that $\overline{\Phi}(f_1) = g_1$.
If $f_1 \in E(\overline{P}_1)$, then $f_1= \{v^1_{s_1}, v^1_{s_1+1}\}$ by \cref{claim:conflict_edge_f_i}.
Hence, as $g_1 \neq g_2$, the disjointness of $G^*$ yields that $v^1_{s_1+1} \neq v^1_{t_2}$.
Consequently, $s_1 + 1 < t_2$.
If $f_1 \in E(\overline{P}_2)$, then $f_1= \{v^2_{t_1 -1}, v^2_{t_1}\}$ by \cref{claim:conflict_edge_f_i}.
Hence, $v^2_{t_1-1} \neq v^2_{s_2}$ because $g_1 \neq g_2$ and $G^*$ is disjoint.
Consequently, $s_2 < t_1 - 1$.
\end{cpf}

\begin{claim}\label{claim:xShorter}
We have $|E(W^{\times})| < |E(\overline{W})|$.
\end{claim}
\begin{cpf}
By the definition of $W^{\times}$, the set $E(W^{\times}) \setminus \{\{v^1_{s_1}, v^2_{t_1}\}, \{v^2_{s_2}, v^1_{t_2}\}\}$ is contained in $E(\overline{W})$.
Hence, there are at most two hyperedges in $W^{\times}$ that do not appear in $\overline{W}$.
Conversely, as $s_1 < t_2$ and $s_2 < t_1$ by \cref{lem:indices}, the edges $\{v^1_{s_1}, v^1_{s_1+1}\}, \{v^2_{s_2}, v^2_{s_2+1}\} \in E(\overline{W})$ are not in $E(W^{\times})$.
Moreover, by \cref{lem:indices}, we have $s_1 + 1 < t_2$ or $s_2 + 1 < t_1$.
In the former case, the edge $\{v^1_{s_1+1}, v^1_{s_1+2}\} \in E(\overline{P}_1)$ is not in $E(W^{\times})$.
In the latter case, the edge $\{v^2_{s_2+1}, v^2_{s_2+2}\} \in E(\overline{P}_2)$ is not in $E(W^{\times})$.
Hence, there are at least three hyperedges in $\overline{W}$ that do not appear in ${W}^{\times}$.
Thus, $|E(W^{\times})| < |E(\overline{W})|$.
\end{cpf}
\smallskip

In the remainder of this section, we will show that $(W^{\times},\Phi^{\times})\qin G^*$ satisfies Properties \ref{property:Paths1} - \ref{property:Paths4}.
Together with \cref{claim:xShorter}, this will contradict the minimality of $(\overline{W},\overline{\Phi})$ and complete the proof of \cref{lem:MinWOTH}.

\begin{claim}
    The hypergraph $W^{\times}$ satisfies Property~\ref{property:Paths1}.
\end{claim}
\begin{cpf}
For $i \in [2]$, the definitions of $P^{\times}_i$ and $\Omega_i$ imply that $V_{\internal}(P^{\times}_i) \subseteq V_{\internal}(\Omega_i) \subseteq V_{\internal}(\overline{P}_1) \cup V_{\internal}(\overline{P}_2)$.
Moreover, $V_{\internal}(\overline{P}_i) \cap \overline{h}  =\emptyset$ for each $i \in[3]$ by Property \ref{property:Paths1} for $\overline{W}$.
By construction, $P^{\times}_i$ is an $\overline{r}$-$\overline{\ell}_i$-path for each $i \in [2]$.
Hence, Property~\ref{property:Paths1} holds for $W^{\times}$.
\end{cpf}
\smallskip

We chose $s_1$ and $t_1$ (respectively, $s_2$ and $t_2$) so that crossing over with $\{v^1_{s_1}, v^2_{t_1}\}$ (respectively, $\{v^2_{s_2}, v^1_{t_2}\}$) ensures $W^{\times}$ contains at most two vertices from $g_1$ (respectively, $g_2$).
We prove this and then use it to show $(W^{\times},\Phi^{\times})$ satisfies Property \ref{property:Paths2}.

\begin{claim} \label{lemma:g_prime_no_conflict}
For each $i \in [2]$, we have
\[
(\Phi^{\times})^{-1}(g_i)=
\left\{
\begin{array}{cl}
	\{\{v^i_{s_i},v^{3-i}_{t_{i}}\}\}=\{g_i\cap V(W^{\times})\}, & \text{if $\{v^i_{s_i}, v^{3-i}_{t_{i}}\} \in E(W^{\times})$} \\[.1cm]
	\emptyset, & \text{otherwise}.
\end{array}
\right.
\]
\end{claim}
\begin{cpf}
We only show the statement for $i = 1$ as the proof for $i=2$ is analogous.
By the definition of $W^{\times}$, we have
\begin{align}
g_1 \cap V(W^{\times}) \subseteq  \phantom{\cup}~& g_1 \cap V(\overline{P}_1[\overline{r}, v^1_{s_1}] + \overline{P}_1[v^1_{t_2}, \overline{\ell}_1])\label{eqg1p1}\\
 \cup  ~&g_1 \cap V(\overline{P}_2[\overline{r}, v^2_{s_2}] + \overline{P}_2[v^2_{t_1}, \overline{\ell}_2])\label{eqg1p2}\\
 \cup ~&g_1 \cap V(\overline{P}_3).\label{eqg1p3}
\end{align}

By \cref{lem:indices}, we have $s_1+1 \le t_2$.
By \cref{claim:g1_neq_g2} and because $G^*$ is disjoint, we have $v^1_{t_2} \notin g_1$.
Hence, \cref{claim:conflict_edge_f_i} yields that the intersection in \eqref{eqg1p1} is $\{v^1_{s_1}\}$. 
By \ref{prop:bar3}, we further have $v^1_{s_1} \in V(W^{\times})$ if and only if $P_2^\times=\Omega_2$.
Moreover, for $f \in E(W^{\times}) \cap E(\overline{P}_1)$, we have $\Phi^{\times}(f) \neq g_1$ by Property \ref{def:quasi_1} for $ \overline{\Phi}\upharpoonright_{E(\overline{W})\cap E(W^{\times})}$. 

Similarly, the intersection in \eqref{eqg1p2} is $\{v^2_{t_1}\}$ 
and we have $v^2_{t_1} \in V(W^{\times})$ if and only if $P_2^{\times} = \Omega_2$.
Moreover, for $f \in E(W^{\times}) \cap E(\overline{P}_2)$, we have $\Phi^{\times}(f) \neq g_1$ by Property \ref{def:quasi_1} for $ \overline{\Phi}\upharpoonright_{E(\overline{W})\cap E(W^{\times})}$.

The intersection in \eqref{eqg1p3} is empty by \eqref{eq:g_1_2_dont_intersect_P_3}.
In particular, for $f \in E(\overline{P}_3)$, we have $\Phi^{\times}(f) \neq g_1$ by Property \ref{def:quasi_1} for $ \overline{\Phi}\upharpoonright_{E(\overline{W})\cap E(W^{\times})}$.

Hence, $\{v^1_{s_1},v^{2}_{t_{1}}\}=g_1\cap V(W^{\times})$ if $\{v^1_{s_1}, v^{2}_{t_{1}}\} \in E(W^{\times})$.
The definition of $\Phi^{\times}$ in \eqref{eq:def_Phi_prime} concludes the proof.
\end{cpf}

\begin{claim}\label{lemT'quasi}
	The pair $(W^{\times},\Phi^{\times})$ satisfies Property \ref{property:Paths2}.
\end{claim}
\begin{cpf}
	By Property \ref{property:Paths2} for $(\overline{W}, \overline{\Phi})$, we have $(\overline{W}, \overline{\Phi}) \qin G^*$.
	Thus, we only need to check Property \ref{def:quasi_1} for $\{v^1_{s_1},v^{2}_{t_{2}}\}$ and $\{v^2_{s_2},v^{1}_{t_{1}}\}$, provided they exist in $E(W^{\times})$, and Property \ref{def:quasi_2} for $g_1,g_2 \in E(G^*)$.
	This all follows from \cref{lemma:g_prime_no_conflict}.
\end{cpf}

We prove Properties \ref{property:Paths3} and \ref{property:Paths4} for $(W^{\times}, \Phi^{\times})$ by leveraging the same properties for $(\overline{W}, \overline{\Phi})$.
To this end, we use the following claim:

\begin{claim}\label{claim:times_conflicts_are_in_overline}
If $e \in E(G^*)$ is a conflict for $(W^{\times},\Phi^{\times})$, then $e$ is a conflict for $(\overline{W}, \overline{\Phi})$.
Moreover, $e \notin \{g_1, g_2\}$ and $e \cap \{v^1_{s_1}, v^2_{s_2}, v^1_{t_2}, v^2_{t_1}\} = \emptyset$.
\end{claim}
\begin{cpf}
As $e$ is a conflict for $(W^{\times}, \Phi^{\times})$, there exists $f \in E(W^{\times})$ such that $\Phi^{\times}(f) = e$ and $f \subsetneq \Phi^{\times}(f) \cap V(W^{\times})$.
It follows from \cref{lemma:g_prime_no_conflict} that $e \notin \{g_1, g_2\}$ because for $g_1$ or $g_2$, no edge $f$ with the mentioned properties can exist.
Hence, $f \notin \{ \{v^1_{s_1}, v^2_{t_1}\}, \{v^2_{s_2}, v^1_{t_2}\}\}$.
Thus, $f \in E(\overline{W})\cap E(W^{\times})$.
By the definition of $\Phi^{\times}$ in \eqref{eq:def_Phi_prime}, we have $\Phi^{\times}(f) = \overline{\Phi}(f)$.
As $V(W^{\times}) \subseteq V(\overline{W})$, we have 
\(
f \subsetneq \Phi^{\times}(f) \cap V(W^{\times}) \subseteq \overline{\Phi}(f) \cap V(\overline{W}).
\) 
Thus, $e$ is a conflict for $(\overline{W}, \overline{\Phi})$.

By Property \ref{prop:bar1}, we have $e \in E_{\ge 4}(G^*)$.
As $e \notin \{g_1, g_2\}$ and $G^*$ is disjoint, we have $e \cap \{v^1_{s_1}, v^2_{s_2}, v^1_{t_2}, v^2_{t_1}\} = \emptyset$.
\end{cpf}

\begin{claim}
The pair $(W^{\times},\Phi^{\times})$ satisfies Property \ref{property:Paths3}. 
\end{claim}
\begin{cpf}
Let $e\in E_{\geq 4}(G^*)$ be a conflict for $(W^{\times}, \Phi^{\times})$.
By \cref{claim:times_conflicts_are_in_overline}, $e$ is a conflict for $(\overline{W}, \overline{\Phi})$.
We assume $e \cap V(P^{\times}_2) \neq \emptyset$ and show $e \cap V(\overline{P}_3) = \emptyset$.

We have 
\(
V(P^{\times}_2) \subseteq V(\Omega_2) = V(\overline{P}_1[\overline{r},v^1_{s_1}])  \cup V(\overline{P}_2[v^2_{t_1}, \overline{\ell}_2])
\)
by the definition of $P^{\times}_2$ in \eqref{eq:CrossOverPath2}.
Thus, $e$ intersects $V(\overline{P}_1[\overline{r},v^1_{s_1}]) $ or $V(\overline{P}_2[v^2_{t_1}, \overline{\ell}_2])$.

If $e \cap V(\overline{P}_1[\overline{r},v^1_{s_1}]) \neq \emptyset$, then the minimality of $s_1$ implies that $v^1_{s_1} \in e$, which contradicts \cref{claim:times_conflicts_are_in_overline}.
If $e \cap V(\overline{P}_2[v^2_{t_1}, \overline{\ell}_2]) \neq \emptyset$, then $e \cap V(\overline{P}_3) = \emptyset$ by Property \ref{property:Paths3} for $(\overline{W}, \overline{\Phi})$.
\end{cpf}
\begin{claim}\label{lemT'SatisfiesP5}
The pair $(W^{\times},\Phi^{\times})\qin G^*$ satisfies Property \ref{property:Paths4}.
\end{claim}
\begin{cpf}
Let $e\in E_{\ge 4}(G^*)$ be a conflict for $(W^\times, \Phi^\times)$ and let $i \in [3]$.
By \cref{claim:times_conflicts_are_in_overline}, $e$ is also a conflict for $(\overline{W}, \overline{\Phi})$.
If $i = 3$, then Property \ref{property:Paths4} holds because it holds for $\overline{W}$ and because $\Phi^\times \upharpoonright_{E(\overline{P}_3)}= \overline{\Phi}\upharpoonright_{E(\overline{P}_3)}$.
The proofs of Property \ref{property:Paths4} are analogous for $i=1$ and $i=2$, so we prove the statement only for $i = 1$.

If $e\cap V(\overline{P}_{2}[\overline{r},v^{2}_{s_{2}}]) \neq \emptyset$, then the minimality of $s_2$ implies that $v^2_{s_2} \in e$, which contradicts \cref{claim:times_conflicts_are_in_overline}.
Hence, $e\cap V(\overline{P}_{2}[\overline{r},v^{2}_{s_{2}}])=\emptyset$.
By definition of $P^{\times}_1$, we have $e \cap V(P^{\times}_1) \subseteq e \cap V(\overline{P}_1[v^1_{t_2}, \overline{\ell}_1])$.
The desired statement then follows from \ref{property:Paths4} for $(\overline{W}, \overline{\Phi})$ and because 
$\Phi^\times(f)=\overline{\Phi}(f)$ for all $f\in E(\overline{P}_1[v^{1}_{t_2},\overline{\ell}_1])$.
\end{cpf}
%

%%%%%%%%%%%%%%%%%%%%%%%%%%%%%%
%%%%%%%%%%%%%%%%%%%%%%%%%%%%%%
\section[Proof of Theorem~\ref{theorem:non_unimodular_laminar}]{Proof of \cref{theorem:non_unimodular_laminar}}\label{sec:proof_main_thm}
%%%%%%%%%%%%%%%%%%%%%%%%%%%%%%
%%%%%%%%%%%%%%%%%%%%%%%%%%%%%%
%
We repeat the statement for convenience. 

\TheoremNonUnimodularLaminar*
For the proof of \cref{theorem:non_unimodular_laminar}, we need to recap and introduce some notation related to laminar set families.
\begin{definition}
	A set family $\mathcal{L}$ is {\bf laminar} if for any $e,f\in \mathcal{L}$, $e\cap f= \emptyset$, $e\subseteq f$ or $f\subseteq e$.
	For $e\in \mathcal{L}$, we denote by $\children_{\mathcal{L}}(e)$ the inclusion-wise maximal proper subsets of $e$ in $\mathcal{L}$. 
	We call $f\in \mathcal{L}\setminus \{e\}$ a {\bf descendant} of $e$ if $f\subseteq e$ and an {\bf ancestor} of $e$ if $e\subseteq f$. 
	We denote the set of descendants of $e$ by $\descendants_{\mathcal{L}}(e)$. 
	We call $e\in\mathcal{L}$ with $\children_{\mathcal{L}}(e)=\emptyset$ a {\bf leaf} of $\mathcal{L}$ and denote the set of leaves of $\mathcal{L}$ by $\leaves(\mathcal{L})$.
\end{definition}

The laminar set families we work with always contain pairwise distinct sets (as they will form part of the edge set of a hypergraph whose incidence matrix has a non-zero determinant). 
However, it is more convenient to allow for multisets rather than having to establish the distinctness of the sets explicitly. 
For multisets, we define the parent-child relationship on equal sets in such a way they form a chain and all strict subsets are descendants of the lowest set in the chain.

\begin{definition}
	Given a laminar set family $\mathcal{L}$ and $e\in\mathcal{L}$, we define the \textbf{disjointification} of $e$ to be $\dis_{\mathcal{L}}(e)\coloneqq e\setminus \bigcup\children_{\mathcal{L}}(e)$.
\end{definition}

In case the family $\mathcal{L}$ is clear from context, we may omit the subscript $\mathcal{L}$ and simply write $\children(e)$, $\descendants(e)$ or $\dis(e)$.
We observe that the union of the children of $e\in\mathcal{L}$ is equal to the union of all of its descendants, which we will use implicitly in some of the subsequent arguments.
\begin{observation}
	Let $\mathcal{L}$ be a laminar set family and let $e\in\mathcal{L}$.\\
	 Then $\dis(e)=e\setminus \bigcup \descendants(e)$.	
\end{observation}

In this section, we prove \cref{theorem:non_unimodular_laminar} by induction on $|V|+|E|$.
To this end, let $G=(V,E)$ be a non-unimodular laminar hypergraph and assume that the statement is true for all non-unimodular laminar hypergraphs $G'=(V',E')$ with $|V'|+|E'|<|V|+|E|$. 
By \cref{thm:camion}, $G$ contains a non-unimodular Eulerian hypergraph $G'$ with $|\supp \mbfs{M}(G')|\equiv 2 \mod 4$.  
We may assume that $G=G'$ as we can apply the induction hypothesis to $G'$ otherwise.
Hence, 
\begin{equation}
	\text{ $G$ is Eulerian and $|\supp \mbfs{M}(G)|\equiv 2 \mod 4$.} \label{eq:G_Eulerian}
\end{equation}

If $G$ is disjoint, then \cref{theorem:size_four} tells us that $G$ contains a tree house or an odd cycle.
Hence, we may assume that $G$ is not disjoint in the following. Let $\mathcal{L}\coloneqq E_{\ge 4}(G)$ and let $e_\circ\in \mathcal{L}$ be a $\subseteq$-maximal hyperedge that is not a leaf of $\mathcal{L}$.  
Define a new hyperedge $f_\circ\coloneqq \dis_{\mathcal{L}}(e_\circ)$ and obtain a new laminar hypergraph $H$ from $G$ by deleting $\children_{\mathcal{L}}(e_\circ)$ and replacing $e_\circ$ by $f_\circ$. 
More precisely, we let 
\begin{equation}\mathcal{L}'\coloneqq \mathcal{L}\setminus (\children_{\mathcal{L}}(e_\circ)\cup \{e_\circ\})\cup \{f_\circ\} \text{ and } H\coloneqq (V,E_2(G)\dot{\cup} \mathcal{L}').\label{eq:def_H_laminar}\end{equation}
Note that $E(G)=E_2(G)\cup E_{\ge 4}(G)$ because $G$ is Eulerian.
 Further observe that $E_{\ge 4}(H)\subseteq \mathcal{L}'$.
  The set family $\mathcal{L}'$ is laminar because $\mathcal{L}$ is laminar and maximality of $e_\circ$ implies $f_\circ \cap e=\emptyset$ for every $e\in\mathcal{L}'\setminus \{f_\circ\}$. 
  Hence, $H$ is a laminar hypergraph.
\begin{lemma}
	$H$ is Eulerian and $|\supp (\mbfs{M}(H))|\equiv 2\mod 4$.
\end{lemma}
\begin{proof}
	As $G$ is Eulerian, all hyperedges of $G$ have even size. Using that $f_\circ=e_\circ\setminus \bigcup\children(e_\circ)$ must have even size as well, all hyperedges of $H$ have even size.
	
	Again using that $G$ is Eulerian, we know that all vertices of $G$ have even degree. The degree of a vertex in $\bigcup\children_{\mathcal{L}}(e_\circ)$ decreases by $2$ when passing from $G$ to $H$, the degree of all other vertices remains the same.
	Hence, $H$ is Eulerian.
	
	Finally, we calculate
	\begin{align*}
		|\supp (\mbfs{M}(H))|&=\sum_{e\in E_2(G)} |e|+\sum_{e\in\mathcal{L}'} |e|\\
		&=\sum_{e\in E_2(G)} |e|+\sum_{e\in\mathcal{L}} |e|-|e_\circ|-\sum_{e\in \children_{\mathcal{L}}(e_\circ)}|e|+|f_\circ|\\
		&=|\supp (\mbfs{M}(G))|-2\cdot \sum_{e\in \children_{\mathcal{L}}(e_\circ)}|e|\\
		&\equiv |\supp (\mbfs{M}(G))|\equiv 2 \mod 4,
	\end{align*}
	where the last identity follows from the fact that all hyperedges of $G$ have even size.
\end{proof}
As $|V(H)|=|V|=|V(G)|$ and $|E(G)|-|E(H)|=|\children_{\mathcal{L}}(e_\circ)|>0$ (because $e_\circ$ is not a leaf), our induction hypothesis implies that there is $J=H[U,F]\pin H$ which constitutes an odd cycle, a tree house, or an avocado. Our goal is to transform $J$ into one of the respective structures in $G$.

We first observe that $H$ constitutes a quasi-subhypergraph of $G$.
\begin{equation}
\text{	Define $\Phi\colon E(H)\rightarrow E(G)$ by $\Phi(e)=e$ if $e\ne f_\circ$ and $\Phi(f_\circ)=e_\circ$.} \label{eq:def_Phi_laminar}
\end{equation}
\begin{proposition}\label{prop:H_quasi}
We have	$(H,\Phi)\qin G$, $\Phi$ is injective, and $\Phi(E(H))=E(G)\setminus \children_{\mathcal{L}}(e_\circ)$. Moreover, $e_\circ$ is the unique conflict and we have $f_\circ = \dis_{\mathcal{L}}(e_\circ)$.
\end{proposition}
\begin{proof}
Injectivity of $\Phi$, as well as $\Phi(E(H))=E(G)\setminus \children_{\mathcal{L}}(e_\circ)$ follow immediately from \eqref{eq:def_H_laminar} and \eqref{eq:def_Phi_laminar}. The latter further implies \ref{def:quasi_1} for $(H,\Phi)$, while \ref{def:quasi_2} follows by injectivity of $\Phi$.
Using that $e=\Phi(e)=\Phi(e)\cap V(H)$ for $e\ne f_\circ$ and \[f_\circ =\dis_{\mathcal{L}}(e_\circ)= e_\circ\setminus \children_{\mathcal{L}}(e_\circ)\subsetneq e_\circ = \Phi(f_\circ)\cap V(H)\] tells us that $e_\circ$ is the unique conflict.

\end{proof}
By \cref{prop:H_quasi,lem:restriction}, we also have $(J,\Phi\upharpoonright J)\qin G$. Moreover, if $e_\circ\notin\Phi(E(J))$, then \cref{cor:restriction} implies $(J,\Phi\upharpoonright J)\pin G$ and we are done. 
Hence, we may assume in the following that $e_\circ$ constitutes the unique conflict for $(J,\Phi\upharpoonright J)$. In particular, by injectivity of $\Phi$, we must have $f_\circ\cap V(J)\in E(J)$. We further know that $f_\circ\cap V(J)=\dis_{\mathcal{L}}(e_\circ)\cap V(J)=(e_\circ \setminus \bigcup\children_{\mathcal{L}}(e_\circ))\cap V(J)$, where $\children_{\mathcal{L}}(e_\circ)$ forms a disjoint set family. Finally, we know that $e_\circ\in\mathcal{L}$ is a maximal set. \cref{lemma:odd_cycle,lemma:odd_tree_house,lemma:avocado} conclude the proof of \cref{theorem:non_unimodular_laminar}. While \cref{lemma:odd_tree_house} has been proven in \cref{subsubsec:no_odd_tree_house}, we prove \cref{lemma:odd_cycle,lemma:avocado} in the following subsections.

\begin{restatable}{lemma}{LemmaOddCycle}\label[lemma]{lemma:odd_cycle}
	Let $G$ be an Eulerian laminar hypergraph and let $\mathcal{L}\coloneqq E_{\ge 4}(G)$. Suppose that there is a quasi-subhypergraph $(C,\Phi)\qin G$ with the following properties:
	\begin{enumerate}[{\it({OQ}1)}, leftmargin = *]
		\item \label{item:laminar_odd_cycle_1}$C$ is an odd cycle.
		\item \label{item:laminar_odd_cycle_2}$\Phi$ is injective.
		\item \label{item:laminar_odd_cycle_3}There is a unique conflict $e_c\in E(G)$, which further satisfies $e_c\in\mathcal{L}$. 
		\item \label{item:laminar_odd_cycle_4}Let $f_c\in E(C)$ with $\Phi(f_c)=e_c$. There exists a disjoint set family $\mathcal{D}\subseteq \descendants_{\mathcal{L}}(e_c)$ such that $f_c=(e_c\setminus \bigcup \mathcal{D})\cap V(C)$. Moreover, for every $e\in E(C)\setminus \{f_c\}$ with $\Phi(e)\in\mathcal{L}$, we have $\descendants_{\mathcal{L}}(\Phi(e))\cap \mathcal{D}=\emptyset$.
	\end{enumerate}
	Then $G$ contains an odd cycle, a tree house, or an avocado.
\end{restatable}
We point out that if $J$ is an odd cycle, then $(J,\Phi\upharpoonright J)$ satisfies \ref{item:laminar_odd_cycle_4}. To this end, we choose $\mathcal{D}\coloneqq \children_{\mathcal{L}}(e_\circ)$. By maximality of $e_\circ$, $e_\circ$ is the unique ancestor of one of its children. As $\Phi$ is injective, the last part of \ref{item:laminar_odd_cycle_4} is satisfied. 

For completeness, we restate \cref{lemma:odd_tree_house}.
\LemmaOddTreeHouse*
We point out that if $J$ is a tree house, then $(J,\Phi\upharpoonright J)$ satisfies \ref{item:laminar_odd_tree_house_4}, so \cref{lemma:odd_tree_house} is indeed applicable.
To see this, note that $f_c=f_\circ\cap V(J)$ and suppose that $f_c\ne h$. As $|h|\ge 4$ and $\Phi$ is injective, we have $\Phi(h)\in \mathcal{L}\setminus \{e_\circ\}$. As $e_\circ\in\mathcal{L}$ is maximal, $f_\circ=\dis_{\mathcal{L}}(e_\circ)$ must be disjoint from $\Phi(h)$. In particular, also $f_c\cap h=\emptyset$.

\begin{restatable}{lemma}{LemmaAvocado}\label[lemma]{lemma:avocado}
	Let $G$ be an Eulerian laminar hypergraph and let $\mathcal{L}\coloneqq E_{\ge 4}(G)$. Suppose that there is a quasi-subhypergraph $(A,\Phi)\qin G$ with the following properties:
	\begin{enumerate}[{\it({AQ}1)}, leftmargin = *]
		\item \label{item:laminar_avocado_1}$A$ is an avocado.
		\item \label{item:laminar_avocado_2}$\Phi$ is injective.
		\item \label{item:laminar_avocado_3}There is a unique conflict $e_c\in E(G)$, which constitutes a maximal set in $\mathcal{L}$. 
		\item \label{item:laminar_avocado_4}Let $f_c\in E(A)$ with $\Phi(f_c)=e_c$. We have $f_c=\dis_{\mathcal{L}}(e_c)\cap V(A)$.
	\end{enumerate}
	Then $G$ contains an odd cycle, a tree house, or an avocado.
\end{restatable}

\subsection[The odd cycle case (proof of Lemma~\ref{lemma:odd_cycle})]{The odd cycle case (proof of \cref{lemma:odd_cycle})}

This section is dedicated to the proof of \cref{lemma:odd_cycle}.
We will assume in the following that
\begin{equation}
	\text{among all $(C',\Phi')\qin G$ satisfying \ref{item:laminar_odd_cycle_1}-\ref{item:laminar_odd_cycle_4}, $C$ is of minimum length.}\label{eq:C_minimal}
\end{equation}
We will further assume that
\begin{equation}
	\text{$G$ does not contain an odd cycle}. \label{eq:no_odd_cycle}
\end{equation}
We can make this assumption without loss of generality as we are done otherwise.
Let $P\coloneqq C-f_c$ and denote the endpoints of $P$ by $v$ and $w$.
\begin{proposition}
We have $e_c\notin\Phi(E(P))$ and $(P\,\Phi\upharpoonright_P)\pin G$.\label{prop:P_conflict_free}
\end{proposition}
\begin{proof}
	As $\Phi$ is injective and $f_c\notin E(P)$, we have $e_c\notin \Phi(E(C))$. As $e_c$ is the unique conflict, $(P,\Phi\upharpoonright_P)\pin G$ by \cref{cor:restriction}.
\end{proof}
\begin{lemma}\label{lemma:consecutive_vertices_odd_distance}
	Let $e\in E(G)$ and let $Q$ be a subpath of $P$ with endpoints $x,y\in e$ such that no inner vertex of $Q$ is contained in $e$. Then $Q$ is odd.
\end{lemma}
\begin{proof}
	As $Q\pin P$, \cref{prop:P_conflict_free} implies  $(Q,\Phi\upharpoonright_{Q})\pin G$. If $Q$ has length $1$, we are done. Otherwise, no edge of $Q$ has both endpoints in $e$, implying $e\notin \Phi(E(Q))$. Now, \cref{lemma:parities}~\ref{lemma:parities:1} and \eqref{eq:no_odd_cycle} conclude the proof.
\end{proof}
\begin{corollary}\label{cor:odd_number_of_vertices_from_e_c}
	$|e_c\cap V(C)|$ is odd.
\end{corollary}
\begin{proof}
	Let $u_0=v,u_1,\dots,u_m,u_{m+1}=w$ be the sequence of vertices in $e_c\cap V(C)$ in the order in which they appear on $P=C-f_c$. By \cref{lemma:consecutive_vertices_odd_distance}, $P[u_i,u_{i+1}]$ is odd for every $i\in \{0,\dots,m\}$. As $P$ has even length, $m$ must be odd.
\end{proof}
\begin{lemma}
	There is a unique $g\in \mathcal{D}$ with $V(C)\cap g\ne \emptyset$.
\end{lemma}
\begin{proof}
	By \ref{item:laminar_odd_cycle_3}, we know that $e_c$ constitutes a conflict, i.e., $f_c\subsetneq e_c\cap V(C)$. By \ref{item:laminar_odd_cycle_4}, we can hence infer that there is at least one $g\in\mathcal{D}$ with $V(C)\cap g\ne \emptyset$. Let $\{g_1,\dots,g_k\}\coloneqq \{g\in\mathcal{D}\colon g\cap V(C)\ne \emptyset\}$ and assume towards a contradiction that $k\ge 2$. By definition of $\mathcal{D}$, $g_1,\dots,g_k$ are pairwise disjoint. Let $v=u_0,u_1,\dots,u_{m+1}=w$ be the sequence of vertices in $V(C)\cap e_c$ in the order in which they appear on $P$. 
	By \ref{item:laminar_odd_cycle_4} and since the $g_i$ are disjoint, for $j\in [m]$, there exists a unique $\alpha(j)\in [k]$ such that $u_j\in g_{\alpha(j)}$. 
	We partition the index set $[m]$ into intervals $J_1,\dots,J_\ell$ such that for each interval $J_h=[a_h,b_h]$, there is $i\in [k]$ with $\alpha(j)=i$ for every $j\in J_h$, but $u_{a_h-1},u_{b_h+1}\notin g_{i}$. 
	This is possible because $u_0=v$ and $u_{m+1}=w$ are not contained in any $g_i$.
	\begin{claim}
		For every $h\in [\ell]$, $b_h-a_h$ is odd, i.e., $J_h$ contains an even number of indices.
	\end{claim}
	\begin{cpf}
		Suppose towards a contradiction that $b_h-a_h$ is even for some $h\in [\ell]$. Let $i\coloneqq \alpha(a_h)$, let $j_{min}\coloneqq a_h-1$ and $j_{max}\coloneqq b_h+1$. By our assumption that $k\ge 2$, there are at least two distinct intervals, i.e., $\ell\ge 2$. In particular, $[j_{min},j_{max}]\subsetneq [0,m+1]$. By \cref{lemma:consecutive_vertices_odd_distance}, $P[u_j,u_{j+1}]$ has odd length for every $j\in [j_{min},j_{max}-1]\cap \mathbb{Z}$. Hence, $Q\coloneqq P[u_{j_{min}},u_{j_{max}}]$ has even length because $j_{max}-j_{min}=b_h-a_h+2$ is even.
		Let $C'\coloneqq Q+\{u_{j_{min}},u_{j_{max}}\}$ and define $\Phi'\colon E(C')\rightarrow E(G)$ via $\Phi'\upharpoonright E(Q)=\Phi\upharpoonright E(Q)$ and $\Phi'(\{u_{j_{min}},u_{j_{max}}\})=e_c$.
		
		By construction, $C'$ is an odd cycle that is strictly shorter than $C$.
		In particular, \ref{item:laminar_odd_cycle_1} holds for $C'$.
		By \cref{prop:P_conflict_free}, $e_c\notin \Phi(E(P))$, so $\Phi'$ is injective. As further $\{u_{j_{min}},u_{j_{max}}\}\subseteq e_c$, $(C',\Phi')\qin G$, so \ref{item:laminar_odd_cycle_2} holds for $(C',\Phi')$. 
		
		As $Q\pin P$, \cref{prop:P_conflict_free} tells us that $(Q,\Phi'\upharpoonright Q)=(Q,\Phi\upharpoonright Q)\pin G$. On the other hand, $G$ does not contain an odd cycle by \eqref{eq:no_odd_cycle}, so $e_c$ must constitute the (unique) conflict of $(C',\Phi')$.
	Hence, property \ref{item:laminar_odd_cycle_3} holds for $(C',\Phi')$.
	
	 Finally, for \ref{item:laminar_odd_cycle_4}, define $\mathcal{D}'\coloneqq \{g_i\}$. As $\mathcal{D}'\subseteq \mathcal{D}$ and $E(C')\setminus \{\{u_{j_{min}},u_{j_{max}}\}\}=E(Q)\subseteq E(P)= E(C)\setminus \{f_c\}$, we have \[\descendants_{\mathcal{L}}(\Phi'(e))\cap \mathcal{D}'=\descendants_{\mathcal{L}}(\Phi(e))\cap \mathcal{D}'\subseteq \descendants_{\mathcal{L}}(\Phi(e))\cap \mathcal{D}=\emptyset\] for $e\in (E(C')\setminus \{\{u_{j_{min}},u_{j_{max}}\}\})\cap \mathcal{L}$.
		We further have \[(e_c\setminus g_i)\cap V(C')=\{u_j\colon j\in [j_{min},j_{max}]\cap \mathbb{Z}\}\setminus g_i=\{u_{j_{min}},u_{j_{max}}\}\] as required.
		Hence, $(C',\Phi')$ contradicts our assumption \eqref{eq:C_minimal}.
	\end{cpf}
	As a consequence, we know that $|e_c\cap V(C)|=m+2=2+\sum_{h=1}^\ell (b_h-a_h+1)$ is even. But this contradicts \cref{cor:odd_number_of_vertices_from_e_c}.
\end{proof}
Hence, let $g\in\mathcal{D}$ be unique with the property that $g\cap V(C)\ne \emptyset$.
\begin{lemma}
	If $|g\cap V(C)|=1$, then $G$ contains a tree house.
\end{lemma}
\begin{proof}
	Let $V(C)\cap g=\{r\}$. By \ref{item:laminar_odd_cycle_4}, we have
	\[\{v,w\}=f_c=V(C)\cap \left(e_c\setminus \bigcup \mathcal{D}\right)=(V(C)\cap e_c)\setminus (V(C)\cap g)=(V(C)\cap e_c)\setminus \{r\}.\] In particular, $r$ is not an endpoint of $f_c=\{v,w\}$. Hence, $r$ is an inner vertex of $P$ and $e_c\cap V(C)=\{v,w,r\}$. By \cref{lemma:consecutive_vertices_odd_distance}, both $P[r,v]$ and $P[r,w]$ are odd. As $G$ is Eulerian and $g\cap V(C)=g\cap V(P)=\{r\}$ has odd size, there is $\ell_3\in g\setminus V(C)$. 
	
	Let $J\coloneqq P+\ell_3$ arise from $P$ by adding the isolated vertex $\ell_3$.
	\begin{claim}
		We also have $(J,\Phi\upharpoonright_P)\pin G$.
	\end{claim}
	We remark that $E(P)=E(J)$, so $\Phi\upharpoonright_P$ is indeed a function with domain $E(J)$.
	\begin{cpf}
		By \cref{prop:P_conflict_free}, we have $(P,\Phi\upharpoonright P)\pin G$.  We need to show that for $e\in E(P)=E(C)\setminus \{f_c\}$, we have $\ell_3\notin \Phi(e)$ because then, $e=\Phi(e)\cap V(P)=\Phi(e)\cap V(J)$. 
	Suppose towards a contradiction that $e\in E(P)$ and $\ell_3\in \Phi(e)$. As $\ell_3\notin V(P)=V(C)$, $|\Phi(e)|\ge |e|+1=3$, so $\Phi(e)\in\mathcal{L}$. Moreover, $\ell_3\in \Phi(e)\cap g\ne \emptyset$, so $g$ is a descendant, equal to, or an ancestor of $\Phi(e)$. We have $|\Phi(e)\cap V(P)|=|e|=2$, but $|g\cap V(P)|=1$. Hence, we must have $g\in \descendants_{\mathcal{L}}(\Phi(e))\cap \mathcal{D}$, contradicting \ref{item:laminar_odd_cycle_4}. 
	\end{cpf}
	We have $g\notin \Phi(E(P))$ because $g\cap V(P)=\{r\}$, and $e_c\notin \Phi(E(P))$ by \cref{prop:P_conflict_free}.
	   Hence, $J\oplus g\oplus e_c\pin G$ by \cref{lem:addition} (applied twice). As $g\cap V(J)=\{\ell_3,r\}$ and $e_c\cap V(J)=\{r,v,w,\ell_3\}$, $J\oplus g\oplus e_c$ is a tree house.
\end{proof}
\begin{lemma}
	If $|g\cap V(C)|\ge 2$, then $G$ contains an avocado.
\end{lemma}
\begin{proof}
	By \ref{item:laminar_odd_cycle_4}, we have $V(P)\cap e_c=\{v,w\}\cup (V(C)\cap g)$. By \cref{lemma:consecutive_vertices_odd_distance}, using that $v$ and $w$ are the endpoints of $P$ and that $P=C-f_c$ is even, we know that $P$ contains an odd number of vertices from $e_c$, and, hence, from $g$. By the assumption of the lemma, $P$ contains at least $3$ vertices from $g$. Using $V(P)=V(C)$, we get
	\begin{equation}
	|V(C)\cap e_c|\equiv |V(C)\cap g|\equiv 1\mod 2\text{ and } |V(C)\cap g|\ge 3. \label{eq:at_least_3_vertices_from_g}	
	\end{equation}
	
	 Let $E_\circ\coloneqq \{e\in E(P)\colon \Phi(e)\in\descendants_{\mathcal{L}}(g)\}$. As for $e\in E_\circ$, $\Phi(e)\cap V(C)=\Phi(e)\cap V(P)=e$ by \cref{prop:P_conflict_free}, and no two edges of $P$ contain each other, no two of the sets $(\Phi(e))_{e\in E_\circ}$ contain each other. As $\mathcal{L}$ is a laminar family, the sets $(\Phi(e))_{e\in E_\circ}$ are pairwise disjoint. By \eqref{eq:at_least_3_vertices_from_g},
	  \[\left|V(C)\cap \left(g\setminus \bigcup_{e\in E_\circ} \Phi(e)\right)\right|=|V(C)\cap g|-\sum_{e\in E_\circ} |e|=|V(C)\cap g|-2\cdot |E_\circ|\]
	   is odd. As $G$ is Eulerian, $|(g\setminus \bigcup_{e\in E_\circ} \Phi(e))|$ is even. Hence, there must be a vertex $x\in (g\setminus \bigcup_{e\in E_\circ} \Phi(e))\setminus V(C)$.
	   \[\]Let $A\coloneqq G[V(C)\cup \{x\},\Phi(E(P))\cup \{g,e_c\}]$.
	   
	We claim that $A$ is an avocado (with $V=V(C)$, $t=x$, $E=\Phi(E(P))[V(C)\cup \{x\}]$, $s=g\cap (V(C)\cup\{x\})$ and $f=e_c\cap (V(C)\cup\{x\})$. 
	
	We recall that by \cref{prop:P_conflict_free}, we have, for every $e\in E(P)$,  $e=\Phi(e)\cap V(P)=\Phi(e)\cap V(C)$.
	 If we had $x\in \Phi(e)$, then as $x\notin V(C)$, $|\Phi(e)|\ge |e|+1\ge 3$, so $\Phi(e)\in\mathcal{L}$. 
	As $x\in g$, $\Phi(e)$ must be a descendant, equal to, or an ancestor of $g$. As $|\Phi(e)\cap V(C)|=|e|=2$, but $|g\cap V(C)|\ge 3$ by \eqref{eq:at_least_3_vertices_from_g}, $\Phi(e)$ must be a descendant of $g$, so $e\in E_\circ$. But this implies $x\notin \Phi(e)$, a contradiction.
	Hence, we have
	\begin{equation}
		\text{	$e=\Phi(e)\cap (V(C)\cup \{x\})$ for every $e\in E(P)$.}\label{eq:edges_of_P_no_conflict}
	\end{equation} 
	
	We check the properties \ref{item:avocado_1}-\ref{item:avocado_3} from \cref{def:avocado} one by one.

	\ref{item:avocado_1}: By construction, $x\in g\subseteq e_c$. We have argued that $|g\cap V(C)|\ge 3$, so as $x\in g\setminus V(C)$, $|s|=|g\cap (V(C)\cup\{x\})|\ge 4$. 
	Finally, as $E=E(P)$ and $x\notin V(P)=V(C)$, $x$ is not incident to any edge in $E$.
	
	\ref{item:avocado_2}: By \eqref{eq:edges_of_P_no_conflict}, we have $(V,E)=P$. 
	Moreover, as $g$ is the unique hyperedge in $\mathcal{D}$ intersecting $V(C)$ and $x\in g$, we have	$(e_c\setminus g)\cap (V(C)\cup \{x\})=(e_c\setminus g)\cap V(C)=(e_c\setminus \bigcup \mathcal{D})\cap V(C)=f_c$ by \ref{item:laminar_odd_cycle_4}. 
	Hence, we recover the odd cycle $C$.
	
	\ref{item:avocado_3}: This follows from \cref{lemma:consecutive_vertices_odd_distance}. 
\end{proof}

\subsection[The avocado case (proof of Lemma~\ref{lemma:avocado})]{The avocado case (proof of \cref{lemma:avocado})}

This section is dedicated to the proof of \cref{lemma:avocado}, which we restate here for convenience.
\LemmaAvocado*
Let $A=(V\cup \{t\},E\cup\{s,f\})$ as in \cref{def:avocado}.
We first prove some useful statements about (the interactions between) different hyperedges in \cref{prop:Phi_s_f_in_L,lemma:f_c_edge,prop:f_s_no_conflict,prop:no_ancestor}.
We then use them to derive the existence of an odd cycle, a tree house or an avocado in $G$ in \cref{lemma:case_even_subpath,lemma:case_disjoint,lemma:case_non_disjoint}.
\begin{proposition}
	We have $\Phi(s)\in \mathcal{L}$ and $\Phi(f)\in\mathcal{L}$.\label{prop:Phi_s_f_in_L}
\end{proposition}
\begin{proof}
	This follows from $4\le |s|\le |\Phi(s)|$ and $4\le |s|\le |f|\le |\Phi(f)|$.
\end{proof}
\begin{lemma}
	We have $f_c\in E$. In particular, $e_c\notin \{\Phi(s),\Phi(f)\}$. \label{lemma:f_c_edge}
\end{lemma}
\begin{proof}
	We need to show that $f_c\notin \{s,f\}$, the second part of the statement then follows from injectivity of $\Phi$. 
	Suppose towards a contradiction that $f_c\in \{s,f\}$ and let $g\in \{s,f\}\setminus \{f_c\}$. 
	Using $s\subseteq f$, \ref{item:laminar_avocado_4} and $g\subseteq \Phi(g)$, we have
	\[\emptyset\neq s=s\cap f=f_c\cap g=\dis_{\mathcal{L}}(e_c)\cap V(A)\cap g\subseteq \dis_{\mathcal{L}}(e_c)\cap \Phi(g).\]
	By \cref{prop:Phi_s_f_in_L}, we have $\Phi(g)\in\mathcal{L}$, and injectivity of $\Phi$ implies that $e_c=\Phi(f_c)\ne \Phi(g)$.
	Thus, $\emptyset\ne \dis_{\mathcal{L}}(e_c)\cap \Phi(g)$ tells us that $\Phi(g)$ is an ancestor of $e_c$ in $\mathcal{L}$.
	But this contradicts the fact that $e_c$ is a maximal set in $\mathcal{L}$ by \ref{item:laminar_avocado_3}.
\end{proof}
\begin{proposition}\label{prop:f_s_no_conflict}
	We have $s=V(A)\cap \Phi(s)$, $f=V(A)\cap \Phi(f)$ and $f\setminus s= (\Phi(f)\setminus \Phi(s))\cap V(A)$. In particular, $\Phi(s)\in\descendants_{\mathcal{L}}(\Phi(f))$.
\end{proposition}
\begin{proof}
	The first three statements follow from \cref{lemma:f_c_edge} and the fact that $e_c$ is the unique conflict  by \ref{item:laminar_avocado_3}.
	 By \cref{prop:Phi_s_f_in_L}, $\Phi(s),\Phi(f)\in\mathcal{L}$. As $\emptyset\ne s\subseteq \Phi(s)$ and $s\subseteq f\subseteq\Phi(f)$, and $\Phi$ is injective, one of $\Phi(s)$ and $\Phi(f)$ is a descendant of the other. As $\emptyset \ne f\setminus s\subseteq \Phi(f)\setminus \Phi(s)$, $\Phi(s)$ must be a descendant of $\Phi(f)$.
\end{proof}
\begin{proposition}\label{prop:no_ancestor}
	Let $e\in E\setminus \{f_c\}$. Then $\Phi(e)$ is not an ancestor of $\Phi(s)$ or $\Phi(f)$ in $\mathcal{L}$.
\end{proposition}
\begin{proof}
	As $e\ne f_c$, $\Phi$ is injective and $e_c$ is the unique conflict by \ref{item:laminar_avocado_3}, we know that $\Phi(e)$ is not a conflict. In particular, $|\Phi(e)\cap V(A)|=|e|=2$, whereas $|\Phi(s)\cap V(A)|\ge |s|\ge 4$ and $|\Phi(f)\cap V(A)|\ge |f|\ge 4$. Hence, $\Phi(e)$ cannot be an ancestor of $\Phi(f)$ or $\Phi(s)$.
\end{proof}
Let $C=(V,E\dot{\cup} \{f\setminus s\})$ be the odd cycle from \cref{def:avocado}.
 Let $f_c\eqqcolon \{v,w\}$ and let $P\coloneqq C-f_c$. Let $v=u_0,u_1,\dots,u_m, u_{m+1}=w$ be the sequence of vertices from $e_c\cap V$ in the order in which they appear on $P$.
\begin{lemma}\label{lemma:case_even_subpath}
	If there exists $i\in \{0,\dots,m\}$ such that $P[u_i,u_{i+1}]$ has even length, then $G$ contains an odd cycle, a tree house, or an avocado.
\end{lemma}
\begin{proof}
	Let $Q\coloneqq P[u_i,u_{i+1}]$, $C'\coloneqq Q+\{u_i,u_{i+1}\}$ and define $\Phi'\colon E(C')\rightarrow E(G)$ via $\Phi'(e)=\Phi(e)$ for every $e\in E(Q)\cap E$, $\Phi'(f\setminus s)=\Phi(f)$ (if $f\setminus s\in E(Q)$) and $\Phi'(\{u_i,u_{i+1}\})=e_c$.
	By construction, $C'$ is an odd cycle.
	As $\Phi$ is injective, $f\ne f_c$ by \cref{lemma:f_c_edge}, and neither $f_c$ nor $f$ are contained in $E(Q)\cap E$, $\Phi'$ is injective. As $(A,\Phi)\qin G$, $f\setminus s\subseteq f\subseteq \Phi(f)$ and $\{u_i,u_{i+1}\}\subseteq e_c$, we have $(C',\Phi')\qin G$.
	
	For $e\in E(Q)\cap E$, we have $e=\Phi(e)\cap V(A)$ because $e_c=\Phi(f_c)$ is the unique conflict for $(A,\Phi)$, and $\Phi$ is injective. As $e\subseteq V(C')\subseteq V(A)$ and $\Phi'(e)=\Phi(e)$, we also have $e=\Phi'(e)\cap V(C')$, i.e., $\Phi'(e)$ is not a conflict for $(C',\Phi')$.
	
	By construction, we have $\{u_i,u_{i+1}\}=e_c\cap V(C')$, so $e_c$ is not a conflict for $(C',\Phi')$, either.
	
	If $f\setminus s\notin E(Q)$ or $\Phi(f)$ does not constitute a conflict for $(C',\Phi')$, then $C'\pin G$ and $G$ contains an odd cycle.
	
	In case $f\setminus s\in E(Q)$ and $\Phi(f)$ constitutes a conflict for $(C',\Phi')$, it is the unique one.
	Our goal is to show that $(C',\Phi')\qin G$ meets the requirements of \cref{lemma:odd_cycle}.
	We have already verified \ref{item:laminar_odd_cycle_1} and \ref{item:laminar_odd_cycle_2}, and \ref{item:laminar_odd_cycle_3} follows from the above and \cref{prop:Phi_s_f_in_L}.
	We are left with verifying \ref{item:laminar_odd_cycle_4}.
	 By \cref{prop:f_s_no_conflict}, we get $f\setminus s=V(A)\cap (\Phi(f)\setminus \Phi(s))$. Using $f\setminus s\subseteq V(C')\subseteq V(A)$, we obtain $f\setminus s=(\Phi(f)\setminus \Phi(s))\cap V(C')$.
	By \cref{prop:f_s_no_conflict}, we get $\Phi(s)\in \descendants_{\mathcal{L}}(\Phi(f))$, so we can choose the family $\mathcal{D}$ for \ref{item:laminar_odd_cycle_4} to be $\{\Phi(s)\}$.

	It remains to show that $\Phi(s)$ is not a descendant of $\Phi'(e)$ for $e\in E(Q)\setminus \{f\setminus s\}\cup \{\{u_i,u_{i+1}\}\}$ with $\Phi'(e)\in\mathcal{L}$. 
	By \cref{prop:no_ancestor}, we know that $\Phi(s)$ cannot be a descendant of $\Phi(e)$ for $e\in E(Q)\setminus \{f\setminus s\}\subseteq E\setminus \{f_c\}$. 
	Finally, suppose $\Phi(s)$ were a descendant of $e_c=\Phi'(\{u_i,u_{i+1}\})$.
	 By \cref{prop:f_s_no_conflict,lemma:f_c_edge} and because $e_c$ is maximal in $\mathcal{L}$ by \ref{item:laminar_avocado_3}, this implies that $\Phi(f)$ is a descendant of $e_c$ as well. But this contradicts the fact that $|\Phi(f)\cap V(C')|\ge |f\setminus s|+1\ge 3$ (as $\Phi(f)=\Phi'(f\setminus s)$ is a conflict), but $|e_c\cap V(C')|=2$.
	
	All in all, $(C',\Phi')\qin G$ satisfies all requirements of \cref{lemma:odd_cycle}, which yields the existence of an odd cycle, a tree house, or an avocado in $G$.
\end{proof}
Hence, we may assume in the following that
\begin{equation}
	\text{$P[u_i,u_{i+1}]$ has odd length for every $i\in \{0,\dots,m\}$.}\label{eq:all_paths_odd}
\end{equation}
Using that $P=C-f_c$ has even length, this implies that $|V\cap e_c|$ is odd.
Finally, as $v$ and $w$, the endpoints of $f_c\subseteq e_c$, are neighboring vertices on $C$, we know:
\begin{align}
	&\text{Let $Q$ be a subpath of $C$ such that both endpoints of $Q$ are in $e_c$, }\notag\\
	&\text{but no inner vertex of $Q$ is. Then $Q$ has odd length.}\notag
\end{align}
	As $e_c$ is the unique conflict by \ref{item:laminar_avocado_3} and as $\Phi(f)\ne e_c$ by \cref{prop:f_s_no_conflict}, we have $\Phi(f)\cap V(A)=f$.
By \ref{item:avocado_3}, and the above, we conclude the following:
\begin{align}
	&\text{Let $e\in \{e_c,\Phi(f)\}$ and let $Q$ be a subpath of $C$ such that both endpoints }\label{eq:odd_distance_both_hyperedges}\\
	&\text{of $Q$ are in $e$, but no inner vertex of $Q$ is. Then $Q$ has odd length.}\notag\\
	&\text{Both $|V\cap e_c|$ and $|V\cap \Phi(f)|$ are odd.} \label{eq:odd_number_of_vertices_from_both}
\end{align}
\begin{lemma}\label{lemma:case_disjoint}
	If $e_c\cap \Phi(f)=\emptyset$, then $G$ contains an odd cycle.
\end{lemma}
\begin{proof}
	We first prove the following claim, and then explain how to use it to construct an odd cycle.
	\begin{claim}
		There exist vertex-disjoint subpaths $Q_1$ and $Q_2$ of $P$ such that:
		\begin{itemize}
			\item Each of $Q_1$ and $Q_2$ has one endpoint in $e_c$ and one endpoint in $\Phi(f)$.
			\item No internal vertex of $Q_1$ or $Q_2$ is contained in $e_c$ or $\Phi(f)$.
			\item $|E(Q_1)|+|E(Q_2)|$ is odd.
		\end{itemize}
	\end{claim} 
	\begin{cpf}
		Recall that $\Phi(f)$ and $e_c$ are disjoint by the assumption of the lemma.
		Let $v=v_0,v_1,\dots,v_k,v_{k+1}=w$ be the sequence of vertices in $\Phi(f)\cup e_c$ in which they appear on $P$.
		By \eqref{eq:odd_number_of_vertices_from_both}, $k$ is even.
		Let 
		\begin{align*}I_1&\coloneqq \{i\in \{0,\dots,k\}\colon v_i\in e_c \text{ and } v_{i+1}\in \Phi(f)\} \text{ and }\\
			I_2&\coloneqq \{i\in \{0,\dots,k\}\colon v_i\in \Phi(f)\text{ and } v_{i+1}\in e_c\}.\end{align*}
		Note that as both $v$ and $w$ are contained in $e_c$, $|I_1|=|I_2|$ equals the number of `blocks' of consecutive indices $i$ with $v_i\in \Phi(f)$.
		For $i\in \{0,\dots,k\}\setminus (I_1\cup I_2)$, we know that $P[v_i,v_{i+1}]$ is odd by \eqref{eq:odd_distance_both_hyperedges}. As $k$ is even, we have
		\begin{align*}0&\equiv |E(P)|\equiv \sum_{i\in \{0,\dots,k\}\setminus (I_1\cup I_2)} \underbrace{|E(P[v_i,v_{i+1}])|}_{\equiv 1\mod 2}+\sum_{i\in I_1\cup I_2} |E(P[v_i,v_{i+1}])|\\
			&\equiv  \underbrace{(k+1)}_{\equiv 1 \mod 2}-\underbrace{(|I_1|+|I_2|)}_{=2\cdot |I_1|}+\sum_{i\in I_1\cup I_2} |E(P[v_i,v_{i+1}])|\\
			&\equiv 1+\sum_{i\in I_1\cup I_2} |E(P[v_i,v_{i+1}])|\mod 2.\end{align*}
		As $|I_1|+|I_2|=2\cdot |I_1|$ is even, the paths $P[v_i,v_{i+1}]$ with $i\in I_1\cup I_2$ cannot be all even or all odd. Hence, there exist $i\ne j\in I_1\cup I_2$ such that $P[v_i,v_{i+1}]$ is even and $P[v_j,v_{j+1}]$ is odd. We are done unless we have $i+1=j$ or $j+1=i$. Hence, assume that we have three consecutive indices $i,i+1$ and $i+2$ such that one of the paths $P[v_i,v_{i+1}]$ and $P[v_{i+1},v_{i+2}]$ is even and the other one is odd, and such that $v_i$ and $v_{i+2}$ are contained in the same one among the two hyperedges $e_c$ and $\Phi(f)$, whereas $v_{i+1}$ is contained in the other one.
		
		As $|V\cap \Phi(f)|=|f|-1=|s|+1\ge 5$ by \ref{item:avocado_1} and \ref{item:avocado_2}, we know that there exists $j\in \{1,\dots,k\}\setminus \{i,i+1,i+2\}$ with $v_j\in \Phi(f)$. If $j\le i-1$, then since $v_0\in e_c$, there is $\ell\le j-1\le i-2$ such that $v_\ell\in e_c$ and $v_{\ell+1}\in \Phi(f)$. As $P[v_i,v_{i+1}]$ and $P[v_{i+1},v_{i+2}]$ have different parities, on of the pairs $P[v_\ell,v_{\ell+1}]$ and $P[v_i,v_{i+1}]$ or $P[v_\ell,v_{\ell+1}]$ and $P[v_{i+1},v_{i+2}]$	has the desired properties. The case $j\ge i+2$ can be handled analogously, using $v_{k+1}\in e_c$.
	\end{cpf}
	Let $Q_1$ and $Q_2$ be the two paths promised by the claim and let $J\coloneqq Q_1+Q_2$. Note that no edge of $J$ has both endpoints in $\Phi(f)$ or in $e_c$.
	In particular, $J\pin (V,E\setminus \{f_c\})\pin A$ and $e_c,\Phi(f)\notin \Phi(E(J))$.
	 As $e_c$ is the unique conflict of $(A,\Phi)$, this implies $(J,\Phi\upharpoonright J)\pin G$ by \cref{cor:restriction}.
	  By \cref{lem:addition}, $(J,\Phi\upharpoonright J)\oplus e_c\oplus \Phi(f)\pin G$ is an odd cycle.
\end{proof}
\begin{lemma}\label{lemma:case_non_disjoint}
	If $e_c\cap \Phi(f)\ne \emptyset$, then $G$ contains an odd cycle, a tree house, or an avocado.
\end{lemma}
\begin{proof}
	By \cref{prop:Phi_s_f_in_L,prop:f_s_no_conflict} and since $e_c\in \mathcal{L}$ is maximal by \ref{item:laminar_avocado_3}, we know that \begin{equation}
	\text{$\Phi(f)$ is a descendant of $e_c$.}\label{eq:Phi_f_descendant_e_c}	
	\end{equation}
	By \eqref{eq:odd_number_of_vertices_from_both}, $|(e_c\setminus \Phi(f))\cap V|$ is even.
	 Moreover, $(e_c\setminus \Phi(f))\cap V\supseteq \dis_{\mathcal{L}}(e_c)\cap V=f_c\ne \emptyset$ by \ref{item:laminar_avocado_4}.
	Thus, the vertices in $(e_c\setminus \Phi(f))\cap V$ subdivide the odd cycle $C$ into a positive even number of paths.
	Not all of them can be odd.
	 Hence, there exists a subpath $Q$ of $C$ of positive even length such that both endpoints of $Q$ are contained in $e_c\setminus \Phi(f)$, but no inner vertex of $Q$ is contained in $e_c\setminus \Phi(f)$.
	 We distinguish two cases:
	
	\textbf{Case 1:} $Q$ contains the edge $f\setminus s\in E(C)$.
	 The path $Q$ splits into two paths when deleting the edge $f\setminus s$ and as $Q$ is even, exactly one of these two paths is even, call it $Q'$.
	  Then $Q'$ has one endpoint $x\in e_c\setminus \Phi(f)$ and its other endpoint $y\in f\setminus s=(\Phi(f)\setminus \Phi(s))\cap V(A)$ by \cref{prop:f_s_no_conflict}.
	  In particular, $x\ne y$.
	   By our choice of $Q$, no inner vertex of $Q'$ is contained in $e_c\setminus \Phi(f)$. 
	   Moreover, by definition of $Q'$, no  inner vertex of $Q'$ is contained in $f\setminus s=(\Phi(f)\setminus \Phi(s))\cap V(A)$.
	    Hence, all inner vertices of $Q'$ that are contained in $e_c$ are contained in $\Phi(s)$, so
	    \begin{equation}
	    \{x,y\}=V(Q')\cap (e_c\setminus \Phi(s)) \label{eq:x_y_difference}
	    \end{equation}
	Let $C'\coloneqq Q'+\{x,y\}$ and define $\Phi'\colon E(C')\rightarrow E(G)$ via $\Phi'\upharpoonright E(Q')=\Phi\upharpoonright E(Q')$ and $\Phi'(\{x,y\})=e_c$. Note that $E(Q')\subseteq E$.
	Our goal is to show that $(C',\Phi')$ satisfies the requirements of \cref{lemma:odd_cycle}.
	As $Q'$ is an even path with distinct endpoints, $C'$ is an odd cycle, i.e., \ref{item:laminar_odd_cycle_1} holds.

	 As $f_c\subseteq \dis_{\mathcal{L}}(e_c)\subseteq e_c\setminus \Phi(f)$ by \ref{item:laminar_avocado_4} and \eqref{eq:Phi_f_descendant_e_c}, we have $f_c\notin E(Q')$ because no edge of the even path $Q'$ has both endpoints in $e_c\setminus \Phi(f)$. As $\Phi$ is injective, $e_c\notin \Phi(E(Q'))$. Moreover, we get that $\Phi'$ is injective, i.e., \ref{item:laminar_odd_cycle_2} holds. As $\{x,y\}\subseteq e_c$, we have $(C',\Phi')\qin G$.

	 By \cref{cor:restriction} and \ref{item:laminar_avocado_3}, we have $(Q',\Phi'\upharpoonright Q')=(Q',\Phi\upharpoonright Q')\pin G$.
	  Hence, $e_c$ is the only potential conflict of $(C',\Phi')$. 
	  If $e_c$ is not a conflict, then $C'\pin G$ is an odd cycle.
	Otherwise, $e_c$ is the unique conflict and $e_c\in\mathcal{L}$ by \ref{item:laminar_avocado_3}.
	So \ref{item:laminar_odd_cycle_3} holds.

	 For \ref{item:laminar_odd_cycle_4}, we set $\mathcal{D}\coloneqq \{\Phi(s)\}$. Note that $\Phi(s)$ is a descendant of $e_c$ by \cref{prop:f_s_no_conflict} and \eqref{eq:Phi_f_descendant_e_c}. By \eqref{eq:x_y_difference}, we have $\{x,y\}=(e_c\setminus \Phi(s))\cap V(Q')=(e_c\setminus \Phi(s))\cap V(C')$. By \cref{prop:no_ancestor}, there is no $e\in E(Q')\subseteq E\setminus \{f_c\}$ such that $\Phi(e)$ is an ancestor of $\Phi(s)$. This concludes the proof of \ref{item:laminar_odd_cycle_4}.
	By \cref{lemma:odd_cycle}, $G$ contains an odd cycle, a tree house, or an avocado.

	\textbf{Case 2:} $Q$ does not contains the edge $f\setminus s\in E(C)$. Then $E(Q)\subseteq E$.
	 Moreover, as $f_c\subseteq \dis_{\mathcal{L}}(e_c)\subseteq e_c\setminus \Phi(f)$ by \ref{item:laminar_avocado_4} and \eqref{eq:Phi_f_descendant_e_c} and only the endpoints of the even path $Q$ are contained in $e_c\setminus \Phi(f)$, we have $f_c\notin E(Q)$.
	 By injectivity of $\Phi$, $e_c\notin \Phi(E(Q))$.
	  By \ref{item:laminar_avocado_3} and \cref{cor:restriction}, $(Q,\Phi\upharpoonright Q)\pin G$.
	   Denote the endpoints of $Q$ by $x$ and $y$.
	   Let $C'\coloneqq Q+\{x,y\}$ and define $\Phi'\colon E(C')\rightarrow E(G)$ via $\Phi'\upharpoonright E(Q)=\Phi\upharpoonright E(Q)$ and $\Phi'(\{x,y\})=e_c$.
	   Again, our goal is to show that $(C',\Phi')$ meets the requirements of \cref{lemma:odd_cycle}.
	    By definition, $C'$ is an odd cycle, so \ref{item:laminar_odd_cycle_1} holds.

	    As $\Phi$ is injective and $e_c\notin \Phi(E(Q))$, $\Phi'$ is injective, so \ref{item:laminar_odd_cycle_2} holds. As $\{x,y\}\subseteq e_c$, we have $(C',\Phi')\qin G$.

	     As $(Q,\Phi'\upharpoonright Q)=(Q',\Phi\upharpoonright Q)\pin G$ is conflict-free, $e_c$ is the only potential conflict for $(C',\Phi')$.
	      If $e_c$ is not a conflict, then $C'\pin G$ is an odd cycle.
	Otherwise, $e_c$ is the unique conflict and we have $e_c\in \mathcal{L}$ by \ref{item:laminar_avocado_3}. So \ref{item:laminar_odd_cycle_3} holds.

	For \ref{item:laminar_odd_cycle_4}, we let $\mathcal{D}\coloneqq \{\Phi(f)\}$.
	 By \eqref{eq:Phi_f_descendant_e_c}, $\Phi(f)$ is a descendant of $e_c$.
	 By construction of $Q$, $\{x,y\}=(e_c\setminus \Phi(f))\cap V(Q)=(e_c\setminus \Phi(f))\cap V(C')$.
	By \cref{prop:no_ancestor}, there is no $e\in E(Q)\subseteq E\setminus \{f_c\}$ such that $\Phi(e)$ is an ancestor of $\Phi(f)$. 
	This concludes the proof of \ref{item:laminar_odd_cycle_4}.
	By \cref{lemma:odd_cycle}, $G$ contains an odd cycle, a tree house, or an avocado.
\end{proof}
%%%%%%%%%%%%%%%%%%%%%%%%%%%%%%
%%%%%%%%%%%%%%%%%%%%%%%%%%%%%%
\section{Generalization to mixed hypergraphs}\label{secSigned}
%%%%%%%%%%%%%%%%%%%%%%%%%%%%%%
%%%%%%%%%%%%%%%%%%%%%%%%%%%%%%

In this section, we generalize \Cref{theorem:size_four} to disjoint mixed hypergraphs, which generalize mixed graphs from~\cite{S2003V1}.
A mixed graph has undirected edges, each consisting of two tail vertices, and directed edges, each consisting of one head and one tail vertex.
Our generalization to mixed hypergraphs allows a hyperarc to have a set of tails and a set of heads.
We remark that a more general notion has been studied in ~\cite{Rusnak_2013} under the name \emph{oriented hypergraphs}.

\begin{definition}
	A {\bf mixed hypergraph} is a pair $D=(V,A)$, where $V \eqqcolon V(D)$ is a finite vertex set and  $A \eqqcolon A(D) \subseteq (2^V \times 2^V) \setminus \{(\emptyset, \emptyset)\}$ is a finite hyperarc multiset.
	Each $a\in A$ has the form $a=(S(a),T(a))$ for disjoint sets $S(a), T(a) \subseteq V(H)$. 
	The {\bf underlying hypergraph} of $D$ is $(V,E)$, where $E \coloneqq \{S(a) \cup T(a) : a \in A \}$.
\end{definition}

The definition of a partial subhypergraph (\cref{def:partial}) generalizes naturally to the mixed setting.
We write $H \pin D$ if there exists a function $\iota$ such that $(H,\iota)$ is a partial subhypergraph of $D$.

For $n\in\mathbb{Z}_{\ge 1}$, we denote the set of hyperarcs of size at least $n$ by $A_{\geq n}(D)\coloneqq \{a\in A(D): |S(a)\cup T(a)|\geq n\}$. 
We say that $D$ is {\it disjoint} if the sets $(S(a)\cup T(a))_{a\in A_{\geq 4}(D)}$ are pairwise disjoint.
The {\it incidence matrix} $\mbfs{M}(D)\in \{0,\pm 1\}^{V(D)\times A(D)}$ of $D$ is defined component-wise by
\[	
\mbfs{M}(D)_{v,a} \coloneqq
\left\{
\begin{array}{rl}
	\phantom{-}1 &\text{if}~ v\in S(a)\\
	-1 &\text{if}~ v\in T(a)\\
	\phantom{-}0 &\text{if}~ v\not\in S(a)\cup T(a)
\end{array}\right..
\]
We say that $D$ is {\it unimodular} if $\mbfs{M}(D)$ is TU.

In order to generalize \Cref{theorem:size_four}, we must generalize the notions of odd cycle and tree house to the mixed setting.
This relies on the observation that each hyperarc $a$ with $|S(a)| = |T(a)| = 1$ can be split into two hyperarcs $a', a''$ having $T(a')=T(a'') = \emptyset$; see \Cref{fig:DtoU}.
This splitting operation has also been considered in~\cite{Rusnak_2013}.

\begin{definition}\label{def:Split}
	Let $D = (V,A)$ be a mixed hypergraph, and let $a \in A$ with $S(a) = \{u'\}$ and $T(a) = \{u''\}$. 
	{\bf Splitting $a$} is the operation that creates a mixed hypergraph with vertex set $V \cupdot \{w\}$ and hyperarc set $A \setminus \{a\} \cupdot \{a', a''\}$, where $S(a') = \{u',w\}$, $S(a'') = \{u'',w\}$, and $T(a') = T(a'') = \emptyset$.
\end{definition}

\begin{figure}[ht]
	\footnotesize
	\centering
	\begin{tabular}{c@{\hskip .2 cm}c@{\hskip .75 cm}c@{\hskip .2 cm}c}
		\begin{tikzpicture}[scale = 1, mynode/.style={circle, draw, fill, minimum size = .15 cm, inner sep=0pt}, baseline = 15]
			\node[mynode, label={[below = 3 pt] $v_0$}] (v0) at (0,0){};
			\node[mynode, label={[below = 3 pt] $v_1$}] (v1) at (1,0){};
			\node[mynode, label={[above = 3 pt] $v_2$}] (v2) at (1, 1){};
			\node[mynode, label={[above = 3 pt] $v_3$}] (v3) at (0, 1){};
			\draw[arrows = {-Stealth[scale=1.5]}] (v0)--(v1);
			\draw[arrows = {-Stealth[scale=1.5]}] (v1)--(v2);
			\draw[arrows = {-Stealth[scale=1.5]}] (v3)--(v2);
			\draw[] (v0)--(v3);
		\end{tikzpicture}
		&
		$
		\begin{bNiceArray}{rrrr}[first-col]
			v_0& 1 & 0 & 0 & 1\\
			v_1& -1 & 1 & 0 & 0\\
			v_2& 0 & -1 & -1 & 0\\
			v_3 & 0 & 0 & 1 & 1
		\end{bNiceArray}
		$
		&
		\begin{tikzpicture}[scale = 1, mynode/.style={circle, draw, fill, minimum size = .15 cm, inner sep=0pt}, baseline = 15]
			\node[mynode, label={[below = 3 pt] $v_0$}] (v0) at (0,0){};
			\node[mynode, label={[below = 3 pt] $v_1$}] (v1) at (1,0){};
			\node[mynode, label={[above = 3 pt] $v_2$}] (v2) at (1, 1){};
			\node[mynode, label={[above = 3 pt] $v_3$}] (v3) at (0, 1){};
			\node[mynode, label={[below = 3 pt] $v_{0,1}$}] (v01) at (0.5, 0){};
			\node[mynode, label={[left = 0 pt] $v_{1,2}$}] (v12) at (1, 0.5){};
			\node[mynode, label={[above = 2 pt] $v_{2,3}$}] (v23) at (0.5, 1){};
			\draw[] (v0)--(v1) -- (v2) -- (v3) -- (v0);
		\end{tikzpicture}
		&
		$
		\begin{bNiceArray}{rrrrrrr}[first-col]
			v_0&1 & 0 & 0 & 0 & 0 & 0 & 1\\
			v_{0,1}&1 & 1 & 0 & 0 & 0 & 0 & 0\\
			v_1&0 & 1 & 1 & 0 & 0 & 0 & 0\\
			v_{1,2}&0 & 0 & 1 & 1 & 0 & 0 & 0\\
			v_2&0 & 0 & 0 & 1 & 1 & 0& 0 \\
			v_{2,3}&0 & 0 & 0 & 0 & 1& 1 & 0\\
			v_3&0 & 0 & 0 & 0 & 0& 1 & 1\\
		\end{bNiceArray}
		$
	\end{tabular}
	\caption{A mixed hypergraph is on the left. 
		The mixed hypergraph obtained by splitting the directed arcs is on the right.}
	\label{fig:DtoU}
\end{figure}

The splitting operation treats a hyperarc $a$ with $|S(a)| = |T(a)| = 1$ as an undirected path of length $2$.
In contrast, a hyperarc $a$ with $|S(a)|=2$ and $T(a)=\emptyset$, or $S(a)=\emptyset$ and $|T(a)|=2$, behaves like one edge.
We capture this idea with the notion of parity:

\begin{definition}\label{def:parity}
	Let $a$ be a hyperarc with $|S(a)\cup T(a)|=2$. 
	The {\bf parity} of $a$ is 
	\[
	\mathrm{par}(a) \coloneqq
	\left\{
	\begin{array}{ll}
		0 &\text{if}~ |S(a)|=|T(a)|=1\\
		1 &\text{if}~ |S(a)|=2 \text{ or } |T(a)|=2
	\end{array}\right..
	\]
\end{definition}

\Cref{def:parity} motivates our generalization of odd cycle.
We say $D$ is a {\it mixed path} (respectively, a {\it mixed cycle}) if its underlying hypergraph is a path (respectively, a cycle).
If splitting all hyperarcs of a mixed cycle $D$ yields an odd cycle, then $D$ should be a mixed odd cycle. 
\begin{definition}\label{def:directedOddCycle}
	For a mixed path or cycle $P$, the {\bf parity} of $P$ is
	\(
	\sum_{a\in A(P)}\mathrm{par}(a) \mod 2.
	\)
	We call $P$ {\bf odd} (respectively, {\bf even}) if its parity is odd (respectively, even).
\end{definition}
We remark that a hypergraph $D$ is a mixed odd cycle if and only if its incidence matrix is an unbalanced hole.

In a tree house, each of the three odd paths creates an even cycle with the proper hyperedge.
We use this to extend the definition to the mixed setting:
\begin{definition} \label{def:directed_odd_tree_house}
	A {\bf mixed tree house} is a mixed hypergraph $D$ consisting of a hyperarc $h$ with $S(h)\cup T(h)=\{r,\ell_1,\ell_2,\ell_3\}$ of size $4$ and, for each $i\in [3]$, a mixed $r$-$\ell_i$-path $P_i$ such that the vertex sets $V(P_i)\setminus\{r\},i\in [3]$ are pairwise disjoint.
	For each $i \in [3]$, the mixed cycle $D[V(P_i), A(P_i) \cup \{h\}] \pin D$ is even.
\end{definition}

See \Cref{figure:signed_odd_tree_house} for an example of a mixed tree house.

\begin{figure}
	\footnotesize
	\centering
	\begin{tabular}{c@{\hskip 1.5 cm}c}
		\begin{tikzpicture}[scale = 1, mynode/.style={circle, draw, fill, minimum size = .15 cm, inner sep=0pt}, baseline = 15]
			\node[mynode, draw = red, fill = red, label={[left = 3 pt] \color{red}$r$}, label={[below = 3 pt] \color{red}$+$}] (A) at (0,0){};
			\node[mynode, draw = red, fill = red, label={[right= 3 pt] \color{red}$\ell_1$}, label={[below = 3 pt] \color{red}$-$}] (B) at (1,0){};
			\node[mynode, draw = red, fill = red, label={[right= 3pt] \color{red}$\ell_2$}, label={[above] \color{red}$+$}] (C) at (1,1){};
			\node[mynode, draw = red, fill = red, label={[left= 3pt] \color{red}$\ell_3$}, label={[above] \color{red}$+$}] (D) at (0,1){};
			\node[mynode, label={[right= 3pt] $v$}] (E) at (0.5,0.5){};
			\draw[arrows = {-Stealth[scale=1.5]}] (A)--(B);
			\draw[arrows = {-Stealth[scale=1.5]}] (E)--(A);
			\draw[] (E)--(C);
			\draw[arrows = {-Stealth[scale=1.5]}] (D)--(A);
			\draw[arrows = {-Stealth[scale=1.5]}] (A)--(D);
		\end{tikzpicture}
		&
		$
		\begin{bNiceArray}{rrrrr}[first-row, first-col]
			& a_1 & a_2 & a_3 & a_4 & h\\
			r & 1 & -1 & 0 & -1 & 1\\
			\ell_1 & -1 & 0 & 0 & 0 & -1\\
			v&  0 & 1 & 1 & 0 & 0 \\
			\ell_2 & 0 & 0 & 1 & 0 & 1\\
			\ell_3 & 0 & 0 & 0 & -1 & 1 
		\end{bNiceArray}
		$
	\end{tabular}
	\caption{A mixed tree house $D' $ and its incidence matrix $\mbfs{M}(D')$. 
		For the proper hyperarc $h$, the signs are indicated in red above and below the vertices.
	}\label{figure:signed_odd_tree_house}
\end{figure}

\Cref{thm:camion} is a special case of Camion's characterization of TU matrices.
We state the full version now.
A mixed hypergraph $D$ is {\it Eulerian} if every row and column of $\mbfs{M}(D)$ contains an even number of non-zero entries.

\begin{theorem}[{Camion \cite[Theorem 2]{Camion}}] 
	A mixed hypergraph $D$ is unimodular if and only if $\sum_{a \in A(H)} \sum_{v \in V(H)} \mbfs{M}(H)_{v,a}  \equiv 0 ~\mathrm{mod}~ 4$ for each Eulerian $H \pin D$ with $|V(H)|=|A(H)|$. \label{cor:camion_directed}	
\end{theorem}

We now prove~\Cref{thm:DirectedUniEquivalence}.

\begin{proof}[Proof of~\Cref{thm:DirectedUniEquivalence}]
	By definition, a mixed odd cycle $C$ has an odd number of odd-parity arcs. 
	If an arc $a \in A(C)$ has odd parity, then $\sum_{v \in V(C)} \mbfs{M}(C)_{v,a}\equiv 2 \mod 4$.
	Otherwise, $\sum_{v \in V(C)} \mbfs{M}(C)_{v,a} = 0$.
	Hence, $C$ is non-unimodular by \Cref{cor:camion_directed}.
	Consider a mixed tree house $W$ defined by paths $P_1$, $P_2$, and $P_3$, and proper hyperarc $h$ with $S(h)\cup T(h) = \{r, \ell_1, \ell_2, \ell_3\}$.
	The cycle $W[V(P_i), A(P_i) \cup \{h\}]$ is even for each $i \in [3]$.
	Thus, 
	\begin{align*}
		\alpha_i \coloneqq \mbfs{M}(W)_{r,h} + \mbfs{M}(W)_{\ell_i, h} + \sum_{a \in A(P_i)} \sum_{v \in V(W)} \mbfs{M}(W)_{v,a} & \equiv 0 \mod 4, && \text{and}\\
		\sum_{a \in A(W)} \sum_{v \in V(W)} \mbfs{M}(W)_{v,a} = \sum_{i=1}^3 \alpha_i - 2 \mbfs{M}(W)_{r,h} & \equiv 2 \mod 4.
	\end{align*}
	Hence, $W$ is non-unimodular by \Cref{cor:camion_directed}.
	Thus, if $D$ is unimodular, then it does not contain a mixed odd cycle or mixed tree house.

	Assume that $D$ is non-unimodular.
	As $\pin$ is transitive, we may assume that $D$ is disjoint, Eulerian, $|V(D)| = |A(D)|$, and $|\det \mbfs{M}(D)| \ge 2$.
	We reduce $D$ to a non-unimodular hypergraph $G$ with $T(a) = \emptyset$ for all $a \in A(G)$, while preserving the existence of mixed odd cycles and mixed tree houses (or the lack thereof), and the absolute value of the determinant.

	First, multiply rows of $\mbfs{M}(D)$ by $-1$ so that $T(a) = \emptyset$ for all $a \in A_{\ge 4}(D)$, which is possible because $D$ is disjoint.
	For a mixed cycle $C \pin D$ and a vertex $v \in V(C)$, multiplying row $\mbfs{M}(D)_{v, \cdot}$ by $-1$ swaps the parities of the two arcs in $A(C)$ incident to $v$, thus preserving the parity of $C$.
	If $D$ contained a mixed odd cycle (respectively, a mixed tree house) before this reduction, then it also contains a mixed odd cycle (respectively, a mixed tree house) after this reduction, and vice versa.

	Second, multiply column $\mbfs{M}(D)_{\cdot, a}$ by $-1$ for each $a \in A(D)$ such that $S(a) = \emptyset$ and $|T(a)| = 2$.
	This operation preserves the parity of $a$, thus preserving the parity of any cycle that contains $a$.
	After this reduction, each hyperarc $a \in A(D)$ satisfies either $T(a) = \emptyset$ or $|S(a)| = |T(a)| = 1$.
	
	Third, split any hyperarc $a \in A(D)$ satisfying $|S(a)| = |T(a)| = 1$.
	To help discuss the effects of the splitting operation, let $u', u'', a', a''$, and $w$ be as in \Cref{def:Split}, and let $D'$ be the mixed hypergraph obtained after splitting $a$.
	The splitting operation induces a 1-1 correspondence between cycles $C \pin D$ and $C' \pin D'$ that preserves cycle parity.
	To see this correspondence, observe that because $w$ has degree $2$, $C'$ must contain either both $a'$ and $a''$ or neither of them.
	If $C'$ contains both $a'$ and $a''$, then $C$ can be obtained by replacing $a'$ and $a''$ by $a$.
	Otherwise, $C$ corresponds to $C'$.
	By definition, the parity of $a$ equals the sum of the parities of $a'$ and $a''$.
	Thus, $C$ and $C'$ have the same parity.
	
	The first two reductions preserve $|\det\mbfs{M}(D)|$ as they multiply a row or column by $-1$.
	Similarly, the splitting operation preserves $|\det\mbfs{M}(D)|$.
	To see why, let $\mbfs{B}$ be the matrix obtained from $\mbfs{M}(D')$ by replacing column $\mbfs{M}(D')_{\cdot, a'}$ with $\mbfs{M}(D')_{\cdot, a'} - \mbfs{M}(D')_{\cdot, a''}$.
	The matrices $\mbfs{M}(D')$ and $\mbfs{B}$ differ by an elementary column operation, so $|\det \mbfs{M}(D')| = |\det \mbfs{B}|$.
	The only non-zero entry in row $\mbfs{B}_{w, \cdot}$ is in column $a''$, and deleting row $\mbfs{B}_{w, \cdot}$ and column $\mbfs{B}_{\cdot, a''}$ gives $\mbfs{M}(D)$.
	Using Laplace expansion along $\mbfs{B}_{w, \cdot}$ gives $|\det\mbfs{M}(D)| =  |\det \mbfs{B}| = |\det \mbfs{M}(D')| $.

	After the third reduction, $T(a) = \emptyset$ for all $a \in A(D)$.
	Hence, it suffices to show that every non-unimodular disjoint mixed hypergraph in which every hyperarc $a$ has $T(a) = \emptyset$, contains a mixed odd cycle or a mixed tree house.
	The incidence matrices of such mixed hypergraphs are incidence matrices of non-mixed hypergraphs.
	Thus, the result follows directly from~\Cref{theorem:size_four}.
\end{proof}

As incidence matrices of odd cycles have a determinant of $2$, the previous proof reveals the following, which has already been shown in~\cite{ZASLAVSKY198247}.

\begin{observation}[\cite{ZASLAVSKY198247}]\label{lem:MOC_MOTH_det}
	If $D$ is a mixed odd cycle, then $|\det \mbfs{M}(D)| = 2$.
\end{observation}

\section[On the Cornu{\'e}jols-Zuluaga Conjecture]{On the Cornu{\'e}jols-Zuluaga Conjecture}\label{secCorn}
%%%%%%%%%%%%%%%%%%%%%%%%%%%%%%
%%%%%%%%%%%%%%%%%%%%%%%%%%%%%% 

%%%%%%%%%%%%%%%%%%%%%%%%%%%%%%
%%%%%%%%%%%%%%%%%%%%%%%%%%%%%%
\subsection[Almost TU matrices of laminar hypergraphs]{A list of almost TU incidence matrices of mixed disjoint and laminar hypergraphs}\label{ssecalmostTU}
%%%%%%%%%%%%%%%%%%%%%%%%%%%%%%
%%%%%%%%%%%%%%%%%%%%%%%%%%%%%% 

(Mixed) odd cycles, (mixed) tree houses, and avocados are the forbidden substructures that prevent a mixed disjoint or (non-mixed) laminar hypergraph from being unimodular.
In this section, we show that their incidence matrices are almost TU.
We first prove that they cannot contain each other.
\begin{proposition}\label{lem:containment_odd_cycle}
Let $C$ be a mixed odd cycle.
Then, $C$ does not contain a mixed odd cycle, mixed tree house, or avocado as a proper partial subhypergraph.
\end{proposition}
\begin{proof}
Neither $C$ nor any partial subhypergraph of $C$ contain a hyperedge of size at least $4$.
As a mixed tree house and an avocado each have a hyperedge of size at least $4$, $C$ cannot contain either.
It remains to show that $C$ does not properly contain a mixed odd cycle, which is immediate, as removing any vertex or edge from $C$ leaves a disjoint union of mixed paths.
\end{proof}

\begin{proposition}\label{lem:containment_tree_house}
Let $T$ be a mixed tree house with mixed paths $(P_i)_{i\in [3]}$ and hyperarc $h$ with $S(h)\cup T(h)=\{r,\ell_1,\ell_2,\ell_3\}$.
Then, $T$ does not contain an mixed odd cycle, mixed tree house, or avocado as a proper partial subhypergraph.
\end{proposition}
\begin{proof}
Note that the only Eulerian proper partial subhypergraphs of $T$ are the mixed even cycles of the form $C_i \coloneqq T[V(P_i), E(P_i) \cup \{h\}]$.
As mixed odd cycles, mixed tree houses, and avocados are Eulerian, $T$ cannot contain any of them.
\end{proof}

\begin{lemma}\label{lem:containment_avocado}
Let $A = (V \cup \{t\}, E \cup \{s,f\})$ be an avocado.
Then, $A$ does not contain a mixed odd cycle, mixed tree house, or avocado as a proper partial subhypergraph.
\end{lemma}
\begin{proof}
As $\mbfs{M}(A)$ is a $\{0,1\}$-matrix, if $A$ contains a mixed odd cycle or mixed tree house, then the respective hypergraph also forms a non-mixed odd cycle/ tree house.

Suppose, towards a contradiction, that $H = A[U,F]$ is a proper partial subhypergraph of $A$ that is an odd cycle, a tree house, or an avocado.
Then every vertex of $U$ has even degree at least $2$ in $H$, and every hyperedge of $F[U]$ has even size at least $2$.
We may further assume that 
\begin{equation}
e\cap U \ne \emptyset \text{ for every $e\in F$}\label{eq:non_empty_intersection}
\end{equation}
because we can replace $F$ by $F'\coloneqq \{e\in F\colon e\cap U\ne \emptyset\}$
otherwise.

First, we assume that $F$ contains at most one of $s$ or $f$.
In particular, $H$ contains at most one proper hyperedge, so it is either an odd cycle or a tree house. 
As the maximum degree of a vertex in $A[V\cup \{t\}, E\cup\{s\}]$, as well as $A[V\cup \{t\}, E\cup\{f\}]$, is at most $3$ and $H$ is a partial subhypergraph of one of the two, $H$ cannot be a tree house (since a tree house has a vertex of degree $4$). 
Hence, $H$ is an odd cycle.
As $A[V\cup \{t\}, E]$ is acyclic (it is a path together with an isolated vertex), we must have $\{s,f\} \cap F \neq \emptyset$.
Thus, $H$ consists of a subpath $Q$ of $(V,E)$ closed to a cycle by the hyperedge $e\in \{s,f\}$. 
Moreover, $V(Q)$ intersects $e$ precisely in its endpoints. 
However, by \ref{item:avocado_3}, consecutive vertices from $s$ and $f$, respectively, have an odd distance on the path $(V,E)$.
Hence, $H$ cannot be an odd cycle.

Second, assume that $F$ contains both $s$ and $f$.
As odd cycles, tree houses, and avocados do not contain parallel hyperedges, we have $s[U] \ne f[U]$.
Thus, $U\cap (f\setminus s)$ is non-empty. 
Let $f\setminus s=\{v,w\}$ with $v\in U\cap (f\setminus s)$, and $P=(V,E)$ be the even $v$-$w$-path in $A$; see \ref{item:avocado_2}.
	\begin{claim}
		We have $V \subseteq U$ and $E \subseteq F$.
	\end{claim}
	\begin{cpf}
		It suffices to show that $E\subseteq F$.
		Indeed, as every hyperedge of $H$ has size at least $2$, \eqref{eq:non_empty_intersection} implies that both endpoints of every $e\in E\cap F$ must be contained in $U$.
		
		Assume, towards a contradiction, that $E\not\subseteq F$ and let $e=\{x,y\}\in E$ be the first edge of $P$ not contained in $F$ when traversing it from $v$ to $w$ and suppose that $x$ is closer to $v$ on $P$ than $y$. 
		Then $x\in U$ because either $x=v$ or the predecessor edge of $e$ is contained in $F$ and contains $x$. 
		If $x=v$, then it has degree $1$ because $f\in F$ but the first edge of $P$ is not.
		If $x$ is an inner vertex of $P$, then it is contained in the predecessor edge of $e$.
		As $x \notin \{v,w\} = f \setminus s$, either both or none of the hyperedges $s$ and $f$ contain $x$.
		Hence, $x$ has odd degree in $H$, contradicting that every vertex in $H$ has even degree.
	\end{cpf} 
	By~\ref{item:avocado_2} and~\ref{item:avocado_3}, we have that $|f\cap V|$ is odd. 
	As $H$ is Eulerian, $|f\cap U|$ is even, so we must have $t\in U$. 
	Hence, $V(A)=V\cup \{t\}\subseteq U$ and $E(A)=E\cup \{s,f\}\subseteq F$, contradicting that $H$ is a proper partial subhypergraph of $A$.
\end{proof}

Since mixed odd cycles, mixed tree houses, and avocados do not properly contain one another, \cref{thm:DirectedUniEquivalence,theorem:non_unimodular_laminar} directly implies that they are exactly the mixed disjoint or laminar hypergraphs with almost TU incidence matrices.

\begin{proposition}\label{prop:almost_TU_matrices}
Let $G$ be a mixed disjoint or laminar hypergraph. 
Then $\mbfs{M}(G)$ is almost TU if and only if $G$ is a mixed odd cycle, a mixed tree house, or an avocado. \label{prop:odd_cycle_tree_house_almost_TU}
\end{proposition}
\begin{proof}
First, assume that $G$ is an mixed odd cycle, a mixed tree house, or an avocado.
Then $\mbfs{M}(G)$ is non-TU by \cref{thm:DirectedUniEquivalence,theorem:non_unimodular_laminar}.
By \cref{lem:containment_odd_cycle,lem:containment_tree_house,lem:containment_avocado}, $G$ does not contain a mixed odd cycle, a mixed tree house, or an avocado as a proper partial subhypergraph.
Hence, again by \cref{thm:DirectedUniEquivalence,theorem:non_unimodular_laminar}, every partial subhypergraph of $G$ is unimodular.
Thus, $\mbfs{M}(G)$ is almost TU.

Next, assume $\mbfs{M}(G)$ is almost TU.
As $G$ is not unimodular, \cref{thm:DirectedUniEquivalence,theorem:non_unimodular_laminar} state there exists some $H\pin G$ that is a mixed odd cycle, a mixed tree house, or an avocado.
Every proper submatrix of $\mbfs{M}(G)$ is TU, so $H=G$.
\end{proof}

%%%%%%%%%%%%%%%%%%%%%%%%%%%%%%
%%%%%%%%%%%%%%%%%%%%%%%%%%%%%%
\subsection{The conjecture holds for incidence matrices of mixed disjoint and laminar hypergraphs}\label{ssecalmostTUDisjoint}
%%%%%%%%%%%%%%%%%%%%%%%%%%%%%%
%%%%%%%%%%%%%%%%%%%%%%%%%%%%%% 

By \cref{prop:almost_TU_matrices}, the only mixed disjoint or laminar hypergraphs with an almost TU incidence matrix
are mixed odd cycles, mixed tree houses, and avocados.
In this section, we establish \cref{cor:CZConj} for the incidence matrix $\mbfs{A}$ of each such (mixed)
hypergraph, together with its transpose $\mbfs{A}^\top$.

For each of the three families, we exhibit a TU matrix $\mbfs{L}$ such that $\mbfs{L}\mbfs{A}$ is an
unbalanced hole.
This directly establishes \cref{cor:CZConj} for $\mbfs{A}$.
To handle $\mbfs{A}^\top$, we exhibit a TU matrix $\mbfs{R}$ such that $\mbfs{A}\mbfs{R}$ is an
unbalanced hole.
Using that $(\mbfs{A}\mbfs{R})^\top = \mbfs{R}^\top\mbfs{A}^\top$, that $\mbfs{R}^\top$
is TU whenever $\mbfs{R}$ is, and that the transpose of an unbalanced hole is again an unbalanced
hole, the existence of $\mbfs{R}^\top$ establishes \cref{cor:CZConj} for $\mbfs{A}^\top$.

The three families call for different techniques.
For mixed odd cycles, $\mbfs{A}$ is itself an unbalanced hole, so $\mbfs{L} = \mbfs{R} = \mbfs{I}$ suffices for both directions at once.
For mixed tree houses, we construct $\mbfs{R}$ directly and obtain $\mbfs{L}$ from it by transposing, using the fact that the transpose of a mixed tree house's incidence matrix is again the incidence matrix of a mixed tree house.
For avocados, no such shortcut is available.
Indeed, it is easy to see that the transpose of an avocado's incidence matrix is not even the incidence matrix of a laminar hypergraph.
Thus, we construct $\mbfs{L}$ and $\mbfs{R}$ separately.

\smallskip

\paragraph{Mixed odd cycles}
Incidence matrices of odd cycles are known to be unbalanced holes.
\begin{observation}\label{lem:conj_oddCycle}
    Let $C$ be a mixed odd cycle.
    Then $\mbfs{M}(C)$ is an unbalanced hole.
\end{observation}
\begin{proof}
    $\mbfs{M}(C)$ is a hole and the entries of $\mbfs{M}(C)$ sum to $2 \bmod 4$ by definition of the parity of a mixed cycle.
\end{proof}

\smallskip

\paragraph{Mixed tree houses}
The following property of mixed even cycles is used when proving~\cref{conjCorn} for mixed tree houses.
It generalizes the idea that the edges of a (non-mixed) even cycle can be signed so that no two adjacent edges have the same sign and follows from \cite{ZASLAVSKY198247} (Theorem 5.1 (g)).

\begin{observation}[\cite{ZASLAVSKY198247}]\label{lem:MEC}
	If $C$ is a mixed even cycle, then there exists a vector $\mbfs{u} \in \{\pm 1\}^{A(C)}$ such that $\mbfs{M}(C) \mbfs{u} =\mbfs{0}$.
	Consequently, $\det \mbfs{M}(C) = 0$.
\end{observation}
\begin{lemma}\label{lem:conj_treeHouse}
	Let $D$ be a mixed tree house.
	Then there exists a TU matrix $\mbfs{R}$ such that $\mbfs{M}(D)\cdot\mbfs{R}$ is an unbalanced hole.
\end{lemma}
\begin{proof}
	Let $h\in A(D)$ with $S(h) \cup T(h) = \{r, \ell_1, \ell_2, \ell_3\}$ be the hyperarc of size $4$ defining $D$.
	For $i \in [3]$, let $P_i$ be the mixed $r$-$\ell_i$-path defining $D$.
	By~\Cref{def:directed_odd_tree_house}, $C := D[V(P_1), A(P_1) \cup \{h\}]\pin D$ is a mixed even cycle.
	Thus, by~\Cref{lem:MEC}, there exists $\mbfs{u} \in \{\pm 1\}^{A(C)}$ such that $(\mbfs{M}(P_1)\mbfs{u})_v = -\mbfs{M}(D)_{v,h}$ for $v \in \{r, \ell_1\}$, and $(\mbfs{M}(P_1)\mbfs{u} )_v = 0$ otherwise.

	Define $\mbfs{x} \in \{0,\pm 1\}^{A(D)}$ to be $\mbfs{x}_f \coloneq \mbfs{u}_f$ if $f \in A(P_1)$,  $\mbfs{x}_h \coloneq 1$, and $\mbfs{x}_f \coloneq 0$ otherwise.
	The vector $\mbfs{x}$ encodes column operations that zeros the entries $\mbfs{M}(D)_{r,h}$ and $\mbfs{M}(D)_{\ell_1,h}$.
	Let $a$ be the arc on $P_2$ with $r \in S(a) \cup T(a)$, and let $\sigma \in \{\pm 1\}$ satisfy $\sigma(\mbfs{M}(P_1)\mbfs{u})_r =  -\mbfs{M}(D)_{r,a}$.
	Define $\mbfs{y} \in \{0,\pm 1\}^{A(D)}$ to be $\mbfs{y}_f \coloneq \sigma\mbfs{u}_f$ if $f \in A(P_1)$,  $\mbfs{y}_{a} \coloneq 1$, and $\mbfs{y}_f \coloneq 0$ otherwise.
	The vector $\mbfs{y}$ encodes column operations that zeros the entry $\mbfs{M}(D)_{r,a}$ while changing $\mbfs{M}(D)_{\ell_1,a}$ from $0$ to $\pm 1$.

	Consider $\mbfs{R} \in \{0,\pm 1\}^{A(D) \times A(D)}$, where $\mbfs{R}_{h, \cdot} := \mbfs{x}$, $\mbfs{R}_{a, \cdot} := \mbfs{y}$, and $\mbfs{R}_{f, \cdot} $ is the standard unit vector for all $f \in  A(D) \setminus \{h,a\}$.
	In terms of the mixed hypergraph $D$, the vector $\mbfs{x}$ removes $r$ and $\ell_1$ from $h$, and $\mbfs{y}$ replaces the first arc (connecting $r$ and $v^2_1$) in the mixed path $P_2$ with an arc between $\ell_1$ and $ v^2_1$.
	Hence, $\mbfs{R}$ encodes column operations, so $|\det \mbfs{R}| = 1$.
	Moreover, $\mbfs{M}(D)\mbfs{R}$ is the incidence matrix of a mixed cycle.
	This mixed cycle is odd because $|\det \mbfs{M}(D)\mbfs{R}| = |\det \mbfs{M}(D)| \cdot |\det \mbfs{R}|  = |\det \mbfs{M}(D)| = 2$, where the final equality follows from~\Cref{lem:MOC_MOTH_det}.
	
	Finally, we check that $\mbfs{R}$ is TU.
	Remove all unit columns or rows from $\mbfs{R}$ as they do not affect the TU property of a matrix.
	After removing them, we are left with the $2$-column matrix $[\mbfs{u}~\sigma\mbfs{u}]$, which is TU because $\mbfs{u} \in \{\pm 1\}^{A(P_1)}$ and $\sigma \in \{\pm 1\}$.
\end{proof}

To obtain $\mbfs{L}$, it remains to show that tree houses are closed under transposing their incidence matrix.

\begin{lemma}
For a mixed tree house $D$ with incidence matrix $\mbfs{M} \coloneqq \mbfs{M}(D)$, $\mbfs{M}^\top$ is also the incidence matrix of a mixed tree house. \label{lemma:transpose_odd_tree_house}
\end{lemma}
\begin{proof}
Let $D^\top$ be the hypergraph having $\mbfs{M}^\top$ as its incidence matrix.
By \cref{prop:almost_TU_matrices}, $\mbfs{M}$ is almost TU, so $\mbfs{M}^\top$ is almost TU.
Also, $\mbfs{M}^\top$ only has one column with more than two non-zero entries.
Thus, $D^\top$ is disjoint, and it is not a mixed cycle.
By \cref{prop:almost_TU_matrices}, we can infer that $D^\top$ is a mixed tree house.
\end{proof}

\smallskip

\paragraph{Avocados}
\begin{lemma}\label{lem:conj_avocado}
	Let $A$ be an avocado.
	Then there exists a TU matrix $\mbfs{L}$ such that $\mbfs{L}\cdot \mbfs{M}(A)$ is an unbalanced hole.
\end{lemma}
\begin{proof}
	Let $A=(V\cup \{t\},E\cup \{s,f\})$ be an avocado, $C=(V,E\dot{\cup}\{f\setminus s\})$ be the odd cycle from \ref{item:avocado_2}, and $f\setminus s\eqqcolon \{v,w\}$.
	We define the matrix $\mbfs{L}\in \{0,\pm 1\}^{V(A)\times V(A)}$ as follows:
	for $x\in V(A)$, we define $\mbfs{L}_{x,x}=1$, and
	for $x\in (s\cup \{w\}) \setminus \{t\}$, we let $\mbfs{L}_{x,t}=-1$.
	We set all remaining entries to $0$.
	Then, every column of $\mbfs{L}$ other than column $t$ is a standard basis vector, so $\mbfs{L}$ is TU.
	
	Let $\mbfs{H} \coloneqq \mbfs{L}\cdot \mbfs{M}(A)$.
	For $x\in V\setminus f$, the row $\mbfs{H}_{x,\cdot}=\mbfs{M}(A)_{x,\cdot}$ encodes that $x$ is incident to its two incident edges in $E$ and nothing else.
	For $x\in s\setminus\{t\}$, $\mbfs{H}_{x,\cdot}$ arises by subtracting $\mbfs{M}(A)_{t,\cdot}$ from $\mbfs{M}(A)_{x,\cdot}$.
	This means that $x$ remains incident to its two incident edges in $E$, but is removed from the hyperedges $s$ and $f$.
	We have $\mbfs{H}_{v,\cdot}=\mbfs{M}(A)_{v,\cdot}$, so $v$ remains incident to its incident edge in $E$ and the hyperedge $f$.
	Similarly, we have $\mbfs{H}_{t,\cdot}=\mbfs{M}(A)_{t,\cdot}$, so $t$ remains incident to $f$ and $s$.
	Finally, $\mbfs{H}_{w,\cdot}$ arises by subtracting $\mbfs{M}(A)_{t,\cdot}$ from $\mbfs{M}(A)_{w,\cdot}$.
	In doing so, $w$ is removed from $f$, but added to $s$ with a negative sign.
	It remains incident to its incident edge in $E$. 

	By the above, $\mbfs{H}$ arises from the incidence matrix of the even cycle
	$\widetilde{C}\coloneqq (V\cup\{t\},\,E\cup\{\{t,v\},\{t,w\}\})$ by switching the entry in row $w$ and column $\{t,w\}$ from $1$ to $-1$.
	Hence, $\mbfs{H}$ is a hole.
	As $\widetilde{C}$ is even, the sum of entries of $\mbfs{M}(\widetilde{C})$ is divisible by $4$, and switching a $1$ to $-1$ changes this sum by $2$.
	Hence, the entries of $\mbfs{H}$ sum to $2\mod 4$.
\end{proof}

\begin{lemma}\label{lem:conj_avocado_transpose}
	Let $A$ be an avocado.
	Then there exists a TU matrix $\mbfs{R}$ such that $\mbfs{M}(A)\cdot \mbfs{R}$ is an unbalanced hole.
\end{lemma}
\begin{proof}
	Let $A=(V\cup \{t\},E\cup \{s,f\})$ be an avocado and let $C=(V,E\dot{\cup}\{f\setminus s\})$ be the odd cycle from \ref{item:avocado_2}. Let $f\setminus s\eqqcolon \{v,w\}$ and let $P=(V,E)$.
	Let $u_0=v,u_1,\dots,u_{2k-1},u_{2k}=w$ be the sequence of vertices from $f$ in the order in which they appear on $P$; note that there is an odd number of these vertices by \ref{item:avocado_3}. 
	Moreover, \ref{item:avocado_3} tells us that $P[u_i,u_{i+1}]$ is odd for $i\in\{0,\dots,2k-1\}$. We call an edge $e\in E$ \emph{odd} if there is $j\in \{0,\dots,k-2\}$ such that $e$ is the first, third, fifth, $\dots$ edge on $P[u_{2j+1},u_{2j+2}]$, and we denote the set of odd edges by $E_{odd}$. 
	Similarly, we call an edge \emph{even} if there is $j\in \{0,\dots,k-2\}$ such that $e$ is the second, fourth, sixth, $\dots$ edge on $P[u_{2j+1},u_{2j+2}]$, and we denote the set of even edges by $E_{even}$. 
	Finally, let $e^*$ be the first edge on $P[u_{2k-1},u_{2k}]$, note that $e^*\notin E_{odd}\cup E_{even}$. 
	We are now ready to define the matrix $\mbfs{R}\in \{0,\pm 1\}^{E(A)\times E(A)}$.
	For each $e\in E(A)$, we set $\mbfs{R}_{e,e}=1$. 
	Moreover, we set $\mbfs{R}_{s,f}=-1$, $\mbfs{R}_{e,s}=-1$ for $e\in E_{odd}$, $\mbfs{R}_{e,s}=1$ for $e\in E_{even}$, $\mbfs{R}_{s,e^*}=-1$, $\mbfs{R}_{e,e^*}=1$ for $e\in E_{odd}$ and $\mbfs{R}_{e,e^*}=-1$ for $e\in E_{even}$. %
	All remaining entries are set to $0$. 
	See \cref{fig:cornuejols_avocado_1}.
	
	Let $\mbfs{H}\coloneqq \mbfs{M}(A)\cdot \mbfs{R}$. 
	We have $\mbfs{H}_{\cdot,f}=\mbfs{M}(A)_{\cdot,f}-\mbfs{M}(A)_{\cdot,s}$, so the column $\mbfs{H}_{\cdot,f}$ encodes the edge $f\setminus s$. 
	For $e\in E\setminus \{e^*\}$, we have $\mbfs{H}_{\cdot,e}=\mbfs{M}(A)_{\cdot,e}$, so the column $\mbfs{H}_{\cdot,e}$ encodes the edge $e$.
	We have $\mbfs{H}_{\cdot,s}=\mbfs{M}(A)_{\cdot,s}-\sum_{e\in E_{odd}}\mbfs{M}(A)_{\cdot,e}+\sum_{e\in E_{even}}\mbfs{M}(A)_{\cdot,e}$. 
	For every inner vertex of one of the paths $P[u_{2j+1},u_{2j+2}]$, we get a $-1$ from the incident edge in $E_{odd}$ and a $+1$ from the corresponding edge in $E_{even}$, so the respective entry cancels. 
	The endpoints of the path $P[u_{2j+1},u_{2j+2}]$ are contained in $s$ and have one incident edge in $E_{odd}$ and no incident edge in $E_{even}$, so the respective entries cancel. 
	Hence, $\mbfs{H}_{\cdot,s}$ encodes the edge $\{u_{2k-1},t\}$ which arises from $s$ by deleting all of the vertices $u_{2j+1},u_{2j+2}$ with $j\in\{0,\dots,k-2\}$.
	Finally, $\mbfs{H}_{\cdot,e^*}=\mbfs{M}(A)_{\cdot, e^*}-\mbfs{M}(A)_{\cdot, s}+\sum_{e\in E_{odd}}\mbfs{M}(A)_{\cdot,e}-\sum_{e\in E_{even}}\mbfs{M}(A)_{\cdot,e}=\mbfs{M}(A)_{\cdot, e^*}-\mbfs{H}_{\cdot,s}$. 
	Recall that $e^*$ is the first edge of $P[u_{2k-1},u_{2k}]$, let $e^*\eqqcolon \{u_{2k-1},y\}$.
	Then $\mbfs{H}_{\cdot,e^*}$ encodes an edge from $y$ to $t$, with a $1$ in row $y$ and a $-1$ in row $t$.
	By the above, $\mbfs{H}$ arises from the incidence matrix of the even cycle $\widetilde{C}$ 
	obtained from $C$ by replacing the edge $e^* = \{u_{2k-1},y\}$ by the two edges
	$\{u_{2k-1},t\}$ (encoded by column $s$) and $\{t,y\}$ (encoded by column $e^*$) by switching the entry in row $t$ and column $e^*$ from $1$ to $-1$.
	As in the proof of \cref{lem:conj_avocado}, we can infer that $\mbfs{H}$ is an unbalanced hole.

	Finally, we show that $\mbfs{R}$ is totally unimodular. For $e\in E\setminus \{e^*\}$, $\mbfs{R}_{\cdot,e}$ is a unit column. 
	Moreover, $\mbfs{R}_{f,\cdot}$ and $\mbfs{R}_{e^*,\cdot}$ are unit rows. 
	As adding unit rows or columns does not affect total unimodularity, we may restrict our attention to $\mbfs{R}[E(A)\setminus \{e^*,f\},\{s,f,e^*\}]$. 
	In this submatrix, the column corresponding to $f$ is a unit column, so we may further restrict ourselves to $\mbfs{R}[E(A)\setminus \{e^*,f\},\{s,e^*\}]$. 
	In this matrix, the two columns arise from each other by multiplication with $-1$. 
	As both columns only have $\{0,\pm 1\}$ entries, the matrix is TU.
\end{proof}

\begin{figure}
	\begin{subfigure}[c]{0.4\textwidth}
		\centering
		\begin{tikzpicture}[mynode/.style={draw, fill, circle, inner sep = 0pt, minimum size = 2mm}, snode/.style={draw=none, fill=brown, fill opacity = 0.5, circle, inner sep = 0pt, minimum size = 5mm}]
			
			\foreach \i in {0,...,4}{
				\draw[thick] ({2*cos(\i*360/11-90)},{2*sin(\i*360/11-90)})--({2*cos((\i+1)*360/11-90)},{2*sin((\i+1)*360/11-90)});
			}
			\foreach \i in {6,...,11}{
				\draw[thick] ({2*cos(\i*360/11-90)},{2*sin(\i*360/11-90)})--({2*cos((\i+1)*360/11-90)},{2*sin((\i+1)*360/11-90)});
			}
			\draw[line width = 1mm, green!70!black] 	({2*cos(5*360/11-90)},{2*sin(5*360/11-90)}) to node[above, midway] {\textcolor{green!70!black}{$f\setminus s$}}({2*cos((5+1)*360/11-90)},{2*sin((5+1)*360/11-90)});
			
			\foreach \i in {7,9,0}{
				\draw[line width = 1mm, blue, opacity = 0.5] ({2*cos(\i*360/11-90)},{2*sin(\i*360/11-90)})-- ({2*cos((\i+1)*360/11-90)},{2*sin((\i+1)*360/11-90)});
				\node at ({2.2*cos((\i+0.5)*360/11-90)},{2.2*sin((\i+0.5)*360/11-90)}){\textcolor{blue}{$-$}};
			}
			
			\foreach \i in {8}{
				\draw[line width = 1mm, red, opacity = 0.5] ({2*cos(\i*360/11-90)},{2*sin(\i*360/11-90)})-- ({2*cos((\i+1)*360/11-90)},{2*sin((\i+1)*360/11-90)});
				\node at ({2.2*cos((\i+0.5)*360/11-90)},{2.2*sin((\i+0.5)*360/11-90)}){\textcolor{red}{$+$}};
			}
			
			\foreach \i in {0,1,4,7,10}
			{
				\node[snode] at ({2*cos(\i*360/11-90)},{2*sin(\i*360/11-90)}){};	
			}
			
			\foreach \i in {0,...,10}{
				\node[mynode] at ({2*cos(\i*360/11-90)},{2*sin(\i*360/11-90)}){};
			}
			
			\node at ({2*cos(5*360/11-90)},{2*sin(5*360/11-90)-0.3}){$w$};
			\node at ({2*cos(6*360/11-90)},{2*sin(6*360/11-90)-0.3}){$v$};
			
			\node at ({2.2*cos((4+0.5)*360/11-90)},{2.2*sin((4+0.5)*360/11-90)}){$e^*$};
			
		\end{tikzpicture}
	\end{subfigure}
	\begin{subfigure}[c]{0.55\textwidth}
		\centering
		$
		\begin{bNiceArray}{rrrrrr}[first-row,first-col]
			\CodeBefore
			\rowcolor{red!30!white}{4-6}
			\rowcolor{blue!30!white}{7-9}
			\Body
			
			& f & s & e^* & \multicolumn{3}{c}{E\setminus \{e^*\}}\\
			f& 1 & 0 & 0 & 0 & \dots & 0\\
			s& -1 & 1 & -1 & 0 & \dots & 0\\
			e^*& 0 & 0 & 1 & 0 & \dots & 0\\
			& 0 & 1 & -1 & & & \\
			{\color{red}E_{even}} & \vdots & \vdots & \vdots & & &\\
			& 0 & 1 & -1 & & & \\
			& 0 & -1 & 1 & & & \\
			{\color{blue}E_{odd}} & \vdots & \vdots & \vdots & \multicolumn{3}{c}{id}\\
			& 0 & -1 & 1 & & & \\
			& 0 & 0 & 0 & & & \\
			 E_{rest} & \vdots & \vdots & \vdots & & &\\
			& 0 & 0 & 0 & & & \\
		\end{bNiceArray}
		$
	\end{subfigure}
	\caption{Left: The cycle $C$ with vertices from $s$ indicated in brown and edges in $E_{odd}$ and $E_{even}$ drawn in blue and red. Right: The matrix $\mbfs{R}$, where $E_{rest}\coloneqq E\setminus (E_{odd}\cup E_{even}\cup \{e^*\})$ .}\label{fig:cornuejols_avocado_1}
\end{figure}

%%%%%%%%%%%%%%%%%%%%%%%%%%%%%%
%%%%%%%%%%%%%%%%%%%%%%%%%%%%%%
\section{Future directions}\label{secFuture}
Our paper provides a complete description of laminar hypergraphs with almost TU incidence matrices.
Laminar hypergraphs correspond to $\{0,1\}$-matrices for which the supports of columns with at least four $1$s form a laminar family.
A natural question for future research is how to generalize this result to $\{0,\pm 1\}$ matrices, corresponding to \emph{mixed} hypergraphs.
In \cref{thm:DirectedUniEquivalence}, we extend our characterization for disjoint hypergraphs to the mixed setting.
It would be interesting to see whether mixed laminar hypergraphs admit a similar extension.

Going beyond the laminar case, a natural next step might be to characterize unimodularity for $\{0,1\}$-matrices that arise from the incidence matrix of a graph by adding \emph{two} laminar families of hyperedges.
This direction is motivated by the fact that given two laminar families $\mathcal{E}_1$ and $\mathcal{E}_2$ on the vertex set $V$, the incidence matrix of the hypergraph $(V,\mathcal{E}_1\dot{\cup}\mathcal{E}_2)$ is known to be TU~\cite{S2003V3}.

Going beyond the unimodular case, in the introduction, we have highlighted the intimate connections between recent advances in $\Delta$-modular optimization and a structural understanding of (special classes of) totally $\Delta$-modular matrices.

As a first step in this direction, one may try to characterize disjoint hypergraphs with bounded subdeterminants.
This problem is open even for hypergraphs where all hyperedges have size at most 3 (these are disjoint by \cref{def:disjoint_hypergraph}).
De Loera and Onn~\cite{2004_DeLoera} show that any polytope with constraint matrix $\mbfs{A}$ can be represented using the incidence matrix $\mbfs{M}(G)$ of such a hypergraph $G$.
However, $\Delta(\mbfs{M}(G))$ can be larger than $\Delta(\mbfs{A})$, so we find that studying disjoint hypergraphs with hyperedges larger than $3$ is also an interesting question.

In order to highlight a challenge in characterizing disjoint hypergraphs with a totally $2$-modular incidence matrix, consider $G'$ in~\cref{fig:quasi_sh}.
Note that $G'$ does not contain two vertex-disjoint odd cycles or tree houses, yet $\Delta(\mbfs{M}(G')) > 2$.
This differs from the graph setting, where \cref{thmOCP} states that $\Delta(\mbfs{M}(G)) \le 2$ for a graph $G$ if and only if $G$ does not contain two vertex-disjoint odd cycles.

%%%%%%%%%%%%%%%%%%%%%%%%%%%%%%
%%%%%%%%%%%%%%%%%%%%%%%%%%%%%%
\section*{Acknowledgements}
%%%%%%%%%%%%%%%%%%%%%%%%%%%%%%
%%%%%%%%%%%%%%%%%%%%%%%%%%%%%%
%
The authors thank Matthias Walter for discussions on TU matrices, Jes\'{u}s De Loera for pointing to~\cite{2004_DeLoera}, and Christopher Hojny and Sten Wessel for discussing~\cite{2024_Hojny}.
\section*{Data availability statement}
No data are associated with this article. Data sharing is not applicable to this article.

\bibliographystyle{plainurl}
\bibliography{references}
 
%%%%%%%%%%%%%%%%%%%%%%%%%%%%%

\end{document}